\input amstex
\documentstyle{amsppt}
\magnification=\magstep1
\hsize=5in
\vsize=7.3in
\TagsOnRight  
\topmatter
\title 3-manifold invariants from cosets
\endtitle
\author Feng Xu \endauthor

\address{Department of Mathematics, University of Oklahoma, 601 Elm Ave,
Room 423, Norman, OK 73019}
\endaddress
\email{xufeng\@ math.ou.edu}
\endemail
\abstract{ 
We construct unitary modular  categories  
for a general class of coset conformal field theories based on 
our previous study of these theories in the  algebraic quantum
field theory framework using subfactor theory.
We also consider the calculations of  the corresponding 3-manifold
invariants. It is  shown  that under certain index conditions 
the link invaraints colored by the representations of coset factorize
into the products of the the link invaraints colored by the representations
of the two groups in the coset. 
But the 3-manifold invariants do not behave so simply in general
due to the nontrivial branching and selection rules of the coset. 
Examples  in the parafermion 
cosets and diagonal cosets  show that  3-manifold
invariants of the coset may be finer than the products of the   3-manifold
invariants associated with the two groups in the coset, and these two
invariants do not seem to be simply related
in some cases, for an example, in the cases when there are
issues of ``fixed point resolutions''.  In the later case our framework
provides a representation theory
understanding of the underlying  unitary modular  categories
 which has not been
obtained before.
}
\endabstract
\thanks
\endthanks 
\endtopmatter
\document   
            
\heading  Introduction \endheading
This is the fourth paper in a series of papers on algebraic coset 
conformal field theories (cf. [X1], [X2] and [X3]), devoted to the construction
and calculation of 
of  3-manifold invariants  associated
with  a general class of coset conformal field theories. \par
The history of quantum invariants started from the  striking
discovery of Vaughan Jones of 
a new polynomial invariant of classical knots and links
(cf. [J2]) by using his subfactor theory (cf. [J1]). The  quantum 3-manifold
invariants were predicted by E. Witten (cf. [Wi]) motivated by a quantum
field theory interpretation of Jones polynomial. 
There are several mathematical 
approaches to the construction of quantum 3-manifold
invariants. The original approach [RT] (also cf. [KM]) 
is based on theory of quantum groups at roots of unity, which has
a subtle tensor product structure (cf. [A])
among a distinguished class of finite dimensional modules.  The work of 
[TW1], [TW2] and [W2] relies on ideas from quantum groups and subfactors.
There are also work (cf. [EK1], [EK2] and references therein) 
based on subfactors but there are no general
calculations of the invariants.
The approach we take in this paper (cf. [R] and [FRS]) 
is very different from that of [RT], [TW1], [TW2] and
[W2]. We still have a finite 
set of modules  , but they are all infinite
dimensional. The tensor product in this framework is very natural: they
correspond to compositions of certain endomorphisms (cf . P. 533 of [C]). 
There are several analytical problems in this approach,
one of the most difficult ones is the proof of finite index of 
certain inclusions. For the examples considered in this paper, this 
problem is settled by [Wa], [X1] and [X7]. 
There are several advantages to this approach. The first is that
our construction  is automatically unitary \footnotemark\footnotetext{
The unitarity or positivity question
in the  case of quantum groups at certain roots
of unity is  not trivial  (cf. [X5] and
[W3]).}  since we are using von Neumann algebras. 
The second is 
that
one of the axioms (non-degeneracy of $S$-matrix) of modular 
categories (cf. \S1.1)
which is most difficult to check in other approaches, 
follows  from finiteness of certain index by [X3]. 
The third and perhaps the most
important one is that our framework provides a mathematical understanding
of some unitary modular
categories which has not been obtained by other approaches (cf. \S3.6). 
However the proof of the existence of invariants does not guarantee an 
explicit formula. One of the advantages in the approaches of [RT]
(also cf. [KM] and [CM]) is that explicit formulas are given for any
3-manifold.   The second goal of this paper
is to obtain explicit formulas of 3-manifold invariants 
for a large class of cosets.\par

Now we'd like to explain why  3-manifold invariants 
from  cosets may be interesting. 
Let $G$ be a simply connected  compact Lie group
and let $H\subset G$ be a Lie subgroup . 
Consider
representations of loop group 
$LG$ (cf. \S1.5) with positive energy at level $k$
\footnotemark\footnotetext{ When G is the direct product of simple groups,
$k$ is a multi-index, i.e., $k=(k_1,...,k_n)$, where $k_i\in \Bbb {N}$
corresponding to the level of the $i$-th simple group. The level of $LH$
is determined by the Dynkin indices of $H\subset G$. To save some writing
we write the coset as $H\subset G_k$.}. 
Assume the coset $H\subset G_k$ verifies the conditions of Th. A in
\S1.7. 
Let $M$ be a closed 
oriented 3-manifold, and  
denote by $\tau_{G/H}(M)$ 
(resp. $\tau_G(M), \tau_H(M)$) the 3-manifold invariants associated
to the coset $H\subset G_k$ (resp. groups $G, H$) whose existence is 
implied by Th. A.  
By Th. B in \S2.3, 
it is natural to
believe that $\tau_{G/H}(M)$ should contain the same information as
$\tau_G(M) \overline{\tau_H(M)}$, but the following example shows that
this is not true.\par
Take $G=SU(2)_1\times SU(2)_1$, and $H=SU(2)_2$ is the diagonal subgroup
of $G$ at level $2$.  It is well known that the coset is the critical
Ising model. By an easy adaption of the 
argument in \S7 of [KM], it is shown in Prop. 3.6.3  that
$$
\tau_{G/H}(M) = \sum_{\theta} \exp(-\frac{2\pi i\mu(M_\theta)}{16})
,
$$ where $\mu(M_\theta)$ is the $\mu$-invariant of $M$ with spin
structure $ \theta$, and the sum is over all spin structures (cf. 
\S7 of [KM]). This formula can be compared to $\tau_H(M)$, which is
$\tau_4(M)$ in the notation of [KM], and is given by
$$
\tau_{H}(M) = \sum_{\theta} \exp(-\frac{6\pi i\mu(M_\theta)}{16})
.
$$  Now compare $\tau_{G/H}(M)$ with $\tau_G(M) \overline{\tau_H(M)}$,
which in the notation of [KM] is 
$$
\tau_3(M)^2 \overline{\tau_4(M)}.
$$
By (1) of Th. 6.3 in [KM], if $\tau_3(M)\neq 0$, then there exists
an integer $\mu$ depending on $M$ such that 
$$
\mu(M_\theta)\equiv \mu \ \text{\rm mod} \ 4, \forall \theta, 
$$ and it follows that
$$
\tau_G(M) \overline{\tau_H(M)} = \tau_3(M)^2 \exp(\frac{\mu \pi i}{2})
\tau_{G/H}(M).
$$ Thus when $\tau_3(M)\neq 0$, 
$\tau_{G/H}(M)$ agrees with $\tau_G(M) \overline{\tau_H(M)}$
up to a homotopy invariant 
$\tau_3(M)^2$ and a phase $\exp(\frac{\mu \pi i}{2})$.  
However when $\tau_3(M)= 0$,
$\tau_G(M) \overline{\tau_H(M)}$ is always $0$, but $\tau_{G/H}(M)$ may
not be $0$, for example when $M= {\Bbb R}P^3$ (Note that $ {\Bbb R}P^3=
L(2,1)$). \par
The above example is a special case of Prop. 3.6.1. 
More examples in \S3 show that $\tau_{G/H}(M)$ may be finer than
$\tau_G(M) \overline{\tau_H(M)}$, and these two invariants do not 
seem to  be 
simply related when there are issues of ``fixed point resolutions'', 
a problem   resolved mathematically in [X2]. \par
The calculation in \S3 is based on Th. B of \S2,    
which states  that under certain index conditions 
the link invariants colored by the representations of coset factorize
into the products of the the link invariants colored by the representations
of the two groups in the coset. The proof is based on certain braiding
properties first appeared in [X4] and further analyzed in [BE3].
These properties are also used in a very interesting recent paper [BEK]. 
\par
The main new results of this paper are Th. A, Th. B and formulas
(3.4.5), (3.5.2) and (3.6.1). \par
Here is a more detailed account of the paper. In \S1.1
we recall the definitions of unitary modular category from [Tu]. 
In \S1.2 to \S1.5 we recall some basic definitions from 
algebraic conformal theories and results which will be used in this
paper. These sections are included to set up notations and for the
benefit of reader with less operator algebra background. In \S1.6 
we discuss the nondegeneracy condition of cosets based on 
[X3]. 
In \S1.7 we prove (Th. A) that a general class of cosets naturally
give rise to unitary modular categories by checking all 
the axioms listed in \S1.1. \par
In \S2.1 we recall some definitions and results from [X4] and provide
a relation between [X4] and
[BE1], [BE2]. 
 In \S2.2 the relative braidings introduced in [BE3] are incorporated
into our formalism of cosets (Prop. 2.3.1) by using \S2.1, and   
these properties are used  to  prove 
factorization of framed link invariants under certain index conditions 
(Th. B) in \S2.3. \par 
\S3 are
 applications of the results in \S1 and \S2. In \S3.1 we show
that our invariants in the case of a type $A$ group is the same
as those of [TW1], [W2] in the type $A$ case. The proof is based on
certain properties of our invariants (lemma 1.7.4, 1.7.5) similar
to the cabling idea introduced in [W2].  In \S3.2 we give a different
proof of the symmetry principle ([KM], [KT1]).  In \S3.3 we give
a different proof of the result in [KT2] on the level-rank duality
of type $A$ invariants as an almost immediate
application of Th. B.  In \S3.4  we calculate invariants associated
with the simple current extensions of affine $su(N)$ based on [BE2]
and [BE3].  In \S3.5 we consider the invariants associated with
the parafermion cosets and observe that even in the simplest 
parafermion coset
our invariants may not always be related to the product of invariants
associated to the two groups in the coset. In \S3.6 formula for the
invariants associated with the diagonal cosets of type $A$ are 
obtained. It is shown that for a special class of  diagonal cosets
known as coset $W$ algebras the corresponding invariants are 
related to  the product of invariants
associated to the two groups in the coset by using the symmetry
principle (Prop. 3.6.1 and Cor. 3.6.2). 
However there does not appear to be such a relation
in general, especially when there are nontrivial fixed points
under certain diagram automorphisms. 
In \S3.7 a ``Maverick'' coset
is considered which does not verify the conditions of Th. B, but
we identify the corresponding invariants with that of a 
parafermion coset and a coset $W$ algebra by using the ideas of 
lemma 1.7.5 and lemma 3.1.1. Finally in \S3.8 a question is raised 
about the pertubative aspects of invariants from cosets.\par
\heading \S 1. Unitary modular categories from cosets \endheading
\subheading{\S1.1  Unitary modular categories}
Let us first recall the definitions of unitary   modular tensor category
from I.1, II.5 of [Tu].\par
A ribbon category is a monoidal category 
(assumed to be strict in the following)
equipped with a braiding $c$, a 
twist $\theta$, and a compatible duality $(*, b,d)$
over a ground field $K$. Let us explain the 
conditions on these three ingredients.  A braiding in a monoidal category
consists of a family of isomorphisms 
$$
c=\{ c_{V,W}: V\otimes W\rightarrow W\otimes V \},
$$ where $V,W$ run over all objects in the category, such that for
any three objects $U,V,W$ we have
$$
c_{U, V\otimes W} = (id_V \otimes c_{U,W}) (c_{U,V} \otimes id_W),
c_{U\otimes V,W} = (c_{U,W} \otimes id_V \otimes ) (id_U \otimes 
c_{U,V}) \tag 1.1.1
$$ The naturality of the isomorphisms of $c$ means that for any
morphisms $f: V\rightarrow V', g: W\rightarrow W'$, we have
$$
(g\otimes f) c_{V,W}= c_{V',W'} (f\otimes g) \tag 1.1.2
$$ The Yang-Baxter Equation (YBE), i.e.
$$
(id_W \otimes c_{U,V})(c_{U,W} \otimes id_V)(id_U \otimes c_{V,W})
=( c_{V,W}\otimes id_U)(id_V\otimes c_{U,W} )( c_{U,V}\otimes id_W)
$$ is a consequence of (1.1) and (1.2) (cf. P. 19 of [Tu]). \par
The twist $\theta$ in a monoidal category with a braiding $c$ consists
of a natural family of isomorphisms 
$
\theta= \{ \theta_V: V\rightarrow V \}$, where $V$ runs over all objects
of the category, such that for any two objects $V,W$ we have
$$
\theta_{V\otimes W}= c_{W,V} c_{V,W} (\theta_V \otimes \theta_W) \tag 1.1.3
$$ The naturality of $\theta$ means that for any morphisms $f:
U\rightarrow
V$, we have:
$$
\theta_V f = f\theta_U \tag 1.1.4
$$  \par
Assume that to each object $V$ of the category there are associated
an object $V^*$ and two morphisms
$$
b_V: 1\rightarrow V\otimes V^*, d_V: V^*\otimes V\rightarrow 1,
$$  
where $1$ denotes the unit object in $C$. 
This is called a duality if the following identities are satisfied:
$$
(id_V \otimes d_V) (b_V \otimes id_V) = id_V,
( d_V\otimes id_{V^*}) ( id_{V^*} \otimes b_V) = id_{V^*} \tag 1.1.5
$$ We need one axiom relating the duality morphisms $b_V, d_V$ with
braiding and twist. We say that the duality is compatible with the 
braiding $c$ and the twist $\theta$ if for any object $V$ in the
category we have:
$$
(\theta_V \otimes id_{V^*}) b_V = (id_V \otimes \theta_{V^*}) b_V \tag 1.1.6
$$  \par
By a ribbon category, we mean a monoidal category equipped with 
a braiding $c$, a twist $\theta$, and a compatible duality 
$(*, b,d)$ which verifies (1.1) to (1.6). A ribbon category is 
called strict if its underlying monoidal category is strict, and is 
called abelian if  its underlying monoidal category is abelian.  
For an endomorphism $f: V\rightarrow V$ in a ribbon category, define its
trace $tr(f) \in K$  to be the following composition:
$$
tr(f) = d_V c_{V,V^*} ((\theta_Vf) \otimes id_{V^*}) b_V: 1\rightarrow
1 \tag 1.1.7
$$ \par
In general ribbon categories do not admit direct sums. It turns out that
instead of decompositions of objects into direct sums one may decompose their
identity endomorphisms. This leads to the following notion of domination. 
We say $V$ is dominated by $\{ V_i \}_{i\in I}$ if the images of the
pairings
$$
\{ (g,f)\rightarrow fg: Hom (V,V_i) \otimes_K Hom (V_i, V) \rightarrow
End(V) \}_I
$$ additively generate $End (V)$. \par
A modular category is a pair consisting of an abelian ribbon category
$C$ and a finite family $\{ V_i \}_{i\in I}$ of simple objects 
(i.e., $End(V_i)$ is a free $K$ module of rank 1) in $C$ satisfying 
the following four axioms:\par
(1) (Normalization axiom). There exists $0\in I$ such that $V_0 =1$; \par
(2) (Duality axiom). For any $i\in I$, there exists $i^*\in I$ such that
the object $V_{i^*}$ is isomorphic to $(V_i)^*$; \par
(3) (Axiom of domination). All objects of $C$ are dominated by the 
family $\{ V_i \}_{i\in I}$.\par
To formulate the next and last axiom we need some notation. For 
$i,j\in I$, set (we add script Tu to distinguish this matrix from
the $S$ matrix defined in \S1.5): 
$$
S_{i,j}^{Tu}= tr(c_{V_j, V_i} c_{V_i, V_j}).
$$\par
(4) ( Non-degeneracy axiom). The square matrix $S= (S_{i,j})_{i,j\in I}$
is invertible over $K$. \par
Let $C$ be a strict monoidal category. A conjugation in $C$ assigns
to each morphism $f: V\rightarrow W$ in $C$ a morphism 
$\bar f: W\rightarrow V$ so that the following identities hold:
$$
\overline{\bar{f}} = f, \overline{f\otimes g}= \bar f \otimes \bar g, 
\overline{fg} = \bar g \bar f. 
$$ In the second formula $f,g$ are arbitary morphisms in $C$ and in the
third formula $f,g$ are  arbitrary composible morphisms in $C$. \par
A hermitian ribbon category is a ribbon monoidal category $C$ endowed
with a conjugation $f\rightarrow \bar f$ satisfying the following 
conditions:
For any objects $V,W in C$, we have
$$
\bar c_{V,W} = (c_{V,W})^{-1}, \tag 1.1.8
$$ and for any object $V$ of $C$ 
$$
\bar \theta_V = (\theta_V)^{-1}, \bar b_{V} = d_V c_{V,V^*} (\theta_V
\otimes id_{V^*}), \bar d_{V} = (id_{V^*} \otimes \theta_V^{-1})
  c_{V,V^*}^{-1} b_V) \tag 1.1.9
$$ \par
A unitary modular category is a Hermitian modular category over 
${\Bbb C}$ such that for any morphism $f$ we have
$$
tr(f\bar f) \geq 0 \tag 1.1.10
$$ \par
Suppose $C$ is a strict 
unitary modular category with a set of simple $V_i, i\in I$  
objects as above.  When no confusion arises we will denote
$V_i$ simply by $i$ and $V_i^*$ by $\bar i$.  
 


We will use the following convention of [Tu], $V^{+1}:= V, V^{-1}:= V^*$.
A graphical calculus is introduced in I.1.6 of [Tu]. 
In this calculus, 
the dualities, braidings and twistings are represented by 
special diagrams called tangles. The braiding $c_{V,W}$ is represented by
a positive crossing and  $c_{V,W}^{-1}$ is represented by a negative 
crossing (One usually follows a right-hand rule when drawing these
crossings, cf. figure 1.6 on Page 25 of [Tu]).
A general morphism $f:V\rightarrow W$
which is not a tangle 
is represented by a  diagram called a coupon. 
If a morphism is represented by a diagram, the closure of this 
morphism is simply obtained by closing all the free ends of the 
diagram (cf. Page 43 of [Tu] for a pictorial representation).
For more details on this graphical calculus including proofs using 
these diagrams see 
I.1.6 of [Tu]. 

Suppose $L$ is a framed oriented link in  3-shpere
with $n$ components colored by $i_1,...,i_n$. We will denote
the corresponding  invariant 
as defined on P. 39 of [Tu]   by $L(i_1,...,i_n)\in
{\Bbb C} $ .  It is often convenient to represent $L(i_1,...,i_n)$
as the closure of a tangle $T_L$ which consists of only braidings and
twists. Applying further braidings if necessary, we will  assume that
$T_L \in End (I_1 I_2... I_n)$, where $I_j:= (i_j^{\epsilon_j})^{k_j},
k_j\in {\Bbb N}, j=1,...,n$ and $\epsilon_j=1$ (resp. $\epsilon_j=-1$)
if the $j$-th component of $L$ is oriented anti-clockwise ( resp.
clockwise ). Recall from the above that $i^{-1}=  \bar i.$\par

Define $d(i):=tr(id_i)$. Define (cf. P. 76 and 
P. 115 of [Tu]) $D_C> 0$ (called rank of $C$) such that
$$
D_C^2:= \sum_i d(i)^2. \tag 1.1.11
$$ Note since $V_i$ is irreducible, the twist $\theta_{V_i}$ must act
on  $V_i$ as multiplication by a nonzero  element $\omega(i) \in {\Bbb C}$.
Set (cf. Page 21 of [MW])
$$
k_C^3 :=  D_C^{-1} \sum_i \omega(i) d(i)^2,
$$ 
then 
$$
D_C^{-1} \sum_i \omega(i)^{-1} d(i)^2 = k_C^{-3}. \tag 1.1.12
$$ 
Given a closed oriented 3-manifold $M_L$ obtained by 
surgery on a framed oriented link $L$ (with $n$ components) 
in the 3-sphere. The 3-manifold invariant $\tau_C(M_L) $ is 
defined by the following formula
$$
\tau_C(M_L):= k_C^{3(b_+(L) -b_-(L))} D_C^{-n}\sum_{i_1 \in I,...,i_n\in I}
d(i_1)... d(i_n) L(i_1,...,i_n),  \tag 1.1.13
$$ where $b_+(L)$ and  $b_-(L)$ are the numbers of positive and negative
eigenvalues of the linking matrix of $L$ (this matrix includes the
framings on the diagonal). The normalization of the 
invariant (1.11) differs from that of [Tu] but agrees [MW]. It
is normalized so that it is multiplicative under connected 
summing; in particular $\tau_C(S^3)=1$. \par
For the constraints of unitarity on the invariants (1.1.13), see
Page 224 of [Tu]. \par

Representation theory of certain 
quantum groups at certain roots of unity provides
examples of unitary modular category (cf. Chap. XII of [Tu] and
[W3]). In such examples the most difficult axioms to check are
the non-degeneracy axioms and positivity (1.10) (cf. [W2], [W3]
and [MW]). In \S1.7  we will show that
algebraic coset conformal field theories as formulated in [X1], [X3]
provides a new class of unitary modular categories. \par

\subheading{ \S1.2 Sectors and Correspondences }
Let us first recall some definitions from [L1] and [L2].

Let $M$, $N$ be von Neumann algebras, that we always assume to have
separable preduals, and ${\Cal H}$ a $M-N$ correspondence, namely
${\Cal H}$ is a (separable) Hilbert space, where $M$ acts on the left,
$N$ acts on the right and the actions are normal.
 
We denote by $x\xi y$, $x\in M$, $y\in N$, $\xi \in {\Cal H}$ the relative
actions.
 
The trivial $M-M$ correspondence is the Hilbert space $L^2(M)$ with the  
standard actions given by the modular theory
$$
x\xi y = xJy^*J\xi \, , \qquad
x, y\in M\, ,\qquad
\xi \in L^2(M) \, ,
$$
where $J$ is the modular conjugation of $M$; the unitary correspondence is
well defined modulo unitary equivalence.
 
If $\rho $ is a normal homomorphism of $M$ into $M$ we let
$H_\rho $ be the Hilbert space $L^2(M)$ with actions:
$x\cdot \xi \cdot y\equiv \rho (x) \cdot\xi \cdot y$, $x\in M$,
$y\in M$, $\xi \in L^2(M)$. Denote by $\text{\rm End}(M)$ the semigroup of
the endomorphism of $M$ and $\text{\rm Corr}(M)$ the set of all $M-M$
correspondences. The following proposition is proved in \cite{L2}
(Corollary 2.2 in \cite{L2}, also cf. [C]).
\proclaim{Proposition 1.1.1}  Let $M$ be an infinite factor. There exists a
bijection between the
unitary equivalence classes of $\text{\rm End}(M)$ and
 the
unitary equivalence classes of $\text{\rm Corr}(M)$, i.e., given
$\rho $, $\rho ^\prime \in \text{\rm End}(M)$, $H_\rho $ is
unitarily equivalent to
$H_{\rho ^\prime }$ iff there exists a unitary $u\in M$ with
$\rho ^\prime (x) = u\rho (x)u^*$.
\endproclaim
 
Let $\text{\rm Sect}(M)$ denote the quotient of $\text{\rm End}(M)$ modulo
unitary equivalence in $M$ as in Proposition 1.1.  We call sectors the    
elements of the semigroup $\text{\rm Sect}(M)$; if
$\rho \in \text{\rm End}(M)$ we denote by $[\rho ]$ its class in
$\text{\rm Sect}(M)$.  By Proposition 2.2 $\text{\rm Sect}(M)$ may be
naturally identified with $\text{\rm Corr}(M)^\sim $ the quotient of
$\text{\rm Corr}(M)$ modulo unitary equivalence. It follows from
\cite{L1} and \cite{L2} that $\text{\rm Sect}(M)$, with $M$ a properly
infinite (on Hilbert space ${\Cal H}$) von Neumann algebra, is endowed
with a natural involution $\theta \rightarrow \bar \theta $. 
The tensor product of correspondences correspond to the composition of
sectors. \par
Suppose $\rho \in \text{\rm End}(M)$ is given together with a normal
faithful conditional expectation
$\epsilon:
M\rightarrow \rho(M)$.  We define a number $d_\epsilon$ (possibly
$\infty$) such that (cf. [PP]) :
$$
d_\epsilon^{-2} :=\text{\rm Max} \{ \lambda \in [0, +\infty)|
\epsilon (m_+) \geq \lambda m_+, \forall m_+ \in M_+
\}$$\par
Now assume  $\rho \in \text{\rm End}(M)$ is given together with a normal
faithful conditional expectation
$\epsilon:
M\rightarrow \rho(M)$, and assume  $d_\epsilon < +\infty$. We define
$$
d = \text{\rm Min}_\epsilon \{ d_\epsilon \}.
$$  $d$ is called the statistical dimension of  $\rho$. It is clear
from the definition that  the statistical dimension  of  $\rho$ depends only  
on the unitary equivalence classes  of  $\rho$. We will denote
the statistical dimension of $\rho$ by $d(\rho)$ (In [X4] it is denoted
by $d_\rho$, here we choose $d(\rho)$ to avoid the confusion with 
the duality $d_V$ in (1.1.5)). $d(\rho)^2$ is called the minimal index
of $\rho$. \par
The properties of the statistical dimension can be found in
[L1], [L2] , [L3] and [K].\par
Denote by $\text{\rm Sect}_0(M)$ those elements of
$\text{\rm Sect}(M)$ with finite statistical dimensions.
For $\lambda $, $\mu \in \text{\rm Sect}_0(M)$, let
$\text{\rm Hom}(\lambda , \mu )$ denote the space of intertwiners from
$\lambda $ to $\mu $, i.e. $a\in \text{\rm Hom}(\lambda , \mu )$ iff
$a \lambda (x) = \mu (x) a $ for any $x \in M$.
$\text{\rm Hom}(\lambda , \mu )$  is a finite dimensional vector
space and we use $\langle  \lambda , \mu \rangle$ to denote
the dimension of this space.  $\langle  \lambda , \mu \rangle$
depends
only on $[\lambda ]$ and $[\mu ]$. Moreover we have
$\langle \nu \lambda , \mu \rangle =
\langle \lambda , \bar \nu \mu \rangle $,
$\langle \nu \lambda , \mu \rangle
= \langle \lambda , \mu \lambda \rangle $ which follows from Frobenius
duality (See \cite{L2} or \cite{Y}).  We will also use the following
notations: if $\mu $ is a subsector of $\lambda $, we will write as
$\mu \prec \lambda $  or $\lambda \succ \mu $.   
Recall (cf.[L2]) for each $\rho\in \text{\rm End(M)}$ and its conjugate
$\bar\rho$ with finite minimal index, there exists $R_\rho \in 
\text{\rm Hom}(id , \bar \rho \rho )$ and $ \bar R_\rho \in
\text{\rm Hom}(id , \rho \bar\rho )$ such that
$$
\align
&\bar R_\rho^* \rho(R_\rho) = 1, \   R_\rho^* \bar \rho(\bar R_\rho) = 1 
\\
& ||R_\rho||= ||\bar R_\rho||=\sqrt {d_\rho} \tag 1.2.1
\endalign
$$
The choice of  $R_\rho, \bar R_\rho$ above is unique up to a phase
only if $\rho$ is irreducible.  However  
 the  minimal
left inverse $\phi_\rho$ of $\rho$  defined by
$$
\phi_\rho (m) = R_\rho^* \bar\rho(m) R_\rho
$$ is unique and depends only one $\rho$. \par
Recall from [L3] that 
$$
\phi_{\rho_1\rho_2} =  \phi_{\rho_2}\phi_{\rho_1}.
$$ For any $f\in End(\rho)$, it is convenient to define an
unnormalized trace as
$$
tr_\rho (f):= d(\rho) \phi_\rho (f),
$$ and  this indeed satisfies tracial properties due to the following:
\proclaim{Lemma 1.2.1}
If $f\in Hom (\rho, \mu), g\in Hom (\mu, \rho)$, then
$$
tr_\rho (gf) = 
tr_\mu (fg).
$$
\endproclaim
\demo{Proof}
First assume that $\mu$ is irreducible. We can assume $\mu$ is a subsector
of $\rho$ since otherwise $f=0,g=0$ and the lemma holds trivially. 
Choose $h_i\in Hom(\mu, \rho)$
such that $h_i^* h_j =\delta_{i,j} id_\mu, \sum_{i} h_i h_i^*=id_\rho$.
Then $g= \sum_{i} h_ih_i^*g, f=\sum_{i} f h_ih_i^*$. Note that
$h_i^*g, f h_i \in End(\mu)$, and therefore are complex numbers
since $\mu$ is assumed to be irreducible. 
Also note by [L2] we have $ tr_\rho (h_i h_i^*) = d(\mu)$. 
We have
$$
\align
tr_\rho (gf) &= tr_\rho ( \sum_{i,j} h_ih_i^*g  f h_jh_j^*) \\
&= \sum_{i,j} h_i^*g  f h_j  tr_\rho (  h_ih_j^*) \\
&= \sum_{i} h_i^*g  f h_i tr_\rho (h_i h_i^*)\\
&= \sum_{i} h_i^*g  f h_i d(\mu)\\
&= \sum_{i,j} tr_\mu( f h_i h_i^* h_j  h_j*g) \\
&= tr_\mu (fg).
\endalign
$$ If $\mu$ is not irreducible, let {$\mu_i$} be a set of irreducible
subsectors of $\mu$, $T_i\in Hom (\mu_i, \mu)$ such that
$ T_i^*T_j = \delta_{i,j} id_j, \sum_{i} T_i T_i^* =id_\mu$. Note
that $ T_i^*f\in Hom(\rho, \mu_i), gT_i\in Hom(\mu_i,\rho)$. By the
previous argument we have
$$
tr_\rho (g T_i  T_i^*f) = tr_{\mu_i}(  T_i^*f g T_i)
= d(\mu_i)  T_i^*f g T_i.
$$ So
$$
\align
tr_\rho (gf) &= tr_\rho ( \sum_{i,j} gT_i T_i^* T_j T_j^*  f ) \\
&=  \sum_{i}  tr_\rho(gT_i  T_i^*f) \\
&=d(\mu_i)  T_i^*f g T_i \\
&=\sum_{i,j} tr_\mu (  T_i T_j^*)  T_i^*f g T_j \\
&= tr_\mu (gf).
\endalign
$$
\enddemo
\qed
\par
We will drop the subscript of tr in the following since it can
always be recovered by tracing the domain of morphisms. \par
\subheading {\S1.3   Conformal precosheaves on $S^1$ and 
superselection structures} 
In this section we recall the basic properties enjoyed by the family of
the von Neumann algebras associated with a conformal Quantum Field Theory
on $S^1$. All the propositions in this section
and \S 1.3 are proved in \cite{GL1}.\par
By an {\it interval} in this section 
 we shall always mean an open  connected subset $I$
of $S^1$ such that $I$ and the interior $I^\prime $ of its complement are
non-empty.  We shall denote by  ${\Cal I}$ the set of intervals in $S^1$. \par
We shall denote by  $PSL(2, {\bold R})$ the group of
  conformal transformations on the complex plane
that preserve the orientation and leave the unit circle $S^1$ globally
invariant.  Denote by ${\bold G}$
the universal covering group of $PSL(2, {\bold R})$.  Notice that  ${\bold G}$
is a simple Lie group and has a natural action on the  unit circle $S^1$. \par
Denote by  $R(\vartheta )$  the (lifting to ${\bold G}$ of the) rotation by  
an angle $\vartheta $.   

A conformal
precosheaf ${\Cal A}$ of von Neumann algebras on the intervals of $S^1$
is a map
$$
I\rightarrow {\Cal A}(I)
$$
from ${\Cal I}$ to the von Neumann algebras on a Hilbert space
${\Cal H}$ that verifies the following property:
\vskip .1in
\noindent
{\bf A. Isotony}.  If $I_1$, $I_2$ are intervals and
$I_1 \subset I_2$, then
$$
{\Cal A}(I_1) \subset {\Cal A}(I_2)\, .
$$

{\bf B. Conformal invariance}.  There is a unitary representation $U$ of
${\bold G}$ (the universal covering group of $PSL(2, {\bold R})$) on
${\Cal H}$ such that
$$
U(g){\Cal A}(I)U(g)^* = {\Cal A}(gI)\, , \qquad
g\in {\bold G}, \quad I\in {\Cal I} \, .
$$                                       

{\bf C. Positivity of the energy}.  The generator of the rotation subgroup
$U(R)(\cdot)$ is positive.

{\bf D.  Locality}.  If $I_0$, $I$ are disjoint intervals then
${\Cal A}(I_0)$ and ${\Cal A}(I)$ commute.
 
The lattice symbol $\vee $ will denote `the von Neumann algebra generated
by'.
 
{\bf E. Existence of the vacuum}.  There exists a unit vector
$\Omega $ (vacuum vector) which is $U({\bold G})$-invariant and cyclic for
$\vee _{I\in {\Cal I}}{\Cal A}(I)$.

A  conformal precosheaf is called irreducible if it also satisfies the
following:
 
{\bf F.  Uniqueness of the vacuum (or irreducibility)}.  The only
$U({\bold G})$-invariant vectors are the scalar multiples of $\Omega $. 

By Prop. 1.1 of [GL] a conformal precosheaf satisfies {\it Haag duality},
i.e., $ {\Cal A}(I)'=  {\Cal A}(I')$. \par
A conformal precosheaf  ${\Cal A}$  is called {\it strongly additive}
if it satisfies the following property: let $I_1,I_2$ be the connected
components of the interval $I$ with one internal point removed, then
$$
{\Cal A}(I)= {\Cal A}(I_1)\vee {\Cal A}(I_2),
$$ and ${\Cal A}$ is called {\it split} if ${\Cal A}(I_1)\vee {\Cal A}(I_2)
$ is 
naturally isomorphic to the tensor product of von Neumann algebras 
${\Cal A}(I_1)\otimes {\Cal A}(I_2)$ with $\bar I_1\cap \bar I_2=\emptyset$. 
According to [KLM], ${\Cal A}$ is called {\it completely rational},
shortly $\mu$-rational, if ${\Cal A}$ is strongly additive and
the minimal index of the inclusion ${\Cal A}(E)\subset {\Cal A}(E')'$ 
is finite, where $E=I_1\cup I_2, \bar I_1\cap \bar I_2=\emptyset .$
In this case the index of  ${\Cal A}(E)\subset {\Cal A}(E')'$  is called the 
{\it $\mu$-index} of  ${\Cal A}$, usually denoted by $\mu_{{\Cal A}}$.\par

Let ${\Cal A}$ be  an irreducible conformal precosheaf of
von Neumann algebras as defined above.  Our next goal is to define
 covariant {\it representation} of ${\Cal A}$ and the associated
concepts which will be important later on. Again all the definitions
and propositions are given in [GL] and we  refer the reader
to  [GL] for further details . \par
 
A covariant {\it representation} $\pi $ of
${\Cal A}$ is a family of representations $\pi _I$ of the
von Neumann algebras ${\Cal A}(I)$, $I\in {\Cal I}$, on a
Hilbert space ${\Cal H}_\pi $ and a unitary representation
$U_\pi $ of the covering group ${\bold G}$ of $PSL(2, {\bold R})$, with
{\it positive energy}, i.e. the generator of the rotation unitary subgroup
has positive generator, such that the following properties hold:
$$
\align
I\supset \bar I \Rightarrow \pi _{\bar I} \mid _{{\Cal A}(I)} 
= \pi _I \quad &\text{\rm (isotony)} \\
\text{\rm ad}U_\pi (g) \cdot \pi _I = \pi _{gI}\cdot
\text{\rm ad}U(g) &\text{\rm (covariance)}\, .
\endalign
$$
A unitary equivalence class of representations of ${\Cal A}$ is called
a {\it superselection sector}.
 
Assuming ${\Cal H}_\pi $ to be separable, the representations
$\pi _I$ are normal because the ${\Cal A}(I)$'s are factors .
Therefore for any given $I_0$, $\pi _{I_0^\prime }$ is unitarily
equivalent $\text{\rm id}_{{\Cal A}(I_0^\prime )}$ because
${\Cal A}(I_0^\prime )$ is a type III factor. By identifying
${\Cal H}_\pi $ and ${\Cal H}$, we can thus assume that $\pi $ is
localized in a given interval $I_0 \in {\Cal I}$, i.e.
$\pi _{I_0^\prime } = \text{\rm id}_{{\Cal A}(I_0^\prime )}$ .
By Haag duality we then have $\pi _I ({\Cal A}(I)) \subset {\Cal A}(I)$ if
$I\supset I_0$. In other words, given $I_0 \in {\Cal I}$ we can choose in
the same sector of $\pi $ a {\it localized endomorphism} with localization
support in $I_0$, namely a representation $\rho $ equivalent to $\pi $
such that
$$
I\in {\Cal I}, I \supset I_0 \Rightarrow \rho _I \in
\text{\rm End}\ {\Cal A}(I)\, , \qquad \rho _{I_0^\prime }
= \text{\rm id}_{I_0^\prime } \, .   
$$
\noindent
To capture the global point of view we may consider the {\it universal
algebra} $C^*({\Cal A})$.  Recall that
$C^*({\Cal A})$ is a $C^*$-algebra canonically associated with the
precosheaf ${\Cal A}$ (see \cite{Fre}).   $C^*({\Cal A})$
has the following properties:  there are injective embeddings
$\iota _I: {\Cal A}(I) \rightarrow C^*({\Cal A})$ so that the local
von Neumann algebras ${\Cal A}(I)$, $I\in {\Cal I}$, are identified with
subalgebras of $C^*({\Cal A})$ and generate all together a dense
$*$-subalgebra of $C^*({\Cal A})$, and every representation of the
precosheaf ${\Cal A}$ factors through a representation of $C^*({\Cal A})$.
Conversely any representation of $C^*({\Cal A})$ restricts to a
representation of ${\Cal A}$.  The vacuum representation
$\pi _0$ of $C^*({\Cal A})$ corresponds to the identity representation of
${\Cal A}$ on ${\Cal H}$, thus $\pi _0$ acts identically on the local von
Neumann algebras. We shall often drop the symbols $\iota _I$ and
$\pi _0$ when no confusion arises.
 
By the universality property, for each $g\in PSL(2, {\bold R})$ the
isomorphism $\text{\rm ad}U(g):{\Cal A}(I) \rightarrow {\Cal A}(gI)$,
$I\in {\Cal I}$ lifts to an automorphism $\alpha_g$ of $C^*({\Cal A})$.
We shall  lift the map $g\rightarrow \alpha_g$ to a representation,
still denoted by $\alpha $, of the universal covering group
${\bold G}$ of $PSL(2, {\bold R})$ by automorphisms of
$C^*({\Cal A})$.
 
The covariance property for an endomorphism $\rho $ of
$C^*({\Cal A})$ localized in $I_0$ means that         
$\alpha_g \cdot \rho \cdot \alpha _{g^{-1}}$ is
$$
\text{\rm ad}z_\rho (g)^* \cdot \rho = \alpha _g \cdot \rho
\cdot \alpha _{g^{-1}} \qquad g\in {\bold G}
$$
for a suitable unitary $z_\rho (g) \in C^*({\Cal A})$.
We define
$$\rho_g = \alpha _g \cdot \rho
\cdot \alpha _{g^{-1}} \qquad, g\in {\bold G}.
$$
$\rho_{g,J}$ is the restriction of $\rho_g $ to ${\Cal A}(J)$.
The map
$g\rightarrow z_\rho (g)$ can be chosen to be a localized
$\alpha $-cocycle, i.e.
$$
\align
z_\rho (g) &\in {\Cal A}(I_0 \cup gI_0) \quad
\forall g \in {\bold G}:I_0 \cup gI_0 \in {\Cal I} \\
z_\rho (gh) &= z_\rho (g) \alpha _g (z_\rho (h))\, ,\qquad
g, h \in {\bold G} \, .
\endalign
$$
To compare with the result of \cite{FG},
let us define:
$$\Gamma_\rho(g) = \pi _0 (z_\rho (g)^*).$$ This notation will be
used in \S1.4. \par
 
 
An endomorphism of $C^*({\Cal A})$ localized in an interval $I_0$ is said
to have {\it finite index} if $\rho _I (=\rho |_{{\Cal A}(I)})$ has finite
index, $I_0 \subset I$ .  The index is independent of
 $I$
due to the following (See Prop. 2.1 of \cite{GL1})  
\proclaim{1.3.1 Proposition}  Let $\rho $ be an endomorphism localized in
the interval $I_0$. Then the index
$\text{\rm Ind}(\rho ) := \text{\rm Ind}(\rho _I)$, the minimal index of
$\rho _I$, does not depend on the interval $I\supset I_0$.
\endproclaim

The following proposition is Prop. 2.2 of  \cite{GL1}:
\proclaim{1.3.2 Proposition}  Let $\rho $ be a covariant (not necessarily
irreducible) endomorphism with finite index. Then the representation
$U_\rho $ described before is unique. In particular, any irreducible
component of $\rho $ is a covariant endomorphism.
\endproclaim
 

By the above proposition the {\it univalence} of an endomorphism
$\rho$ is well defined by
$$
\omega(\rho) = U_\rho (2\pi ),
$$  
since
$U_{\rho ^\prime }(g):= \pi _0 (u)U_\rho (g) \pi _0 (u)^*$, where
$\rho ^\prime (\cdot ) := u\rho (\cdot ) u^*$,
$u\in C^*({\Cal A})$,  $\omega(\rho) $ depends only on the superselection
class of $\rho $. \par
When  $\rho $ is irreducible, $\omega(\rho) $ is a complex number of
modulus one since by definition
$\omega(\rho) $ belongs to $\pi (C^*({\Cal A}))^\prime $
and  $\rho $ is irreducible, and we have
$$
 \omega(\rho)= e^{2\pi i \Delta_\rho }
$$
with $\Delta_\rho $ the lowest weight of $U_\rho $.   $\Delta_\rho $ is
also referred to as {\it conformal dimension}.  Formulas for certain 
conformal dimensions may be found in (1.5.10).  \par

Let $\rho _1$, $\rho _2$ be endomorphisms of an algebra ${\Cal B}$.
Recall from \S1.1 that their
intertwiner space is defined by
$$ \text{\rm Hom}
(\rho _1 , \rho _2) = \{ T\in {\Cal B}:
\rho _2 (x)T = T \rho_1 (x), \quad x \in {\Cal B}\}.
$$
In case ${\Cal B}=C^*({\Cal A})$, $\rho _i$ localized in the interval    
$I_i$ and $T\in (\rho _1, \rho _2 )$, then $\pi _0(T)$ is an intertwiner
between the representations $\pi _0 \cdot \rho _i$. If
$I\supset I_1 \cup I_2$, then by Haag duality its embedding
$\iota _I \cdot \pi _0 (T)$ is still an intertwiner in
$(\rho _1, \rho _2)$ and a local operator. We shall denote by
$(\rho _1 , \rho _2)_I$ the space of such local intertwiners
$$
(\rho _1 , \rho _2)_I = (\rho _1 , \rho _2) \cap
{\Cal A}(I)\, .
$$
If $I_1$ and $I_2$ are disjoint, we may cover $I_1 \cup I_2$ by an
interval $I$ in two ways: we adopt the convention that, unless otherwise
specified, a local {\it intertwiner} is an element of
$(\rho _1 , \rho _2)_I$ where $I_2$ follows $I_1$ inside $I$ in the
clockwise sense.
 
We now define the statistics. Given the endomorphism $\rho $ of
${\Cal A}$ localized in $I\in {\Cal I}$, choose an equivalent endomorphism
$\rho _0$ localized in an interval $I_0\in {\Cal I}$ with
$\bar I _0 \cap \bar I = \emptyset $ and let $u$ be a local intertwiner in
$(\rho , \rho _0)$ as above, namely $u\in (\rho , \rho _0)_{\tilde I}$ with
$I_0$ following clockwise $I$ inside $\tilde I$.
 
The {\it statistics operator} $\sigma := u^* \rho (u)
= u^* \rho _{\tilde I}(u)$ belongs to $(\rho _{\tilde I}^2 ,
\rho _{\tilde I}^2)$.                 
Recall that if $\rho $ is an endomorphism of a $C^*$-algebra
${\Cal B}$, a {\it left inverse} of $\rho $ is a completely positive map
$ \phi $ from ${\Cal B}$ to itself such that
$\phi \cdot \rho = \text{\rm id}$.
 
It follows from  Cor.2.12 of [GL] that  there exists
a unique left inverse $ \phi $ of $\rho $ and that the
{\it statistics parameter}
$$
\lambda _\rho :=  \phi (\sigma )
$$
depends only on the sector of $\rho $.

By [GL] The {\it statistical dimension} $d(\rho )$ are given by
$$
d(\rho ) = |\lambda _\rho |^{-1}\,
$$ and the {\it statistics phase} 
$\kappa _\rho $ are  defined by
$$
\qquad \kappa _\rho = \frac {\lambda _\rho }{|\lambda _\rho |}  .
$$
In \cite{GL1}, the following  theorem is proved:
\proclaim{Theorem 1.3.3 (Spin-Statistics Theorem)}\par
If $\rho $ is irreducible and has finite index, then
$\kappa _\rho = \omega(\rho)$.
\endproclaim
We shall need the following easy corollary.
\proclaim{Corollary 1.3.4}
Assume that $\rho$ has finite index. Then
$d(\rho) \lambda _\rho$  is a constant $a$  times identity if
and only if for any irreducible subsector $x$ of $ \rho$, 
$\omega(x)= a$. \par
\endproclaim
\demo{Proof}
Ad (1). By Prop. 6.3 of [DHR] (also cf. [L1]) we have
$$
d(\rho) \lambda _\rho = \sum_{x} d_x \lambda_x E_x = \sum_{x} \kappa _x
E_x ,$$ where the summation is over all different irreducible 
subsectors $x$ of $\rho$, and $E_x$ is the minimal abelian projection
in $End(\rho)$ corresponding to $x$. By the spin-statistics theorem
we have
$$
d(\rho) \lambda _\rho = \sum_{x} \omega(x) E_x 
,$$ and (1) follows immediately. 
\enddemo
\qed
\par
\vskip .1in
\noindent      
\subheading { 1.4 Coherence equations}
 
In this section, we assume $\Delta $ is a set of localized covairant
endomorphism of ${\Cal A}$ with localization support in $I_0$.
Let $h$, $g$ be elements of ${\bold G}$.  We assume
$h I_0 \cap I_0 = \emptyset $, $g I_0 \cap I_0 =  \emptyset  $, $h I_0
\cap g I _0 =  \emptyset  $. Choose $J_1$, $J_2\in {\Cal I}$ such that
$J_1\cup J_2 \subsetneq S^1 $,
$J_1 \supset I_0 \cup g.I_0$, $J_2 \supset I_0 \cup h.I_0$,
$J_1 \cap h.I_0 =  \emptyset $, $J_2 \cap g.I_0 =  \emptyset  \,$ and
$ J_1 \cap J_2 = I_0 $. 
We assume in $J_1$ (resp. $J_2$), $g.I_0$ (resp. $h.I_0$) lies a clockwise
(resp. anti clockwise) from $I_0$.
The proofs of all the results of this section can be found in \S1.4 of
[X4]. \par 
\proclaim{Lemma 1.4.1}  For any $J\supset J_1 \cup J_2$, $J\in {\Cal I}$,
$\gamma, \lambda \in \Delta $ and $x\in  {\Cal A}(J)$,  we have
\roster
\item"{(0)}"
$\Gamma _\lambda (g) \in {\Cal A}(J_1)$.
\item"{(1)}"
$\Gamma _\lambda (g)^* \gamma _{J_1}
(\Gamma _\lambda (g))
\gamma _J \cdot \lambda _J (x)
= \lambda _J \cdot \gamma _J (x)
\Gamma _\lambda (g)^* \gamma _{J_1}
(\Gamma _\lambda (g) $.
\item"{(2)}"  $\Gamma _\lambda (g)^*
\gamma _{J_1} (\Gamma _\lambda (g))
= \lambda _{J_2} (\Gamma _\gamma (h)^*)
\Gamma _\gamma (h)$
\item"{(3)}"  $\Gamma _\gamma (g)^*
\gamma _{J_1} (\Gamma _\lambda (g)) \in {\Cal A}(I_0)$.
\endroster
\endproclaim  
Because the property (1) of Lemma 1.4.1,
$\Gamma _\lambda (g)^* \gamma _{J_1}
(\Gamma _\lambda (g))$ is called the braiding operator. 

We shall use $\sigma_{\gamma\lambda}$ to denote 
$\Gamma _\lambda (g)^* \gamma _{J_1}
(\Gamma _\lambda (g))$.  $\sigma_{\gamma\lambda}$ will be
called the positive braiding operator between 
$\gamma$ and $\lambda$.  We define
$$
\tilde \sigma _{\gamma  \lambda }:= \sigma _{\lambda \gamma }^*
$$ and 
$\tilde \sigma _{\gamma  \lambda }$ will be called the 
negative braiding operator between 
$\gamma$ and $\lambda$.
We are now ready to state the following equations. For simplicity
we will drop the subscript $I_0$ and write $\mu_{I_0}$ as $\mu$
for any $\mu \in \Delta$ in the following.
 
\proclaim{Proposition 1.4.2}  \text{\rm (1)}  Yang-Baxter-Equation (YBE)
$$
\sigma _{\mu  \gamma } \mu (\sigma _{\lambda \gamma })
\sigma _{\lambda \mu } = \gamma (\sigma _{\lambda \mu })
\sigma _{\lambda \gamma }\lambda (\sigma _{\mu  \gamma })\, .
$$        
\text{\rm (2)}  Braiding-Fusion-Equation (BFE)
 
For any $w\in \text{\rm Hom} (\mu \gamma  \delta )$
$$
\align
\sigma _{\lambda \delta } \lambda (w)
= w\mu (\sigma _{\lambda \gamma })
\sigma _{\lambda \mu }\tag a\\
\sigma _{\delta  \lambda }w = \lambda (w)
\sigma _{\mu  \lambda } \mu
(\sigma _{\gamma  \lambda }) \, .\tag b \\
\sigma _{ \delta,\lambda }^* \lambda (w)
= w\mu (\sigma _{ \gamma\lambda }^*)
\sigma _{\mu\lambda }^*\tag c\\
\sigma _{ \lambda \delta }^* \lambda (w)
= w\mu (\sigma _{ \gamma\lambda }^*)
\sigma _{\lambda\mu }^*\tag d
\endalign
$$
\endproclaim  
Suppose $\xi_1 \in I_{\xi _1} \subset J_1$,
$I_{\xi _1} \cap gI_1 \cap I_1 = \emptyset $,
$\xi _2 \in I_{\xi _2} \subset J_2$,  and
$I_{\xi _2} \cap I_1 \cap h.I_1 =  \emptyset $.
Here $g$, $h$, $J_1$, $J_2$ are defined as the beginning of this section.

It follows from (2) of Lemma 1.4.1 that:
$$
\gamma_{I_{\xi _1}^c} (\Gamma _\lambda (g))^*
\Gamma _\lambda (g) = \gamma_{J_2}
(\Gamma _\lambda (h)^*)\Gamma _\lambda (h)
= \sigma _{\lambda \gamma }\, .
$$

Hence $\sigma _{\lambda \gamma } \sigma _{\gamma , \lambda }
= \gamma _{I_{\xi _1}^c} (\Gamma _\lambda (g)^*)
\gamma_{I_{\xi _2}^c} (\Gamma _\lambda (g))$.

$\sigma _{\lambda \gamma }\sigma _{\gamma , \lambda }$ is called
the {\it monodromy operator}. \par
Let $T_e: \delta \rightarrow \gamma \lambda $ be an
intertwiner.  Recall $\omega(\rho) = U_\rho (2 \pi )$ 
is the univalence of a covariant
endomorphism. When $\rho $ is irreducible, $\omega(\rho)$ is a complex number.
We have:
\proclaim{Proposition 1.4.3(monodromy equation)}
If $\omega(\delta) $, $\omega(\gamma) $, $\omega(\lambda) $ are
complex numbers, then
$$                    
T_e^* \gamma_{I_{\xi _1}^c}(\Gamma _\lambda (g)^*)
\gamma_{I_{\xi _2}^c}(\Gamma _\lambda (g))
T_e = T_e^* \sigma _{\lambda \gamma }
\sigma _{\gamma  \lambda }
T_e = \frac {\omega(\delta) }{\omega(\lambda) \omega(\gamma) } \, .
$$
\endproclaim  
\subheading{1.5  Genus 0 and 1 modular matrices }
 Next we will recall some of the results of [Reh] (also cf. [FRS]) and
introduce
notations. \par
Let $[\rho_i]$ denote the equivalence classes of irreducible
superselection \par
sectors  in a finite set.
Suppose this set is closed under
conjugation
and composition. We will denote the conjugate of $[\rho_i]$ by
$[\rho_{\bar i}]$
and identity sector by $[1]$ if no confusion arises, and let
$N_{ij}^k = \langle [\rho_i][\rho_j], [\rho_k]\rangle $. 
The algebra generated by $ [\rho_i]$'s will be called 
{\it fusion rule algebra} and $N_{ij}^k $ will be referred to as
{\it fusion coefficients}. 
We will
denote by $\{T_e\}$ a basis of isometries in $\text {\rm
Hom}(\rho_k,\rho_i\rho_j)$.
Recall from \S1.3 the univalence of $\rho_i$  is  denoted by
$\omega_{\rho_i}$. 
Let $\phi_i$ be the unique minimal
left inverse of $\rho_i$, define:
$$
Y_{ij}:= d_{\rho_i}  d_{\rho_j} \phi_j (\sigma_{\rho_j \rho_i}
\sigma_{\rho_i \rho_j})^*, \tag 1.5.0
$$ where $\sigma_{\rho_j \rho_i}$ is the unitary braiding operator
defined in \S1.4. \par
We list two properties of $Y_{ij}$ (cf. (5.13), (5.14) of [Reh]) which
will be used later: 
$$
\align
Y_{ij} = Y_{ji} & = Y_{i\bar j}^* = Y_{\bar i \bar j}, 
\omega(i)= \omega(\bar i)
\tag 1.5.1 \\
Y_{ij} = \sum_k N_{ij}^k \frac{\omega_i\omega_j}{\omega_k} d_{\rho_k} \tag
1.5.2
\endalign
$$
 Define
$D_- := \sum_i d_{\rho_i}^2 \omega_{\rho_i}^{-1}$.
If the matrix $(Y_{ij})$ is invertible,
by proposition on P.351 of [Reh] $D_-$ satisfies
$|D_-|^2 = \sum_i d_{\rho_i}^2$.
Suppose $D_-= |D_-| \exp(i x), x\in {\Bbb R}$.
Define matrices
$$
S:= |D_-|^{-1} Y, T:=  C Diag(\omega_{\rho_i})
,\tag 1.5.3            
$$ where $C:= \exp(i \frac{x}{3}).$  Then these matrices satisfy the algebra:
$$
\align
SS^{\dag} & = TT^{\dag} =id, \tag 1.5.4  \\
TSTST&= S, \tag 1.5.5 \\
S^2 =\hat{C}, T\hat{C}=\hat{C}T=T, \tag 1.5.6
\endalign
$$
where $\hat{C}_{ij} = \delta_{i\bar j}$ is the conjugation matrix. Moreover
$$
N_{ij}^k = \sum_m \frac{S_{im} S_{jm} S_{km}^*}{S_{1m}}. \tag 1.5.7
$$
(1.5.7) is known as Verlinde formula. \par
We will refer the $S,T$ matrices
as defined in  (1.5.3)  as  {\bf genus 0 modular matrices} since
they are constructed from the fusions rules, monodromies and minimal
indices which can be thought as  genus 0 data associated to
a Conformal Field Theory (cf. [MS]). \par
Now let us consider an example which verifies (1.5.1) to (1.5.7) above.
Let $G= SU(N)$. We denote $LG$ the group of smooth maps
$f: S^1 \mapsto G$ under pointwise multiplication. The
diffeomorphism group of the circle $\text{\rm Diff} S^1 $ is
naturally a subgroup of $\text{\rm Aut}(LG)$ with the action given by
reparametrization. In particular the group of rotations
$\text{\rm Rot}S^1 \simeq U(1)$ acts on $LG$. We will be interested
in the projective unitary representation $\pi : LG \rightarrow U(H)$ that
are both irreducible and have positive energy. This means that $\pi $
should extend to $LG\ltimes \text{\rm Rot}\ S^1$ so that
$H=\oplus _{n\geq 0} H(n)$, where the $H(n)$ are the eigenspace
for the action of $\text{\rm Rot}S^1$, i.e.,
$r_\theta \xi = \exp^{i n \theta}$ for $\theta \in H(n)$ and
$\text{\rm dim}\ H(n) < \infty $ with $H(0) \neq 0$. It follows from
\cite{PS} that for fixed level $k$ which
is a positive integer, there are only finite number of such
irreducible representations indexed by the finite set 
$$
 P_{++}^{h}
= \bigg \{ \lambda \in P \mid \lambda
= \sum _{i=1, \cdots , N-1}
\lambda _i \Lambda _i , \lambda _i \geq 1\, ,
\sum _{i=1, \cdots , n-1}
\lambda _i < h \bigg \}
$$
where $P$ is the weight lattice of $SU(N)$ and $\Lambda _i$ are the
fundamental weights and $h=N+k$.  
Elements of $P_{++}^{h}$ are also denoted by 
$( \lambda_1,..., \lambda_{N-1})$. 
We will use
$1$ to denote the  representation  corresponding to $(1,1,...,1)$, referred
to
as {\it vacuum} representation
. For $\lambda , \mu , \nu \in  P_{++}^{K}$, define
$$
N_{\lambda \mu}^\nu  = \sum _{\delta \in P_{++}^{K} } \frac{S_{\lambda
\delta}
S_{\mu \delta} S_{\nu \delta}^*}{S_{1\delta}} \tag 1.5.8
$$
where $S_{\lambda\delta}$ is given
by the Kac-Peterson formula:     
$$
S_{\lambda \delta} = c \sum _{w\in S_N} \varepsilon _w \exp
(iw(\delta) \cdot \lambda 2 \pi /n). \tag 1.5.9
$$
Here  $\varepsilon _w = \text{\rm det}(w)$ and $c$ is a normalization
constant fixed by the requirement that $(S_{\lambda\delta})$
is an orthonormal system.
It is shown in \cite{Kac} P.288 that $N_{\lambda \mu}^\nu $ are
non-negative
integers. Moreover, define $ Gr(C_K)$
to be the ring whose basis are elements
of $ P_{++}^{K}$ with structure constants $N_{\lambda \mu}^\nu $.
  The natural involution $*$ on $ P_{++}^{K}$ is
defined by $\lambda \mapsto \lambda ^* =$ the conjugate of $\lambda $ as
representation of $SU(N)$. \par  
The irreducible positive energy representations of $ L SU(N)$ at level
$K$ give rise to an irreducible conformal precosheaf ${\Cal A_{G}}$
and                   
its covariant representations (cf. P. 362 of [X1]). Let us 
recall  how ${\Cal A_{G}}$ is defined. Let $\pi^0$ be the vacuum
 representation of $LG$  
on Hilbert space $H^0$. Then $A_{G}(I):= \pi^0(L_IG)''$,
where $\pi^0(L_IG)''$ denote the von Neumann algebra generated by 
all elements of the form $\pi^0(x), x\in LG, x={e} \text{\rm on} I^c$ with
${e}$ the identity element of $G$. 
 
This conformal precosheaf  is strongly additive, split 
(cf. [Wa]) and $\mu$-rational (cf. [X7]).
The unitary equivalent
classes of such  representations are the superselection sectors.
  We will use
$\lambda$ to denote such   representations. \par
For $\lambda$ irreducible, the univalence $\omega_\lambda$ is given by
an explicit formula .
Let us first
define
$$
\Delta_\lambda =
\frac {c_2(\lambda)}{K+N} \tag 1.5.10
$$ where $c_2(\lambda)$ is the value of
Casimir
operator on representation of $SU(N)$ labeled by dominant weight
$\lambda$ (cf. 1.4.1 of [KW]).
 $\Delta_\lambda$ is usually called the conformal dimension.

We have
$\omega_\lambda = \exp({2\pi i} \Delta_\lambda)$.
Note that $\omega_\lambda=\omega_{\bar\lambda}$. 
\par

Define the central charge (cf. 1.4.2 of [KW])
$$
C_G := \frac {K \text {\rm dim(G)}}{K+N} \tag 1.5.11
$$ and $T$ matrix as
$$
T=diag(\dot\omega_\lambda), \tag 1.5.12
$$ where $\dot\omega_\lambda = \omega_\lambda exp (\frac{-2\pi i
C_G}{24})
$.  By Th.13.8 of [Kac] $S$ matrix as defined in (1.5.9) and $T$ matrix
in (1.5.12) satisfy relation (1.5.4), (1.5.5) and (1.5.6).
 Since $S,T$ matrix defined in (1.5.8) and
(1.5.11) are related to the modular properties of characters which are
related to Genus 1 data of CFT (cf. [MS]), we shall call them
{\bf genus 1 modular matrices.}
\par                  
By Cor.1 in \S34 of [Wa],  The fusion ring generated by all
$\lambda \in   P_{++}^{(K)}$
is isomorphic to $ Gr(C_K)$, with structure constants $N_{\lambda
\mu}^\nu$ as defined in (1.5.8).
 One may therefore ask what are the $Y$ matrix
(cf. (1.5.0)) in this case. By using (1.5.2) and the formula for
 $N_{\lambda
\mu}^\nu$, a simple calculation shows:
$$
Y_{\lambda \mu} = \frac{S_{\lambda \mu}}{S_{1\mu}}
,$$ and it follows that $Y_{\lambda \mu}$ is nondegenerate, and $S,T$
matrices
as defined in (1.5.3) are indeed the same $S,T$ matrix defined in (1.5.8) and
(1.5.11),
which is a surprising fact. This fact is  referred to as
{\it genus 0 modular matrices coincide with genus 1 modular matrices}.
\par
The fusion coefficient for $N_{v\lambda}^{\mu}$ takes a particular
simple form when $v:=2\Lambda_1+ \Lambda_2+...+\Lambda_{N-1}$ corresponds
to the defining representation of $SU(N)$. Represent $\lambda$ by
$(\lambda_1,...,\lambda_{N-1})$, we have
$$
\align
v\lambda &= (\lambda_1+1, \lambda_2,..., \lambda_{N-1})
+ (\lambda_1-1, \lambda_2+1,..., \lambda_{N-1}) +
... + \\
&(\lambda_1, \lambda_2,..., 
\lambda_{N-2}-1, \lambda_{N-1}+1) \tag 1.5.13
\endalign
$$ 
where on the righthand side if a weight is not in $P_{++}^h$ it
is defined to be $0$.
\proclaim{Lemma 1.5.1}
Denote by $v:=2\Lambda_1+ \Lambda_2+...+\Lambda_{N-1}, a_i=\Lambda_1+...
+2\Lambda_i +...+ \Lambda_{N-1}, 
s_i= (i+1)\Lambda_1+ \Lambda_2+...+\Lambda_{N-1}, i=1,...,N-1 $. 
Then: \par
(1)
$[a_i]$ appears once and only once in $[v^i]$. All other subsectors 
except $[a_i]$ and $[s_i]$ appear with multiplicity at least 2.\par
(2) If $A$ is a map from the set $P_{++}^h$ to itself, and
$$
S_{A(\lambda)A(\mu)} =S_{{\lambda}{\mu}}, \forall \lambda,\mu\in P_{++}^h,
$$  and $A(a_i) = a_i, i=1,...,N-1$, then
$A=id$.
\endproclaim
\demo{Proof}
Ad (1). When $k=1$ the fusion ring is isomorphic to $Z_N$ and (1)
is trivial. Let us assume that $k\geq 2$. Let us prove a lightly
stronger statement by induction on
$i$. The statement is that for any $  i=1,...,N-1 $, (1) is true and in
addition the subsectors  of $[a_{i-1}v]:=[a_i] +[b_i]$ appear
($[a_{0}]=[id]$) in $[v^i]$, and  if $[s_i]$ appears in $[v^i]$, 
then the subsectors  of 
$[s_{i-1}v]:=[s_i] + [c_i]$ 
($[s_{0}]=[id]$)  appear in $[v^i]$. When $i=1$ the statement
is trivial. Assume it is proved for $1\leq i\leq {N-2}$, let us prove it
for $i+1$.  Note that denote by  $[a_{i}v]:=[a_{i+1}] +[b_{1+i}]$, then
$[b_i v]\succ
b_{1+i}$ by (1.5.13) , and similarly  if $[s_{i+1}]$ appears in $[v^{i+1}]$,
then  $[s_{i}v]:=[s_{i+1}] +[c_{1+i}]$, and $[c_i v]\succ
c_{1+i}$.  It is then easy to see by induction that the stronger 
statement as above is true for $i+1$. By induction (1) is proved.
\par
Ad(2):
By assumption and (1.5.8), we have 
$N_{A(\lambda) A(\mu)}^{A(\nu)} = N_{\lambda\mu}^\nu$. 
It follows that $a\rightarrow A(a)$
is a ring endomorphism of $Gr(C_k)$. Since $a_i, i=1,2,...,{N-1}$ are
generators of  $Gr(C_k)$ by \S34 of [Wa] and $A(a_i) = a_i, i=1,...,N-1$
by our assumption, it follows that $A=id$.
\enddemo
\qed
\par  

\subheading {1.6 Rationality and non-degeneracy of cosets} 
Let $G$ be a simply connected  compact Lie group
and let $H\subset G$ be a Lie subgroup. Let $\pi^i$ be an irreducible
representations of $LG$ with positive energy at level $k$ as in the 
introduction. 
Suppose when restricting to $LH$, $H^i$ decomposes as:
$$          
H^i = \sum_\alpha H_{i,\alpha} \otimes H_\alpha 
,$$ and  $\pi_\alpha$ are irreducible representations of $LH$ on
Hilbert space $H_\alpha$.  The set of $(i,\alpha)$ which appears in
the above decompositions will be denoted by $exp$. \par
Let ${\Cal A_{G/H}}$ 
be the Vacuum Sector of the coset $G/H$ as defined on Page 5 of
[X4]. Let us recall how ${\Cal A_{G/H}}$ is defined. Let $\pi^0$ be
the vacuum representation of $LG$ on Hilbert space $H^0$. Then
$ A_{G/H}(I)$ is naturally isomorphic to
$\pi^0 (L_IH)'\cap \pi^0 (L_IG)''$, where for any algebra $A\subset
B(H^0)$, $A'$ is defined to be its commutant. The coset 
$H\subset G_k$ is called {\it cofinite} if the inclusion
$$
(\pi^0 (L_IH)'')\vee(\pi^0 (L_IH)'\cap \pi^0 (L_IG)'') \subset
\pi^0 (L_IG)''
$$ has finite index and the square root of its minimal index will 
be denoted by $d(G/H)$ (cf. \S3 of [X1]). All the cosets 
considered in this paper are cofinite.\par 
The decompositions above naturally give rise to a class of
covariant representations of ${\Cal A_{G/H}}$, 
denoted by $\pi_{i,\alpha}$ or simply
$(i,\alpha)$ on Hilbert space $H_{i,\alpha}$. 
Note that on $H_{i,\alpha}$, the rotation group acts as 
$e^{i\theta}\rightarrow e^{ i \theta(\Delta_i-\Delta_\alpha)}$
by the nature of coset construction (cf. [GKO]), and hence for any
irreducible sector of $(i,\alpha)$, its univalence is
$$
e^{2\pi i (\Delta_i-\Delta_\alpha)} = \omega (i) \omega (\alpha)^{-1},
$$ hence the sector $(i,\alpha)$ has a uniform univalence with
univalence  $\omega (i) \omega (\alpha)^{-1}$ (cf. definition in \S1.7 
before Cor. 1.7.3). \par

In [X4], certain rationality
results (cf. Th. 4.2) are proved for a class of coset $H\subset G$.
A stronger rationality condition, $\mu$-rational is defined in \S1.3.
We will give a formula of  $\mu$-index for the 
coset in (1.6.2). 
We will see in the following $\mu$-index  is the 
square of the rank of certain unitary modular category. \par
We have:
\proclaim{Lemma 1.6}
Suppose $H\subset G_k$ is cofinite and the conformal precosheaf
  associated with
$H$ and $G$ are  $\mu$-rational. Then the coset conformal precosheaf has
finite $\mu$-index. 
Moreover its $\mu$-index (denoted by $\mu_{G/H}$) is given by
$$
\mu_{G/H} = \frac{d(G/H)^4 \mu_G}{\mu_H}, \tag 1.6.1
$$ where $d(G/H)$ is the cofinite statistical dimension defined
above, and $\mu_G ,\mu_H$ are $\mu$-index of the conformal precosheaf
associated with
$G$ and $H$ respectively.   
\endproclaim
This lemma is proved in [X3].\par
Throughout this section and the rest of the paper, we assume that the 
vacuun sectors associated with $H,G$ are  
$\mu$-rational , and $H\subset G$ is cofinite
(cf. \S 3 of [X1]). By the above lemma and Cor. 9 of [KLM] [KLM], 
the set of irreducible
superselection sectors of the coset is finite, and by Cor. 3.2 of [X3], 
every such irreducible sector appears as an irreducible subsector of some
$(j,\beta) \in exp$.  Moreover  the
$Y$-matrix as defined in \S2.1 is non-degenerate  
by Prop. 2.4 of [X3]. We will see in the
next section that this is the non-degeneracy condition of certain
unitary modular category. \par
Let us give a formula for 
$\mu_{G/H}$ as promised. It follows by Prop. 2.4 of [X3] that
$$
\mu_{G/H}= \sum_x d(x)^2, \tag 1.6.2
$$  where the summation is over all the different 
irreducible sectors $x$
which appears as irreducible sectors of some $(j,\beta) \in exp$.\par

\subheading{ 1.7 Unitary modular category from cosets}
Fix a proper interval $I_0\subset S^1$,
and let ${\Cal A}$ be the conformal precosheaf on $S^1$. 
Define $M:={\Cal A}(I_0)$.  
Suppose   a set of covariant 
representations  $\rho_i, i\in I$ of ${\Cal A}$ , localized on $I_0$ are given 
($\rho_i$ can be 
thought as an element in $End(M)$). 
We will assume that there is  an involution $i\rightarrow \bar
i$ on $I$  and there is a distinguished element $0\in I$ so that 
$$
[\rho_{\bar i}] = [\bar \rho_i], 
\rho_0 = id, \bar 0 =0. \tag 1.7.0
$$ 

Recall  from \S1.2.1 that 
for each $\rho_i\in \text{\rm End}$ and its conjugate
$\rho_{\bar i}$
there exists $R_{\rho_i} \in 
\text{\rm Hom}(id ,  \rho_{\bar i} \rho_i )$ and $ \bar R_{\rho_i} \in
\text{\rm Hom}(id , \rho_i \rho_{\bar i} )$ such that
$$
\bar R_{\rho_i}^* \rho_i(R_{\rho_i}) = 1,  R_{\rho_i}^*  
\rho_{\bar i}(\bar R_\rho) = 1
$$ and $||R_{\rho_i}||= ||\bar R_{\rho_i}||=\sqrt {d_{\rho_i}}$.  
The choice of the pair $(R_{\rho_i}, \bar R_{\rho_i})$ is not 
unique even up to
a phase unless 
$\rho_i$ is irreducible, but
 we fix a choice for each $i\in I$. We will see later 
that the invariants we are interested in are independent of 
such choices (cf. (4) of Lemma 1.7.4.).  We will also choose the
pairs  $(R_{\rho_i}, \bar R_{\rho_i})_{i\in I}$ so that 
$$
R_{\rho_{\bar i}} = \bar R_{\rho_i}, \overline{ R_{\rho_{\bar i}}} = R_{\rho_i}
.$$ Recall the minimal
left inverse $\phi_{\rho_i}$ of $\rho_i$ is defined by
$$
\phi_{\rho_i} (m) = R_{\rho_i}^* \rho_{\bar i}(m) R_{\rho_i}
.$$\par  
Let us start to define a category denoted by  $C({\Cal A})$ from $\rho_i, i\in
I$
subject to (1.7.0).  
The objects of  $C({\Cal A})$ are defined to be any finite compositions of
$\rho_i, i\in I$, considered as elements in $End(M)$. The homorphisms are
defined as intertwinners as in \S1.2 . 
The tensor product of two objects $V, W \in C({\Cal A}) $ are defined by
$$
V\otimes W:= VW
$$ where $VW(m):= V(W(m)), \forall m\in M$, 
i.e.,  $VW$ stands for the composition of endomorphism V with 
endomorphism W. We have omitted the usual $\cdot$ when no confusion
arises. Now let $f\in Hom(V_1,V_2), g\in Hom(W_1, W_2)$, then
$$
f\otimes g:= f V_1(g).
$$
Note this definition makes sense since $V_1\in End (M), g\in M$.
It is then easy to see that
$C({\Cal A})$ is a strict abelian monoidal category. The braidings 
$c_{V,W}, \forall \lambda,\mu\in C({\Cal A})$  are defined 
to be $\sigma_{VW}$  as in the
beginning of \S1.4.  Note that   (1.1.1), (1.1.2)
follows immediately from Prop. 1.4.2 . \par
Next we define duality. For each $\rho_i$, we define
$$
b_{ \rho_i} := \bar R_{\rho_i} , 
d_{ \rho_i} :=  R_{\rho_i}^*.
$$ To save some writing we shall use $i$ to denote  
$\rho_i$ in the following. So $\rho_{i_1}\rho_{i_2}...\rho_{i_n}$
can be writen as $i_1 i_2...i_n$. For any objects $i_1 i_2...i_n$ we define
$$
(i_1 i_2...i_n)^* := \bar i_n... \bar i_2 \bar i_1,
$$ and define
$$
b_{i_1 i_2...i_n} :=i_1i_2...i_{n-1} (b_{i_n})...i_1(b_{i_2}) b_{i_1},
$$ and
$$
d_{i_1 i_2...i_n} := 
\bar i_n \bar i_{n-1}...\bar i_{2} (d_{i_1}) \bar i_n \bar i_{n-1}...
\bar i_{3} (d_{i_2})... d_{i_n}.
$$

\par
For each object $i_1 i_2...i_n \in C({\Cal A})$ and
$f\in End(i_1 i_2...i_n)$ we define unnormalized trace of $f$ 
as in \S1.2:
$$
tr(f):= d(i_n) ... d(i_1)
\phi_n (\phi_{n-1}(...(\phi_1(f))...)) \tag 1.7.1
$$ where $ \phi_k$ is the unique minimal left inverse of $\rho_{i_k},
k=1,...,n$, and
$d(i_k)$ is the statistical dimension of 
$\rho_{i_k},
k=1,...,n$.   Note that (1.7.1) can be also written as
$$
tr(f):=  d_{i_n} ... d_{i_1}\phi_{i_1...i_n} (f)
$$ where $\phi_{i_1...i_n}$ is the unique minimal left inverse of
$i_1...i_n$ since by \S1.2 $\phi_{i_1...i_n} = \phi_n \phi_{n-1}...\phi_1$.
\par
Next let us define twist: 
$$
\theta_{i_1...i_n}:= d_{i_n} ... d_{i_1}\phi_{i_1...i_n} 
(c_{i_1...i_n,
i_1...i_n}).
$$ 
Denote by $\tilde \theta_V :=  V(b_V^*)c_{V,V} V(b_V)$.\par 
Finally note that for any  morphisms $f: U\rightarrow V$ in $C({\Cal A})$,
$f$ is an element of von Neumann algebar $M$.  The conjugation of $f$
is defined to be $f^*: V\rightarrow U$.\par
Let us first note the following lemma which follows immediately
from the definitions:
\proclaim{Lemma 1.7.1}
For the category $C({\Cal A})$ defined above, we have:\par
(1). (1.1.1) and (1.1.2) are true;  \par
(2). (1.1.8) is true, i.e.,
$$
\theta_V^* = \theta_V^{-1};
$$\par
(3). (1.1.5) is true.
\endproclaim
\demo{Proof}
Ad (1): (1.1), (1.2) follows immediately from Prop. 1.4.2 of \S1.4; 
(2)  follows from
Page 244 of [L1].  (3) follows from definitions and
(1.2.1).
\enddemo
\qed 
\par
\proclaim{Lemma 1.7.2}
(1). (1.1.3) is true (note by our definition 
$\theta_V\otimes \theta_W:= V(\theta_W)$)
, i.e.,
$$
\theta_{VW} =  c_{W,V}  c_{V,W}\theta_V V(\theta_W); 
$$
(2).
$$
\tilde \theta_{VW} =   c_{W,V} c_{V,W}
\tilde \theta_V V(\tilde \theta_W); 
$$
(3). $\theta_{V}=\tilde \theta_{V}$ if and only if
$ \theta_{\rho_i} = \tilde \theta_{\rho_i}$;\par 
(4). $\theta_{V}=\tilde \theta_{V}$ if and only if
(1.1.6) is true;  \par
(5).  $ \theta_{\rho_i}= \tilde \theta_{\rho_i} = c_i id_{\rho_i}$ 
if and only if for any irreducible subsector $x$ of $\rho_i$, 
the univalence of $x$ is equal to $c_i$.  
\endproclaim
\demo{Proof}
Ad (1):  By Prop. 1.4.2
$ c_{VW, VW} = V(c_{V,W}) V^2(c_{W,W}) c_{V,V} V(c_{W,V})$,
and so
$$
d(V)\phi_V (c_{VW, VW}) = c_{V,W} V(c_{W,W}) \theta_V c_{W,V}.
$$ By using Braiding-Fusion-Equation of  Prop. 1.4.2  and the fact
that $c_{1,V}= id$ we have
$$
\theta_V c_{W,V} =  c_{W,V} W(\theta_V),
\theta_W c_{V,W} =  c_{V,W} V(\theta_W),
$$ and by Yang-Baxter equation
$$
c_{V,W} V(c_{W,W})c_{W,V} = W(c_{W,V}) c_{W,W}  W(c_{V,W})
$$ and so
$$
\align
\theta_{VW} & = d(V) d(W) \phi_W \phi_V ( c_{VW, VW}) \\
&= d(W) \phi_W ( c_{V,W} V(c_{W,W}) \theta_V c_{W,V}) \\
& =  d(W)
\phi_W (  c_{V,W} V(c_{W,W})  c_{W,V}  W(\theta_V)) \\
&= d(W)\phi_W (W(c_{W,V}) c_{W,W}  W(c_{V,W})   W(\theta_V)) \\
&= c_{W,V} \theta_W  c_{V,W} \theta_V \\
&= c_{W,V} c_{V,W} V(\theta_W) \theta_V \\
&= c_{W,V} c_{V,W}  \theta_V V(\theta_W)
\endalign
$$ which is (1). (2) is proved in similar way.
(3) follows from (1), (2) and the construction of $C({\Cal A})$. 
\S 1.2.
Ad (4): By BFE of Prop. 1.4.2 and (1.1.5)
$$
\align
\tilde \theta_V & = V(b_V^*)V(c_{V,V^*}^{-1}) b_V \\
&= V(b_V^*)V(c_{V,V^*}^{-1}) b_V VV^*(b_V V(d_V)) \\
&= V(b_V^*)V(c_{V,V^*}^{-1})VV^*(b_V V(d_V))  b_V \\
&=  V(d_V)V(b_V^*) V^2 (c_{V^*,V^*})V(b_V) b_V  \\
&=   V(d_V)V(\theta_{V^*})b_V.
\endalign
$$ Hence  
$\tilde \theta_V = \theta_V$ is equivalent to
$$
V(d_V)V(\theta_{V^*})b_V =\theta_V,
$$ and multiplied both sides on the right by $b_V$,
using
$$
V(d_V)V(\theta_{V^*})b_V b_V=  V(d_V)V(\theta_{V^*}) VV^*(b_V)  b_V
= V(d_V)V V^*(b_V) V(\theta_{V^*}) b_V
$$ and (1.1.5)  we get
$$
V(\theta_{V^*}) b_V = \theta_V b_V,
$$ which is (1.1.6).  The other direction is similar. \par
Ad (5): Suppose $ \theta_{\rho_i}= \tilde \theta_{\rho_i} = c_i
id_{\rho_i}$. By Cor. 1.3.2 the univalence of $x$ must be $c_i$. 
Now suppose the univalence of any irreducible subsector $x$ of 
$\rho_i$ is the same constant $c_i$. By Cor. 1.3.2, $
\theta_{\rho_i}= c_i
id_{\rho_i}$, so we just have to show that
$ \tilde \theta_{\rho_i}=\theta_{\rho_i} 
$ to finish the proof. 
By (4) this is equivalent to show that 
$$
\theta_{\rho_{\bar i}} = \theta_{\rho_i}=  c_i,
$$
and by    Cor. 1.3.2 again it is
sufficient to show that for any   irreducible subsector $y$ of 
$[\rho_{\bar i}]= [\overline{\rho_i}] $, the univalence of $y$ is $c_i$, 
which follows from  the third equation in (1.5.1). 
\enddemo
\qed
\par
For convenience we shall say that $\rho_i$ has a {\it uniform univalence}
  if for any irreducible subsector $x$ of $\rho_i$, 
the univalence of $x$ is equal to a constant $c_i$ which depends
only on $\rho_i$, and in
this case  we will define $\exp(2\pi i \Delta_{\rho_i}) := c_i$. Note
that any irreducible endomorphism has a  uniform univalence. \par
\proclaim{Corollary 1.7.3 }
(1):
If each $\rho_i, i\in I$ has  a uniform valence, then the 
trace defined in (1.7.1) agrees with (1.1.7), and the
category
$C({\Cal A})$ is an abelian unitary ribbon category; \par
(2): If the set $\{ \rho_i, i\in I \}$ contains a finite family
$\{ V_j, j\in J \}$ of simple objects satisfying the four 
defining axioms of modular category in \S1.1, then 
$C({\Cal A})$ is a unitary  modular category.
\endproclaim
\demo{Proof}
Ad (1):
Let us first prove the second equation in (1.1.9), i.e., 
$$
b_V^* = d_V c_{V,V^*} \theta_V. 
$$   By (1.1.5) it is sufficient to show that
$$
d_V c_{V,V^*} \theta_V V(d_V^*) =1.
$$
Note $ d_V c_{V,V^*} \theta_V V(d_V^*) = d_V c_{V,V^*}
V(d_V^*)\theta_V$,
and by  Prop. 1.4.2  we have
$$
 d_V c_{V,V^*}
V(d_V^*)= \theta_V^*,
$$ and so
$$
d_V c_{V,V^*} (\theta_V) V(d_V^*) =\theta_V^* \theta_V=1.
$$  \par
To prove the third equation 
in (1.1.9) 
$d_V^*=V^*(\theta_V^{-1}) c_{V^*,V}^{-1} b_V,$
by (1.1.5) it is
sufficient to show that
$$
b_V^*V(V^*(\theta_V^{-1}) c_{V^*,V}^{-1} b_V)=1,
$$ but
$$
b_V^*V(V^*(\theta_V^{-1}) c_{V^*,V}^{-1} b_V) =\theta_V^{-1}
b_V^* V(c_{V^*,V}^{-1}) V(b_V),
$$
$$
b_V^* V(c_{V^*,V}^{-1})V(b_V) = V(b_V^*)c_{V,V} V(b_V),
$$
hence it is enough to show that
$$
 V(b_V^*)c_{V,V} V(b_V) = \theta_V,
$$ i.e., $\tilde  \theta_V = \theta_V$, which is implied by 
(3), (5) of lemma 1.7.2 and our assumptions.\par
Now let us show that (1.7.1) agrees with (1.1.7). 
 We need to check for any
$f:V\rightarrow V$,
$$
d_V c_{V,V^*} \theta_V f b_V = d_V V^*(f) d_V^*
,$$ and by (1.1.9) just proved, this is equivalent to
$$
b_V^* f b_V = d_V V^*(f) d_V^*.
$$  Note by BFE and the last equation of (1.1.9) we have:
$$
b_V^* f b_V =  b_V^* c_{V^*,V} V^*(f) c_{V^*,V}^{-1} b_V
=d_V V^*(\theta_V^{-1} f \theta_V)) d_V^*,
$$ but the last expression is equal to   
$ d_V V^*(f) d_V^*$ by the tracial property lemma 1.2.1. \par
Finally we check (1.1.4), i.e.,
$$
d_V V^*(c_{V,V}) d_V^* f = f d_U V^*(c_{U,U}) d_U^*
$$ for any $f: U\rightarrow V, U,V \in C({\Cal A})$. By using the naturality
of the braiding (1.1.1), (1.1.2), $ d_V c_{V,V} d_V^* f = d_V V^*(f)
V^*(c_{V,U})d_V^* $, $f d_U c_{U,U} d_U^* =  d_U U^*(c_{V,U}) U^*(f)d_U^*$.
Let $f_1 :V^*\rightarrow U^*$ be a homorphism such that
$$
U^*(f) d_U^*= f_1 d_V^*.
$$ Note by (1.1.4) the above defines $f_1$ uniquely. In fact
$$
f_1 = U^*(b_U^*) U^*(f) d_U^*.
$$ Using $f_1$ we have
$$
d_U U^*(c_{V,U}) U^*(f)d_U^*  =d_U U^*(c_{V,U})  f_1 d_V^*
= d_U f_1  V^*(c_{V,U})d_V^* ,
$$ hence to prove (1.1.4) we just have to show  
$$
d_V V^*(f) = d_U f_1,
$$ or equivalently by using (1.1.5)
$$
U^*(b_U^*) U^*(f) d_U^* = d_V V^*(f) V^*(b_U).
$$  Denote by
$$
A:= U^*(b_V^*) U^*(f) d_U^*, B:=  d_V V^*(f) V^*(b_U).
$$ Using the fact that (1.7.1) agrees with (1.1.7) which 
is proved above, we have:
$$
tr (AA^*) = tr (ff^*) = tr (B B^*) ,
$$ and $tr(AB^*) = tr(BA^*) = tr (f^*f)$. But by 
lemma 1.2.1 we have $tr (ff^*)=  tr (f^*f)$, so
$$
tr((A- B)(A-B)^*) = 0,
$$ and since $tr$ is faithful, we have 
$$
A=B,
$$ which completes the proof of (1.1.4). \par
Finally (1.10) follows from the definition. \par
So we have verified (1.1.1) to (1.1.6) and (1.1.9), (1.1.10)
for $C({\Cal A})$, which shows that $C({\Cal A})$ is an abelian
unitary ribbon category. \par
(2) follows from (1) and definitions.
\enddemo \hfill \qed
\par
Now suppose ${\Cal A}$ is ${\Cal A}_{G/H}$, the conformal precosheaf
associated the coset $H\subset G_k$ as defined in \S1.6. We will denote
$C({\Cal A}_{G/H})$ simply by $C(G/H)$ in the following.  
Let us  make a specific choice of $\rho_i, i\in I$. We will choose
 $\rho_i,i\in I$ to be all the irreducible sectors of any 
$(j,\beta)\in exp$ and also all $(j,\beta)\in exp$. 
Note by (2.3.1) if $(j,\beta)\in exp$, then $(\bar j,\bar \beta) \in
exp$ and $[\overline{(j,\beta)}]= [(\bar j,\bar \beta)]$. Hence such a set
verifies (1.7.0).  
Notice that
for any $(j,\beta)\in exp$, it has a uniform univalence which 
is equal to $\exp(2\pi i (\Delta_j - \Delta_\beta))$ by \S1.6. 
The reader may wonder why we include these (possibly reducible) 
sectors as our building blocks. The reason is that it is convenient
to include these objects for the computations of invariants,
cf. \S3.4 and \S3.6 below.
\par
\proclaim{Theorem A}
If a coset $H\subset G_k$  verifies the conditions of
\S 1.6, i.e., $H\subset G_k$ is cofinite, and the conformal
precosheaf associated with $H$ and $G$ are $\mu$-rational,
then  the category $C(G/H)$ constructed from the 
choices of  $\rho_i,i\in I$ above is a unitary modular 
category.
\endproclaim
\demo{Proof}
The  category $C(G/H)$ is a unitary ribbon category
by Cor. 1.7.3, so we just need to check axioms (1)-(4) in \S1.1. 
Denote by $x, x\in X$ the finite set of simple objects 
in $C(G/H)$ obtained as 
irreducible subsectors of all $(j,\beta) \in exp$. Denote the 
identity sector in $X$ by $0$ as in axiom (1). This set obviously
satisfies  axiom (1), and it satisfies axioms (2) and (3) since all
irreducible sectors of the conformal precosheaf ${\Cal A}_{G/H}$ 
appear in $X$ by Cor. 3.2 of [X3]. To prove axiom (4), note that
by definition (1.5.0) 
$$
Y_{xy}= \overline{S^{Tu}_{xy}},
$$ and since the matrix $(Y_{xy})$ is nondegenerate by \S1.6, 
it follows that axiom (4) holds.
\enddemo
\qed
\par
We conjecture that all the cosets $H\subset G_k$ satisfy the 
conditions of Theorem A. Infinite series of cosets which satisfy 
the conditions of Theorem A will be discussed in \S3. For more
such infinite series, see \S3.1 of [X3] and \S3.2 of [X2].\par
Note that if we take $H_1={e}$ to be a trivial group of $G$, then
$H_1\subset G_k$ is automatically cofinite, 
in fact the  conformal
precosheaf associated with $H_1=\{ e \}\subset G$ is simply 
the conformal
precosheaf associated with $G$. 
So under the conditions of Th. A, we have 
a  unitary modular 
category associated  $G$, denoted by $C(G)$. Similarly 
we have  unitary modular 
category associated  $H$, denoted by $C(H)$.
Much of the rest of this paper focus on the relations between
$C(G/H)$ and the categories $C(G)$ and $C(H)$.\par 
Let us  calculate $D_{C(G/H)}, k_{C(G/H)}^3$, as defined for any modular
category in (1.11) and (1.12). 
We will write $D_{C(G/H)}, k_{C(G/H)}^3$ simply as $D_{G/H}, k_{G/H}^3$ 
in the following. 
By comparing (1.6.2) and (1.1.11) we have
$$
D_{G/H}^2 = \mu_{G/H}
$$ and similarly 
$D_{G}^2 = \mu_{G}, D_{H}^2 = \mu_{H}$, and by 
lemma 2.2 we have the following formula for $D_{G/H}$:
$$
D_{G/H} = \frac{D_G d(G/H)^2}{D_H} \tag 1.7.2
$$
By the remark after the proof of Prop. 3.1 in
[X3] we have 
$$
k_{G/H}^3= \frac{k_G^3}{k_H^3} \tag 1.7.3
$$
where $ k_G^3$ (resp. $k_H^3$) is $k_{C(G)}^3$ 
(resp. $k_{C(H)}^3$)
as defined in (1.12) for modular category $C(G)$ (resp. $C(H)$). \par
To calculate $D_{G}, k_G^3$ (or $D_{H},k_H^3$), 
let us assume that $G=SU(N)$ (or $H=SU(k)$) even though
the following argument works for any $G$ which has the property that
its genus 0 modular matrix agrees with genus 1 modular matrix
(cf. \S1.5). 
By comparing (1.11), (1.5.0) and (1.5.9) we immediately have
the following formula:
$$
D_{G}= \frac{1}{S_{11}} \tag 1.7.4
$$
and  comparing (1.12), (1.5.3) and (1.5.12) we 
have the following:
$$
k_G^3 =\exp(\frac{\pi i C_G}{4}) \tag 1.7.5
$$
where $C_G$ the central charge as defined in (1.5.11).\par
Define the central charge $C_{G/H}$ of the coset by 
$$C_{G/H}:= C_G-C_H \tag 1.7.6
$$
It follows from (1.7.3) and (1.7.5) that
$$
k_{G/H}^3= \exp(\frac{\pi i C_{G/H}}{4}) \tag 1.7.7
$$ if $H, G$ are of type $A$.\par
The following two lemmas will be useful in the calculation
of framed oriented link invariants from $C(G/H)$. 
Given a framed oriented link $L$ with
$n$ components, colored by objects $i_1,...,i_n \in C(G/H)$. 
The corresponding link invariants as given by Th. I.2.5 of [Tu]
will be denoted by $L(i_1,...,i_n)$. As in \S1.1,  suppose
$L(i_1,...,i_n)$ is represented by the closure of a tangle $T_L \in
End (I_1 I_2...I_n) $, 
where $I_m:= (i_m^{\epsilon_m})^{k_m},
k_m\in {\Bbb N}, m=1,...,n$ and $\epsilon_m=1$ (resp. $\epsilon_m=-1$)
if the $m$-th component of $L$ is oriented anti-clockwise ( resp.
clockwise ).
 Let $j_k(l), 1\leq l \leq f(k)$ be a set
of susectors of $i_k$, $T_k(l) \in Hom (j_k(l), i_k)$ such that
$$
T_k(l)^* T_k(l')= \delta_{l,l'}, 
\sum_{l}T_k(l) T_k(l)^* =id_{i_k}.
$$  Let $p_k(l) :=T_k(l)T_k(l)^*$.  Then we have the following :
\proclaim{Lemma 1.7.4}
With the notations introduced above, we have \par
(1)
$$
L(j_1(l_1),... j_n(l_n)) = tr ( p_1(l_1) I_1 (p_2(l_2))
...  I_1 I_2 I_{n-1}(p_n(l_n))
T_L);
$$   
(2) 
$$
\sum_{j_1(l_1),... j_n(l_n)} L(j_1(l_1),... j_n(l_n)) =
L(i_1,...,i_n);
$$
(3)
If $j_1(l_1),... j_n(l_n) \in C(G/H)$ such that 
$[j_m(l_m)] =[i_m], m=1,2,...,n$, then
$$
L(j_1(l_1),... j_n(l_n)) = L(i_1,...,i_n); 
$$
(4)
$L(i_1,...,i_n)$ are independent of the choices of $(R_{\rho_i}, 
\bar R_{\rho_i})_{i\in I}$ in the construction of 
$C(G/H)$.
\endproclaim
\demo{Proof}
Ad (1):  The idea is to
insert $ T_k(l_k)^*  T_k(l_k)=id_{j_k(l_k)}$ 
on the $k$-th string of $T_L$ colored
by $j_k(l_k)$. When we pull the coupon represented by $  T_k(l_k)^*$ 
clockwise around the $k$-th component of $L$ and approaches 
$ T_k(l_k)$ from below, we obtain the righthand side of 
(1) by the property of ribbon category. 
(2) follows from (1). (3) follows from (1) since 
$p_m(l_m) = id_{i_m}$. (4) follows since the tangle $T_L$ consists
of braidings and twistings only, and the closure of  $T_L$ are given
by minimal conditional expectations which are independent of the
choices of  $(R_{\rho_i}, 
\bar R_{\rho_i})_{i\in I}$ (cf. \S1.2).
\enddemo
\qed
\par
Now suppose given $C(G/H)$. Let $v$ be an object of $C(G/H)$. We say
that $v$ is a {\it generator} of $C(G/H)$ if there exists a finite
subset $K\subset {\Bbb N}$ such that 
any irreducible object of  $C(G/H)$, as a 
sector  appears as 
a subsector of $[v^m], m\in K$.    
Notice that there is a natural inclusion $End(v^k)\subset End(v^l)$
if $l\geq k$. 
Now suppose $C(G/H)$ (resp.  $C(G_1/H_1)$)
are given with generators $v$ (resp. $v_1$). We say that $C(G/H)$
is {\it compatible}  with  $C(G_1/H_1)$ if:\par
(c1): $v$ and  $v_1$ are of the uniform univalence; \par 
(c2): there exists a 
sequence of trace-preserving 
$*$-isomorphisms $\psi_k: End(v^k)\rightarrow  End(v_1^k)
, \forall k\geq 1$ such that
$$
\align
\psi_k (f) = \psi_l (f)   \ \ \text{\rm if} l\geq k, 
\psi_{k+1} (v(f)) & = v_1 (\psi_{k}(f)) \ \forall f\in End(v^k) ;
\\
\psi_2(c_{v,v}) & = \psi_2(c_{v_1,v_1}); 
\endalign
$$ 
(c3): there exists a  finite
subset $K\subset {\Bbb N}$ such that the irreducible objects of 
$C(G/H)$ (resp.$C(G_1/H_1)$), as sectors are in one-to-one correspondence
with the irreducible subsectors of $ [v^m], m\in K$ 
(resp. $ [v_1^m], m\in K$), i.e., 
$$
\langle[v^m]  , [v^{m'}]\rangle =0, \ \text{\rm if }m\neq m', m,m'\in K,
$$
$$
\langle[v_1^m]  , [v_1^{m'}]\rangle =0, \ \text{\rm if }m\neq m', m,m'\in K,
$$
and the set of 
irreducible subsectors of $ [v^m], m\in K$ 
(resp. $ [v_1^m], m\in K$) are the same as  the set of irreducible objects of 
$C(G/H)$ (resp.$C(G_1/H_1)$) as sectors.  \par
Note that by definition the univalence of $v$ is $\frac{tr (c_{v,v})}
{tr(id_v)} $, and by 
(c2) 
$$ 
tr (c_{v,v}) =  tr (c_{v_1,v_1}),tr(id_v)= tr(id_{v_1}), 
$$ so if 
$C(G/H)$ and $C(G_1/H_1)$ are compatible then 
$v$ and $v_1$ have the same univalence. We can therefore add to
(c1) that $v$ and $v_1$ have the same univalence. \par 
Let $p_i\in End(v^{m_i}),  i\in I, m_i\in K$ 
be the set of minimal projections corresponding to the set of  
irreducible objects (denoted by $c(p_i)$) of   $C(G/H)$. 
Let $q_i := \psi_{m_i} (p_i), m_i\in K$ and  $c(q_i)$ be the simple objects
of $C(G_1/H_1)$. 

\proclaim{Lemma 1.7.5}
If $C(G/H)$ and  $C(G_1/H_1)$ are compatible, then with 
the notations above we have:\par
(1) $d (c(p_i)) = d (c(q_i)), \omega(c(p_i)) = \omega(c(q_i));$ \par
(2) $D_{G/H} = D_{G_1/H_1}; k_{G/H} = k_{G_1/H_1};$ \par
(3)
$\tau_{G/H}(M_L) = \tau_{G_1/H_1}(M_L).$
\endproclaim
\demo{Proof}
Ad (1):
Note $ d (c(p_i))= tr(c(p_i)),  d (c(q_i))= tr(c(q_i)) = tr (\psi_{m_i}
(c(p_i)))
$, so $d (c(p_i)) = d (c(q_i))$ since $\psi_{m_i}$ is trace-preserving. 
By definition and naturality 
of the twist (cf. (1.1.4)) 
$ d (c(p_i)) \omega(c(p_i)) = tr(\theta_{c(p_i)}) =
tr( \theta_{v^{m_i}} p_i)$ and similarly $d (c(q_i)) \omega(c(q_i)) = 
tr( \theta_{v_1^{m_i}} q_i)$, so we just have to show
$\psi_m(\theta_{v^{m_i}})= \theta_{v_1^{m_i}}$. By (1.1.3) $\theta_{v^{m_i}}$
(resp. $\theta_{v_1^{m_i}}$) is a product of $v^j(c_{v,v}),j\in {\Bbb N} $ 
(resp. 
 $v^j(c_{v_1,v_1}), j\in {\Bbb N}$), their inverses, and the phase
$\omega(v)^{m_i}$ (resp. $\omega(v_1)^{m_i}$), and it follows from (c1), (c2)
and the remark after them
that $\psi_m(\theta_{v^{m_i}})= \theta_{v_1^{m_i}}$ and (1) is proved. \par 
(2) follows from (1) and definitions. \par
Ad (3): 
Let $L$ be any framed oriented link with $n$ 
component in 
$S^3$, since $\tau_{G/H}(M_L)$ (resp. $\tau_{G_1/H_1}(M_L)$)
is independent of the orientations of $L$, we can assume that
$L$ is represented by the closure of a tangle $T_L$ with all
the bands oriented downward.   By (2), to prove 
(3) it is enough  to show that for such $L$
$$
L(c(p_1),...,c(p_n)) = L(c(q_1),...,c(q_n))   
. 
$$  By (1) of lemma 1.7.4, with $i_1= v^{m_1},..., i_n=v^{m_n}$, 
$$
L(c(p_1),...,c(p_n)) = tr(p_1 v^{k_1m_1} (p_2) ..., v^{k_1m_1+...+k_{n-1}
m_{n-1}} 
(p_n) T_{L^c} ) 
,$$  where $L^c$ is the 
cabling of $L$ (cf. [W2]) with cabling vector
$c=(m_1,m_2,...m_n)$.  By  using the compatibility conditions  
$$
\align
& tr(p_1 v^{k_1m_1} (p_2) ... v^{k_1m_1+...+k_{n-1}m_{n-1}} (p_n) T_L )  \\
&= tr(\psi_{k_1m_1+...+k_{n}m_n}(p_1 
v^{k_1m_1} (p_2) ... v^{k_1m_1+...+k_{n-1}m_{n-1}}(p_n)) \times \\
&\psi_{k_1m_1+...+k_{n-1}m_{n-1}} (T_{L^c} )),
\endalign
$$ note 
$$
\align
& \psi_{k_1m_1+...+k_{n}m_{n}}
(v^{k_1m_1+...+k_{s-1}m_{s-1}} (p_s)) \\ 
&= \psi_{k_1m_1+...+k_{s-1}m_{s-1}+m_s} 
( v^{k_1m_1+...+k_{s-1}m_{s-1}} (p_s)) \\
&=v_1(\psi_{k_1m_1+...+k_{s-1}m_{s-1}+m_s-1} 
( v^{k_1m_1+...+k_{s-1}m_{s-1}-1} (p_s)))
\\
&= ... = v_1^{k_1m_1+...+k_{s-1}m_{s-1}}(\psi_{m_s} (p_s)) 
\endalign
$$ by repeatedly using
the compatibility conditions, and so we have
$$
L(c(p_1),...,c(p_n)) = tr(q_1  v^{k_1m_1} (q_2) ... 
v^{k_1m_1+...+k_{n}m_{n}} 
(q_n)
\psi_{k_1m_1+...+k_{n}m_{n}}(T_{L^c}) ). 
$$ To finish the proof we just have to show that 
$$
\psi_{k_1m_1+...+k_{n}m_{n}}(T_{L^c} (v,v,...v) ) = T_{L^c} (v_1,v_1,...v_1) 
.$$  Note that  $T_{L^c}$   consists of only 
twistings and braidings. So $T_{L^c} (v,v,...v)$ up to a phase 
factor determined by the univalence of $v$, is equal to products of
$v^j(c_{v,v}), j\in {\Bbb N}$ and their inverses.  
It follows from (c1) (c2) and the remark
after them that
$$
\psi_{nm}(T_{L^c} (v,v,...v) ) = T_{L^c} (v_1,v_1,...v_1), 
$$ and the proof is complete.
\enddemo
\qed

\heading \S2. Factorization of link invariants \endheading
The goal of this section is to calculate link invariants colored by
$(j,\beta)\in exp$ under the condition 
that $d((j,\beta)) = d(j) d(\beta)$. 
In [X1], [X2] the sectors of $[(j,\beta)]$ 
are studied by using the idea of [X4].
 But to
calculate link invariants one needs to analyze the 
braiding properties of the system
in [X4]. This is essentially done in [BE3], and ideas are also contained in
the proof of lemma 3.3 on Page 384 of [X4]. We will use the notations of
 [X1]and  [X4], and set up a dictionary in the following between our
notations and that of [BE1] to [BE3] so the similarity of the argument
is clear. \par
\subheading{2.1 Preliminaries}
In this section we sketch some of the results of \cite{LR} which
will be used in this paper. For the details of all the proofs  , 
we  refer the reader to \cite{LR}.
We have changed some of the notations in \cite{LR} since they have
been used to denote different objects in this paper.\par
Let $N\subset M$ be an inclusion of type III von Neumann algebras
on a Hilbert space $H$. Let $\phi \in H$ be a joint cyclic and separating
vector (which always exists if $N$ is type III and $H$ is separable).
Let $j_N = Ad_{J_N}$ and  $j_M = Ad_{J_M}$ be the modular
conjugations w.r.t.$\phi$ and the respective algebra. Then
$$
\alpha = j_Nj_M |_M \in End(M)
$$
maps $M$ into a subalgebra of $N$. We call $\alpha$ the
canonical endomorphism associated with the subfactor, and denote by
$$
N_1 := j_Nj_M (M) \subset N, M\subset M_1:= j_M j_N(N)
$$
the canonical extension resp. restriction. $\phi$ is again
joint cyclic and separating for the new inclusions above, giving
rise to new canonical endomorphisms $\beta = j_{N_1} j_N \in End(N)$
and $\alpha_1 = j_M j_{M_1} \in End(M_1)$.\par
We have the following formula for canonical endomorphism:
(cf. Prop.2.9 of [LR]):
\proclaim{Proposition 2.1.1}    
Let $N\subset M$ be an inclusion of properly infinite factors,
$\epsilon: M\rightarrow N$ a faithful normal conditional
expectation, and $e_N \in N'$ the associated Jones projection.
The canonical endomorphism $\alpha: M\rightarrow N$ is given by
$$
\alpha = \Psi^{-1} \cdot \Phi
$$
where
$$
\Phi(m) = U m U^* (m\in M)
$$
is the isomorphism of $M$ into $Ne_N$ implemented by an isometry
$U\in M_1 = \langle M, e_N \rangle$ with $ U U^* = e_N$, and
$\Phi$ is the isomorphism of $N$ with $Ne_N$ given by
$$
\Phi(n)= ne_N (n\in N)
.$$                    
Every canonical endomorphism of $M$ into $N$ arises in this
way.
\endproclaim
{\bf Definition}: A net of von Neumann algebras over a
partially ordered set ${\Cal J}$ is an assignment ${\Cal M}:
i \rightarrow M_i$ of von Neumann algebras on a Hilbert space
to $i\in {\Cal J}$ which preserves the order relation, i.e.,
$M_i \subset M_k$ if $i\leq k$. A net of subfactors consists of
two nets ${\Cal N}$ and  ${\Cal M}$ such that for every $i\in {\Cal J}$,
$N_i \subset M_i$ is an inclusion of subfactors.
We simply write ${\Cal N} \subset {\Cal M}$. The net ${\Cal M}$ is called
standard if there is a vector $\Omega \in H$ which is cyclic and separating
for every $M_i$. The net of subfactors  ${\Cal N} \subset {\Cal M}$ is called
standard if   ${\Cal M}$ is standard and   ${\Cal N}$ is standard on a subspace
$H_0 \subset H$ with the same cyclic and separating vector $\Omega \in H_0$.
For a net of subfactors  ${\Cal N} \subset {\Cal M}$, let $\epsilon$ be
a consistent assignment $i \rightarrow \epsilon_i$ of normal
conditional expectations. Consistency means that       
$\epsilon_i = \epsilon_k|M_i$ whenever $i\leq k$. Then we call
$\epsilon$ a normal conditional expectation from $ {\Cal M}$ onto
 ${\Cal N}$. $\epsilon$ is called standard, if it preserves the vector
state $\omega = (\Omega, \cdot \Omega)$ \par

If the index set ${\Cal J}$ is directed, i.e., for $j,k \in {\Cal J}$
there is $m\in {\Cal J}$ with $j,k \leq m$, we associate with a net
${\Cal M}$ of von Neumann algebras the inductive limit $C^*$ algabra
$(\bigcup_{i\in {\Cal J}} M_i)^-$ and denote it by the same
symbol ${\Cal M}$. Then we have (cf. Cor.3.3 of [LR]):
\proclaim{Proposition 2.1.2}
Let  ${\Cal N} \subset {\Cal M}$ be a directed standard net of subfactors
(w.r.t. the vector $\Omega \in H$) over a directed set ${\Cal J}$,
and $\epsilon$ a standard conditional expectation. For every
$i \in {\Cal J}$ there is an endomorphism $\alpha^i$ of the
$C^*$ algebra $ {\Cal M}$ into  ${\Cal N}$ such that
$\alpha |M_j$ is a canonical endomorphism of $M_j$ into $N_j$
whenever $i\leq j$. Furthermore, $\alpha^i$ acts trivially on
$M_i' \cap {\Cal N}$. As $i\in  {\Cal J}$ varies to $k$,  
the corresponding endomorphisms $\alpha^i$ and  $\alpha^k$
are inner equivalent by a unitary in $N_l$ whenever
$i,k \leq l$
\endproclaim
If the index set ${\Cal J}$ is directed, by Cor.4.2 of [LR]
the index is constant in a directed standard net of subfactors
with a standard conditional expectation. We have (cf. Cor.4.3 of [LR]):
\proclaim{Proposition 2.1.3}
Let  ${\Cal N} \subset {\Cal M}$ be a directed standard net of subfactors
(w.r.t. the vector $\Omega \in H$) over a directed set ${\Cal J}$,
and $\epsilon$ a standard conditional expectation. If the index $d^2 =
Ind (\epsilon)$ is finite, then
 
 for any $i \in {\Cal J}$,
there is an isomorphic intertwiner
$v_1:id \rightarrow \alpha$ in $M_i$ which satisfies the following
identity with the isometric intertwiner $w_0: id \rightarrow \beta$
($\beta := \alpha|{\Cal N}$) in
$N_i$:
$$       w_0^* v_1 = d^{-1} id = w_0^* \alpha (v_1)
.$$
$\alpha$ is given on ${\Cal M}$ by the formula
$$
\alpha (m) = d^2 \epsilon (v_1 m v_1^*) (m \in {\Cal M})
$$
Furthermore, every element in ${\Cal M}$ is of the form
$nv_1$ with $n\in {\Cal N}$, namely
$$ m = d^2 \epsilon (m v_1^*)v_1 =  d^2 v_1^* \epsilon (v_1 m)
.$$
\endproclaim  
Now assume that two conformal precosheaf (cf. \S1.3) $N(I), M(I)$
are given, and use $H^0$ to denote the Hilbert space
of the vacuume representation of $M(I)$, and $\Omega\in H^0$
is the vacuum vector (uniquely determined up to a nonzero constant).
Assume that there is a covariant representation $\pi^0$ of  $N(I)$
on $H^0$ such that $ \pi^0(N(I)) \subset M(I)$, and moreover, 
fix any proper interval of $S^1$, 
the 
net $\pi^0(N(I)) \subset M(I)$ for any $I\subset J\subset S^1$ is
a directed standard net of subfactors w.r.t. $\Omega$. An immediate
corollary of the above proposition is the following:
\proclaim{Lemma 2.1.4}
Suppose the net $\pi^0(N(I)) \subset M(I)$ has finite index. Then
if $N(I)$ is strongly additive,  $M(I)$ is also strongly additive.
\endproclaim
\demo{Proof}
We need to show that if $I_1,I_2$ are the connected components
of an interval $I$ with one internal point removed, then
$$
M(I)\subset M(I_1)\vee M(I_2).
$$ By the above proposition we can choose isometry $v_1\in \pi^0(N(I_1))$ such
that any element of $M(I)$ can be written as $nv_1$ for some 
$n\in \pi^0(N(I))$. Since   $N(I)$ is strongly additive by our assumption,
$ \pi^0(N(I))= \pi^0(N(I_1))\vee \pi^0(N(I_2)) \subset 
M(I_1)\vee M(I_2)$, it follows that 
$$
M(I)\subset  M(I_1)\vee M(I_2)
$$ and the proof is complete.
\enddemo
\qed
\par
For the rest of this section, we assume that a net 
 $\pi^0(N(I)) \subset M(I)$ as described before lemma 2.1.4 with 
finite index is given. By Prop. 2.1.1  and Prop. 2.1.3 
we have a canonical endomorphism
$\alpha^{I_0}: M(I) \rightarrow N(I)$ for any $I\supset I_0$.
Recall $\beta^{I_0} := \alpha^{I_0}|N(I)$. For simplicity,
we shall drop the labels and write $\alpha^{I_0}, \beta^{I_0}$
simply as $\alpha, \beta$ in the following.  Define $M:= N(I_0)$. \par 
 
Define $U(\gamma, I): H^0 \rightarrow H_0$ to be a unitary
operator which commutes with the action of $N(I)$  
, where $H_0$ denotes the Hilbert
space of the vacuum representation $\pi_0$ of $N(I)$.   
Notice that such an operator always exists since
$\pi^0|N(I) $ is equivalent to that of
$\pi_0$. We shall think of  $U(\gamma, I)$ as an
element of 
$B(H^0)$ by identifying $H_0$ as a subspace of $H^0$.
Define:
$
\phi_I: B(H^0) \rightarrow B(H_0)
$ by
$$
\phi_I (m) = U(\gamma, I) m U(\gamma, I)^*, m\in  B(H^0)
$$
As on Page 368 of [X4] $\alpha$  can be chosen to be:
$$
\alpha(m) = \phi_I^{-1} (\phi_{I_0^c} (m)). \tag 2.1.1
$$
For any $a\in N(I)$, define:
$$
\gamma_I (a) := \phi_I  \beta  \phi_I^{-1}(a)
.$$ Notice since $ \phi_I^{-1}(N(I)) = \pi^0(N(I))$, $\gamma$ is well defined. 
By using (2.1), we have:
$$
\gamma_I (a) = \phi_{I_0^c}(\phi_I^{-1}(a)) \tag 2.1.2
$$
which is precisely the definition of localized (localized on $I_0$) 
covariant endomorphism
associated with the representation $H^0$ of $M(I)$ .\par
We shall denote $\gamma_{I_0}$ for the fixed interval $I_0$ by $\gamma$.
 Notice $\gamma \in End(M)$. \par
Since $\beta|_{N(I_0)}$ is a canonical endomorphism, we can find
$\rho_1 \in End(N(I_0))$ such that:
$$
\beta|_{N(I_0)} = \rho_1 \bar \rho_1 
$$
with $\bar {\rho_1} (N(I_0)) \subset N(I_0)$
conjugate to $N(I_0) \subset M(I_0)$.  Moreover, 
$$
 {\rho_1} (N(I_0)) = \alpha ( M(I_0))
$$ . Recall $M:=N(I_0)$. 
Define $\rho \in End(M)$  by:
$$
\rho   := \phi_{I_0}  {\rho_1}  \phi_{I_0}^{-1}
$$
Anticipating the 
dictionary below, we denote by  $j_{BE}$  
the localized covariant endomorphism (localized on
$I_0$) of $M(I)$ on $H_{j_{BE}}$. Let us recall how 
 $j_{BE}$ is defined. Let 
$U(j,I): H_{j_{BE}} \rightarrow H^0$
be a unitary map which commutes with the action of $L_IG$. Define:
$
\psi_I: B( H_{j_{BE}}) \rightarrow B(H^0)
$ by $\psi_I (x) = U(j,I) x U(j,I)^*$. Then  $j_{BE}$ is given by:
$$
j_{BE}(m) = \psi_{I_0^c}(\psi_I^{-1}(m)) \ (\forall m\in M(I))
.$$  Let $\gamma_j$ be the reducible representation of $N(I)$ on
$ H_{j_{BE}}$. Then the localized endomorphism, denoted by the same
 $\gamma_j$, is given by:
$$
 \gamma_j = \phi_{I_0^c}  \psi_{I_0^c}  \psi_{I}^{-1} 
 \phi_I^{-1}.     
$$ Define $\sigma_j \in End(M)$ by:
$$
\sigma_j = \rho^{-1} \phi_{I_0^c}  j_{BE}  
\phi_{I_0^c}^{-1}  \rho \tag 2.1.3
$$ Notice by the definition of $\rho_1$, 
$\rho_1(N(I_0)) = \alpha (M(I_0)) = \phi_{I_0}^{-1}  \phi_{I_0^c} (M(I_0))$,
and it follows that $\phi_{I_0^c}^{-1} \rho (M) = M(I_0)$, so  $\sigma_j$
as above is well defined.  Note that our $\sigma_j$ is 
$\sigma_j'$ in [X4] and our $j_{BE}$ is $\sigma_j$ in [X4].\par
For each covariant endomorphism $\lambda$ (localized on $I_0$) 
of $N(I)$ with finite index,  let $\sigma_{\gamma\lambda}
 : \gamma  \lambda \rightarrow
\lambda 
\gamma$ (resp.  $\tilde \sigma_{\gamma\lambda} : 
\gamma  \lambda \rightarrow
\lambda 
\gamma$), be the positive (resp. negative) 
braiding operator as defined in \S1.4 .
Denote by $\lambda_\sigma \in$ End$(M)$ which is defined by
$$\lambda_\sigma(x) := Ad(\sigma_{\gamma\lambda}^*)\lambda(x)= 
\sigma_{\gamma\lambda}^* \lambda(x) \sigma_{\gamma\lambda}
$$ for any $x \in M$.
By (1) of Th. 3.1 of [X4], $\lambda_\sigma  \rho (m) \in \rho(M)$

(similarly $\lambda_{\tilde\sigma}  \rho (m) \in \rho(M)$), and
hence the following definition makes sense.
\proclaim{Definition 2.1.1}
$a_\lambda$ is an endomorphism of $M$ defined by:
$$
a_\lambda(m):= \rho^{-1} (\lambda_\sigma  \rho (m)) , \
\tilde a_\lambda(m):= \rho^{-1} (\lambda_{\tilde\sigma}  \rho (m)),
\forall m\in M.
$$
\endproclaim
The endomorphisms $a_\lambda$ are called braided endomorphisms in [X4]
due to its braiding properties (cf. (2) of Cor. 3.4 of [X4]), and enjoy
a remarkable set of properties (cf. \S3 of [X4]). Motivated by 
[X4], these properties are also studied in a slightly different 
context in [BE1] and [BE2]. We will set up a dictionary between our
notations here and that of [BE1] and [BE2]. In the following the 
notations from  [BE1] and [BE2] will be given a subscript BE. 
The formulas are :
$$
\align
\gamma = \phi_{I_0}  \theta_{BE} \phi_{I_0}^{-1}, \ &
\alpha = \gamma_{BE}, \\
\lambda = \phi_{I_0} \lambda_{BE} \phi_{I_0}^{-1}, \ \ & \sigma_{\lambda
\gamma} =
\phi_{I_0} (\epsilon(\lambda_{BE}, \theta_{BE})), \\
\bar\rho (\sigma_
{\lambda\gamma}) =  \rho^{-1} \phi_{I_0^c} 
(\epsilon (\lambda_{BE}, \theta_{BE})) \tag 2.1.4
\endalign
$$
where is last formula is obtained by the following:
$$
\rho^{-1} \phi_{I_0^c} (\epsilon (\lambda_{BE}, \theta_{BE}))
= \rho^{-1} \gamma ( \sigma_{\lambda\gamma}) = \bar\rho (\sigma_
{\lambda\gamma})
$$ which follows from  (2.1.3) and $\rho \bar\rho= \gamma$.
Now we are ready to set up a dictionary between $a_\lambda \in End (M)$ 
in definition 2.1.1  and
$\alpha_\lambda^{-}$ as in Definition 3.3 and
3.5  of [BE1]. By using definitions
above we have:
$$
\align
a_\lambda  &= 
\rho^{-1} \phi_{I_0^c} \alpha_{\lambda_{BE}}^{-}    
\phi_{I_0^c}^{-1}  \rho \\ 
\tilde 
a_\lambda  &= \rho^{-1} \phi_{I_0^c} \alpha_{\lambda_{BE}}    
\phi_{I_0^c}^{-1}  \rho \tag 2.1.5
\endalign
$$ 
In fact 
$$
\align
\rho^{-1} \phi_{I_0^c} \alpha_{\lambda_{BE}}^{-}    
\phi_{I_0^c}^{-1}  \rho & = \rho^{-1} \phi_{I_0^c}
\alpha^{-1}  Ad (\epsilon^{-}(\lambda_{BE},\theta_{BE})) 
\lambda_{BE} \alpha \phi_{I_0^c}^{-1}
 \rho \\
&= \rho^{-1}  \phi_{I_0} Ad (\epsilon^{-}(\lambda_{BE},\theta_{BE})) 
\lambda_{BE}
\phi_{I_0}^{-1}
\rho \\
&= \rho^{-1} Ad (\phi_{I_0}(\epsilon^{-}(\lambda_{BE},\theta_{BE})))
\phi_{I_0} \lambda_{BE} \phi_{I_0}^{-1} \rho \\
&=\rho^{-1} Ad ( \tilde\sigma_{\lambda\gamma}) \lambda \rho \\
&= a_\lambda
\endalign
$$   
The second equation in (2.1.5) follows similarly.
Notice the similarity of between (2.1.5) and   (2.1.3). 
Formulas (2.1.4) and (2.1.5) above will be referred to as our {\it dictionary} 
between the notations of [X4] and that of [BE1] to [BE3]. \par
\subheading{\S2.2 Relative braidings from [BE3]}
We preserve the setup of \S2.1.1.
\proclaim{Lemma  2.2.1}
Let $p, q\in Hom (a_{\mu_1}, a_{\mu})$. Then
$$
\bar\rho (\sigma_{\mu\lambda})  p \bar\rho (\sigma_{\mu_1\lambda}^*)
= a_\lambda(p), 
\bar\rho (\sigma_{\lambda\mu})  \tilde a_\lambda (q)
 \bar\rho (\sigma_{\lambda\mu_1}^*)
= q
$$
\endproclaim
\demo{Proof}
The proof follows from line 9 (count from the top) on Page 385 of 
[X4] (In the diagram drawing below line 9 one needs to change 
$\mu$ to  $\mu_1$).  This proof can also be translated into the proof
of Lemma 3.25 of [BE3] by using our dictionary (2.1.4) and (2.1.5).
\enddemo
\qed
\par
Now let us define relative braiding as first introduced in \S3.3 of
[BE3]. Let $\tilde \beta, \delta\in End(M)$ be subsectors of
$\tilde a_\lambda$ and $a_\mu$. By lemma 3.3 of [X4], $[\tilde \beta]$ and
$[\delta]$ commute. Denote by $\epsilon_r (\tilde\beta, \delta)$ given
by:
$$
\align
\epsilon_r (\tilde\beta, \delta):&= s^* a_\mu (t^*) 
\bar\rho(\sigma_{\lambda\mu})
\tilde a_\lambda (s) t \in Hom (\beta \delta, \delta\beta) \\ 
\epsilon_r (\delta, \tilde\beta):&= \epsilon_r (\tilde\beta, \delta)^{-1}, 
\tag 2.2.1
\endalign 
$$ with isometries $t\in Hom (\tilde \beta, \tilde a_\lambda)$ and
$s\in Hom (\delta,  a_\mu)$. Note that the first formula above 
is exactly (10) of [BE3]
using (2.1.4), (2.1.5).
Also note from (2.2.1) we have
$$
\epsilon_r (\tilde a_\lambda, a_\mu) = \bar\rho(\sigma_{\lambda\mu}) \
, \ \epsilon_r ( a_\lambda, \tilde a_\mu) = \bar\rho(\tilde
\sigma_{\lambda\mu}) \tag 2.2.2
$$
\par
\proclaim{Lemma 2.2.2}
The operator $\epsilon_r (\beta, \delta)$ defined above does not 
depend on $\lambda, \mu$ and the isometries $s,t$ in the sense that, if
there are isometries $x\in Hom (\beta, \tilde a_\nu)$ and $y\in
Hom( \delta, a_{\delta_1})$, then 
$$
\epsilon_r (\beta, \delta) =  s^* a_{\delta_1}
 (t^*) \bar\rho(\sigma_{\nu\lambda_1})
\tilde a_\nu (y) x
$$
\endproclaim
\demo{Proof}
The proof follows from Lemma 2.2.1 as in the proof of Lemma 3.11 of 
[BE3] by using the dictionary (2.1.4) and (2.1.5).
\enddemo
\qed
\par
The following two lemmas  are Prop. 3.12 and Prop. 3.15 of [BE3],
translated into our notations by using the dictionary (2.1.4) and (2.1.5).
\proclaim{Lemma 2.2.3}
The system of unitaries of Eq. (6) provides a relative braiding between
representative endomorphisms of subsectors of 
$\tilde a_\lambda$ and $a_\mu$ in the sense that , if $\beta, \delta,
\omega, \xi$ are subsectors of $[\tilde a_\lambda], [ a_\mu],
[\tilde a_\nu], [a_{\delta_1}]$, respectively, then we have 
initial conditions 
$$
\epsilon_r (id_M, \delta) = \epsilon_r (\beta, id_M) = 1 
,$$
compositions rules
$$
\epsilon_r ( \beta \omega, \delta) = \epsilon_r ( \beta,\delta)
\beta( \epsilon_r (\omega, \delta)), 
\epsilon_r ( \beta,  \delta \xi) = \delta (\epsilon_r ( \beta, \xi))
 \epsilon_r (\beta, \delta)
,$$ and naturality
$$
\delta (q_+) \epsilon_r (\beta, \delta) =\epsilon_r (\omega, \delta) q_+,
q_-,  \epsilon_r (\beta, \delta) = \epsilon_r (\beta, \xi) \beta (q_-)
$$ whenever 
$q_+ \in Hom(\beta, \omega)$ and $q_-\in Hom (\delta, \xi)$.
\endproclaim
\proclaim{Lemma 2.2.4}
$\rho^{-1} \phi_{I_0^c} (\sigma_{i_{BE} j_{BE}}) = 
\epsilon_r (\tilde \sigma_i, \sigma_j)
$, where $\tilde \sigma_i$ means that the relative braiding 
$\epsilon_r (\tilde \sigma_i, \sigma_j)$ is defined so that 
its first argument is considered to be a subsector of some $\tilde
a_\lambda$ 
(and hence the second argument  is considered to be a subsector of some $
a_\mu$), and $\sigma_{i_{BE} j_{BE}}$ is the braiding operator as defined in
\S1.4 for the conformal precosheaf associated with $\{ M(I), \forall I\}$.    
\endproclaim
\subheading{\S2.3 Factorization of coset link invariants}
Let $H\subset G_k$ be an inclusion which satisfies the assumptions of
\S1.6. Let ${\Cal A}_{G/H}(I)$ be the conformal precosheaf of the coset as described
in \S1.6. 
Recall from \S4.2 of [X4] that $N(I):= {\Cal A}_{G/H}(I)\otimes \pi_0(L_IH)''$ and
$M(I):= \pi^0(L_IG)''$ verifies
all the assumptions of \S2.1 and \S2.2 so we can apply the results of
\S2.1 and \S2.2.  The sectors $\pi_{i,\alpha}$ of ${\Cal A}_{G/H}(I)$ are obtained
in the decompositions of $\pi^i$ of $LG$ with respect to subgroup $LH$,
and we will denote the set of such $(i,\alpha)$ by $exp$.
We will use the following tensor notation
\proclaim{Tensor Notation}
Let $\rho \in End( {\Cal A}_{G/H}(I)\otimes
\pi_0(L_IH)'')$. We will denote $\rho$ by $\rho_1\otimes \rho_2$
if
$$
\rho(x\otimes 1)= \rho_1(x)\otimes 1, \forall x\in {\Cal A}_{G/H}(I),
\rho(1\otimes y)= 1\otimes \rho_2(y), \forall y\in\pi_0(L_IH)'' ,
$$ where $\rho_1 \in End( {\Cal A}_{G/H}(I)),\rho_2 \in End(\pi_0(L_IH)'')$.
\endproclaim
So $(i,\alpha)\otimes \beta$ will be a covariant endomorphism of
$N(I)= {\Cal A}_{G/H}(I)\otimes \pi_0(L_IH)''$ obtained from 
the covariant endomorphisms $(i,\alpha)$ and $\beta$ of
$  {\Cal A}_{G/H}(I)$ and $ \pi_0(L_IH)''$ respectively. Note we have
$$
\omega((i,\alpha)\otimes \beta) = \omega((i,\alpha)) \omega(\beta).
$$  
For simplicity we shall
denote $(i,\alpha)\otimes 1$ (resp. $1\otimes \beta$) as
$(i,\alpha)$ (resp. $\beta$) where $1$  
stands for identity automorphism. \par
As in \S2.1 fix a proper  open interval $I_0\subset  S^1$.
For the rest of this section we will 
fix a choice  a finite set of irreducible endomorphisms of ${\Cal A}(I_0)$
which appear as irreducible 
subsectors of all $(i,\alpha)\in exp$ localized on $I_0$
and fix a choice of  all 
irreducible endomorphisms (denoted by $\beta$) 
of $\pi_0(L_{I_0}H)''$  obtained from 
the  irreducible covariant representations of the conformal precosheaf
associated with $H$ localized on $I_0$. 
Also as in \S2.1 fix a choice of $\sigma_i$, $\gamma$ 
and   
$\rho, \bar\rho$ such that $\rho\bar\rho = \gamma$. \par
Note  $[a_{(i,\alpha)}]$ is a subsector of $[\sigma_i a_{\bar\alpha}]$ by
(2) of Prop. 4.2 of [X1]. In fact the proof of (2) of Prop. 4.2 of [X1]
also shows that if $(i,\alpha)\in exp$, then  
$$
[\overline{(i,\alpha)}]= [(\bar i, \bar\alpha)] \tag 2.3.1
$$ and it follows that  $(\bar i, \bar\alpha) \in exp$. 
Let us give a proof of (2.3.1) following the proof of (2) of Prop. 4.2 of
[X1]. Assume that $[(i,\alpha)]= \sum_j m_j [x_j]$, where the sum is 
finite, $m_j\in {\Bbb N},$ and $x_j$ is irreducible. So
$[\overline{(i,\alpha)}]=\sum_j m_j [\bar x_j]$. Since
$$
\langle a_{x_j\otimes 1},  a_{1\otimes \bar \alpha} \sigma_i \rangle 
= \langle x_j, (i,\alpha) \rangle 
= m_j
$$  
by the proof of (2) of Prop. 4.2 of [X1], and  $[a_{\bar \lambda}]= 
[\bar a_{\lambda}]$ by 
(2) of Cor. 3.5 of [X4],  $[\bar\sigma_i] = [\sigma_{\bar i}]$
by definition, we have
$$
\langle a_{\bar x_j\otimes 1},  a_{1\otimes  \alpha} \sigma_{\bar i} \rangle 
= \langle \bar x_j, (\bar i,\bar \alpha) \rangle =m_j
$$
and it follows  that $[(\bar i,\bar \alpha)]\succ [\overline{(i,\alpha)}]$ for
any $(i,\alpha)\in exp$. This shows that $(\bar i,\bar \alpha)\in exp$
if  $(i,\alpha)\in exp$. 
Replacing $i, \alpha$ by $\bar i, \bar \alpha$ we get
$[( i, \alpha)]\succ [\overline{(\bar i,\bar\alpha)}]$, and so
$[\overline{(i,\alpha)}]\succ [(\bar i,\bar\alpha)]$, and this together
with $[(\bar i,\bar \alpha)]\succ [\overline{(i,\alpha)}]$ proves (2.3.1). \par

\proclaim{Proposition 2.3.1}
Suppose $u_{i,\alpha} \in Hom(a_{(i,\alpha)}, \sigma_i a_{\bar\alpha}),
 u_{j,\beta} \in Hom(a_{(j,\beta)}, \sigma_j a_{\bar\beta}),$
then:
(1) 
$$
\epsilon_r ( \tilde a_{\bar\alpha},  \sigma_j a_{\bar\beta})
\tilde a_{\bar\alpha}(u_{j,\beta}) = u_{j,\beta} \epsilon_r (
 \tilde a_{\bar\alpha},  a_{(j,\beta)});
$$ \par
$$
\epsilon_r ( a_{\bar\alpha}, \tilde \sigma_j\tilde a_{\bar\beta})
a_{\bar\alpha}(u_{j,\beta}) = u_{j,\beta} \epsilon_r (
 a_{\bar\alpha}, \tilde a_{(j,\beta)});
$$ \par
$$
\epsilon_r (\tilde \sigma_i \tilde a_{\bar\alpha}, a_{(j,\beta)})
u_{i,\alpha} = a_{(j,\beta)}(u_{i,\alpha})  \epsilon_r (\tilde 
a_{(i,\alpha)}, a_{(j,\beta)});
$$ \par
(2)
$$
\epsilon_r (
 \tilde a_{\bar\alpha},  a_{(j,\beta)}) = \epsilon_r (
 a_{\bar\alpha}, \tilde a_{(j,\beta)});
$$\par
(3)
Define 
$$
c(\sigma_i a_{\bar\alpha}, \sigma_j a_{\bar\beta}):=
\epsilon_r (\tilde \sigma_i, \sigma_j a_{\beta}) \sigma_i (
\epsilon_r (a_{\bar\alpha}, \tilde \sigma_j \tilde a_{\bar\beta})).
$$ Then
$$
c(\sigma_i a_{\bar\alpha}, \sigma_j a_{\bar\beta}) 
\sigma_i a_{\bar\alpha} (u_{j,\beta})u_{i,\alpha}
= u_{j,\beta} a_{(j,\beta)} ( u_{i,\alpha}) 
\epsilon_r (\tilde a_{(i,\alpha)}, a_{(j,\beta)}).
$$
\endproclaim
\demo{Proof}
Ad (1): The first  equation follows from the naturality
of relative braiding (cf. Lemma 2.2.3). The second equation also follows from
the naturality of relative braiding and the fact that
$$
\align
Hom(a_{(j,\beta)}, \sigma_j a_{\bar\beta}) =
Hom(a_{(j,\beta)}\bar\rho , \sigma_j a_{\bar\beta}\bar\rho)
& = Hom( \bar\rho(j,\beta),  \sigma_j \bar\rho \bar\beta) \\
& = Hom(\tilde a_{(j,\beta)}\bar\rho , \sigma_j \tilde a_{\bar\beta}\bar\rho)
\\
& = Hom(\tilde a_{(j,\beta)}, \sigma_j \tilde a_{\bar\beta}),
\endalign
$$ where we have used Cor. 3.2 and (1) of Th. 3.3 in [X4]. 
The third equation follows similarly. 
\par
Ad(2): By monodromy equation (cf. Prop. 1.4.3) 
$$
\sigma_{\bar\alpha, (j,\beta)} \sigma_{ (j,\beta),\bar\alpha }=
\exp{2\pi i (\Delta_{(j,\beta)\otimes \bar \alpha}   
-\Delta_{\bar\alpha} -\Delta_{ (j,\beta)})}=1,
$$ since 
$$
\Delta_{(j,\beta)\otimes \bar \alpha} = \Delta_{\bar\alpha} +\Delta_
{ (j,\beta)}.
$$ It follows that 
$$ 
\sigma_{\bar\alpha, (j,\beta)} = \tilde \sigma_{\bar\alpha, (j,\beta)}, 
$$ and by applying $\bar\rho$ to both sides and 
use definition (2.2.1) we obtain (2).\par
(3)
We have
$$
\align
\epsilon_r (a_{\bar\alpha}, \tilde \sigma_j \tilde a_{\bar\beta})
a_{\bar\alpha}(u_{j,\beta}) 
&= u_{j,\beta}  
\epsilon_r (a_{\bar\alpha}, \tilde a_{j,\beta})\\
&=  u_{j,\beta}  
\epsilon_r (\tilde a_{\bar\alpha},  a_{j,\beta})\\
&= \epsilon_r (\tilde a_{\bar\alpha},   \sigma_j  a_{\bar\beta})
\tilde a_{\bar\alpha}(u_{j,\beta})
\endalign
$$
where we have used (2) in the second = and (1) in the the first and last =. 
It follows that
$$
\align
c(\sigma_i a_{\bar\alpha}, \sigma_j a_{\bar\beta}) \sigma_i a_{\bar\alpha}
(a_{j,\beta}) &= \epsilon_r (\tilde \sigma_i, \sigma_j a_{\beta}) 
\sigma_i (
\epsilon_r (a_{\bar\alpha}, \tilde \sigma_j \tilde a_{\bar\beta})
a_{\bar\alpha}(u_{j,\beta})
) \\
&= \epsilon_r (\tilde \sigma_i, \sigma_j a_{\beta}) 
\sigma_i ( \epsilon_r (\tilde 
a_{\bar\alpha},  \sigma_j  a_{\bar\beta})
\tilde a_{\bar\alpha}(u_{j,\beta})
) \\
&= \epsilon_r (\tilde \sigma_i, \sigma_j a_{\beta}) 
\sigma_i ( \epsilon_r (\tilde 
a_{\bar\alpha},  \sigma_j  a_{\bar\beta}))
\tilde \sigma_i \tilde a_{\bar\alpha}(u_{j,\beta})
\\
&=  \epsilon_r (\tilde \sigma_i \tilde 
a_{\bar\alpha}, \sigma_j a_{\bar\beta})\tilde \sigma_i \tilde
 a_{\bar\alpha}(u_{(j,\beta)}) \\
&= u_{j,\beta}  \epsilon_r (\tilde \sigma_i \tilde 
a_{\bar\alpha},  a_{(j,\beta)}),  
\endalign
$$
where we have used the equation
derived above  in the second = and (1) and lemma 2.2.3 
in the rest of equalities (note by the convention of lemma 2.2.4, 
$\tilde \sigma_i=\sigma_i$) .
The proof is then complete by using the third  equation in (1).
\enddemo
\qed
\par
Note that in  (3) of the above proposition, the operator
$c(\sigma_i a_{\bar\alpha}, \sigma_j a_{\bar\beta})$
involve relative braidings between 
$\sigma_i$ and $a_\beta$,  $\sigma_j$ and $a_{\bar\alpha}$ 
as defined in \S2.2,
where $\sigma_i, \sigma_j$ are  considered as  subsectors
 of some $\tilde a_\mu, \tilde a_\nu$.
Since 
$$
\epsilon_r (\tilde \sigma_i, a_\beta) 
\epsilon_r ( a_\beta,\tilde \sigma_i )=1,
\epsilon_r (\tilde \sigma_j, a_{\bar\alpha}) 
\epsilon_r (a_{\bar\alpha} ,\tilde \sigma_j )=1
$$ by definition, 
we can think that the braidings between $\sigma_i$ and $a_\beta$,
$\sigma_j$ and $a_{\bar\alpha}$ are
``trivial''. Also note the braidings between $a_{\bar\alpha}$ and 
 $a_{\bar\beta}$ is the image under $\bar\rho$ of the negative
braiding between $ \bar\alpha$ and $\bar\beta$, and the 
 braidings between $\sigma_i$ and 
 $\sigma_j$ is the image under $\rho^{-1} \phi_{I_0^c}$ of the positive
braiding between  $\sigma_i$ and $\sigma_j$. These facts will be 
important in the proof of Th. B below.

\par
Now let us focus on the cases that $[\sigma_i a_{\bar\alpha}]
= [ a_{(i,\alpha)}]$ which holds for a general class of cosets
including diagonal cosets of type A (cf. Th. 4.3 of [X4]). 
By (2) of Prop.4.2 of [X1], this is equivalent to 
$d(i,\alpha)= d(i)d(\alpha)$. 
Fix  a choice of 
unitaries $u_{i,\alpha}\in Hom(a_{(i,\alpha)},\sigma_i a_{\bar\alpha})$ 
for each $(i,\alpha)\in 
exp$. 
The relations
between the braidings among
$ a_{(i,\alpha)}$ 
and the braidings among $\sigma_i a_{\bar\alpha}$ is given by (3) of 
Prop. 2.3.1.  Now we need to relate the dualities. 
Denote by $r_{\sigma_i}:=\rho^{-1} \phi_{I_0^c}(R_{i_{BE}}) \in 
Hom (id, \sigma_{\bar i}\sigma_i) , 
\bar r_{\sigma_i}:=\rho^{-1} \phi_{I_0^c}(\bar R_{i_{BE}}) \in 
Hom (id, \sigma_i \sigma_{\bar i}),$
and 
$
r_{a_\alpha} = \bar\rho(R_\alpha) \in Hom(id, a_{\bar \alpha} a_\alpha),   
\bar r_{a_\alpha} = \bar\rho(\bar R_\alpha) \in 
Hom(id, a_\alpha a_{\bar \alpha} ),
$ where $R$ and $\bar R$ are defined as in (1.2.1). 
Define:
$$
r_{i,\alpha} =      a_{\bar i,\bar\alpha} (u_{i,\alpha}^*)
u_{\bar i, \bar\alpha}^*
\epsilon_r(a_\alpha, \tilde\sigma_{\bar i})
a_\alpha( r_{\sigma_i})r_{a_{\bar\alpha}} ,
$$ and
$$
\bar r_{i,\alpha} =      a_{i,\alpha} (u_{\bar i,\bar \alpha}^*)
u_{ i, \alpha}^*
\epsilon_r(a_\alpha, \tilde\sigma_{\bar i})
\sigma_i ( \bar r_{\alpha}) \bar r_{i}. 
$$ Note by definition and Prop. 4.2 of [X1] 
$$
\align
r_{i,\alpha}\in Hom (id, a_{(\bar i, \bar\alpha)} a_{( i, \alpha)}
)& = Hom (\bar \rho, \bar \rho (\bar i, \bar\alpha)( i, \alpha)) \\
&= \bar\rho( Hom (1, (\bar i, \bar\alpha)( i, \alpha))) 
\endalign
$$
and similarly 
$
\bar r_{i,\alpha} \in \bar\rho( Hom (1,  ( i, \alpha)
(\bar i, \bar\alpha))) $, we may assume that
$$
r_{i,\alpha} = \bar\rho (R_{i,\alpha}),\ 
\bar r_{i,\alpha} = \bar\rho (\bar R_{i,\alpha}).
$$ The next lemma  shows that $R_{i,\alpha}, \bar R_{i,\alpha}$
indeed satisfy the duality conditions (thus justifying the notations).
\proclaim{Lemma  2.3.2}
(1): 
$$
 a_{\bar i,\bar\alpha} (u_{i,\alpha})r_{i,\alpha} =
u_{\bar i, \bar\alpha}^*
\epsilon_r(a_\alpha, \tilde\sigma_{\bar i})
a_\alpha( r_{\sigma_i})r_{a_{\bar\alpha}} ,
$$ 
$$
 a_{i,\alpha} (u_{\bar i,\bar \alpha})\bar r_{i,\alpha} =      
u_{ i, \alpha}^*
\epsilon_r(a_\alpha, \tilde\sigma_{\bar i})
\sigma_i ( \bar r_{\alpha}) \bar r_{i}. 
$$\par
(2):
$$
\bar R_{i,\alpha}^* (i,\alpha)(R_{i,\alpha}) =1_{(i,\alpha)},
R_{i,\alpha}^* (\bar i , \bar \alpha) (\bar R_{i,\alpha}) =
1_{(\bar i, \bar \alpha)}; 
$$ 
$$
||R_{i,\alpha}||=||\bar R_{i,\alpha}||=\sqrt{d((i,\alpha))}; 
$$
\par
(3)
$$
u_{i,\alpha} \bar \rho (\theta_{(i,\alpha)}) = \rho^{-1} \phi_{I_0^c}
(\theta_{i_{BE}})
\bar \rho (\theta_{\bar \alpha}^{-1}) u_{i,\alpha},
$$ where $\theta$'s are the twistings as constructed in \S1.7.
\endproclaim
\demo{Proof}
(1) follows from definitions. \par
Ad (2):
Since  braidings among $\sigma_i$ and $a_{\alpha}$ are relative braiding
and hence 
$$
\epsilon_r (\tilde \sigma_i, a_{\alpha}) \epsilon_r (a_{\alpha},
\tilde \sigma_i) =1,
$$ in another words, the braiding among 
$\sigma_i$ and $a_{\alpha}$  are trivial, by using the naturality 
of the relative braiding (cf. lemma 2.2.3) 
we have
$$
\bar r_{i,\alpha}^* a_{(i,\alpha)} ( r_{i,\alpha}) =
u_{i,\alpha}^* \bar r^*_{\sigma_{\bar i}} \sigma_i
( r_{\sigma_i}) 
\sigma_i(\bar r^*_{a_{\bar\alpha}}
a_{\bar\alpha}(r_{a_\alpha})) u_{i,\alpha}
= id_{a_{(i,\alpha)}},
$$ 
and this implies the first equation in (2) by definitions. 
The second equation is proved similarly. The last equation 
of (2) follows from
the definitions and $d((i,\alpha))= d(i)d(\alpha)$. \par 
(3) follows since both sides
are equal to 
$u_{i,\alpha} \exp(2\pi i (\Delta_i - \Delta_\alpha))$. 
\enddemo
\qed
\par
Note by (2) of lemma 2.3.2 we can choose our pairs $(R_{i,\alpha}, 
\bar R_{i,\alpha}) $ as in lemma 2.3.2 in our construction of 
$C(G/H)$ in \S1.7. Note that a different choices of the pairs $(R_{i,\alpha}, 
\bar R_{i,\alpha})$ will not change the value of link invariants
by (4) of lemma 1.7.4. \par
Denote by $L((i_1,\alpha_1),..., ((i_n,\alpha_n)))$ the link 
invariant of an oriented framed link $L$ with $n$ components colored
by $(i_1,\alpha_1), ...,  (i_n,\alpha_n)\in exp$.  
Since $$
L((i_1,\alpha_1),..., (i_n,\alpha_n)) =\bar\rho 
(L((i_1,\alpha_1),..., (i_n,\alpha_n)))
$$ 
and
$\bar\rho$ is an endomorphism,  we can think of 
$$
L((i_1,\alpha_1),..., (i_n,\alpha_n))
$$ as a link whose 
$k$-th component is colored by $a_{i_k, \alpha_k}, k=1,...,n$, and replacing
the operators corresponding to the tangles in $L$ by their
image under $\bar\rho$. We can do the same with links colored 
by $\alpha$'s. For the link $L$ whose $k$-th component is colored by
$i_k$, using the isomorphism $\rho^{-1}\phi_{I_0^c}$ we can
think of $L(i_1,..,i_n)$ as  a link whose 
$k$-th component is colored by $\sigma_{i_k}, k=1,...,n$, and replacing
the operators corresponding to the tangles in $L$ by their
image under $\rho^{-1}\phi_{I_0^c}$.   
\proclaim{Theorem B}
Assume a coset $H\subset G_k$ verifies the conditions of Theorem A. 
Let $L$ be an oriented framed link in three sphere with $n$ components.
If
$d((i_k, \alpha_k))= d(i_k) d(\alpha_k) , (i_k,\alpha_k)\in exp,
k=1,...,n, $ then
$$
L((i_1,\alpha_1),..., (i_n,\alpha_n)) =
L(i_1,...,i_n) \overline{L(\alpha_1,... \alpha_n)}.
$$
\endproclaim
\demo{Proof}
It is sufficient to show 
$$
\bar\rho (L((i_1,\alpha_1),..., (i_n,\alpha_n)))
= \rho^{-1}\phi_{I_0^c}(L(i_1,...,i_n)) 
\bar\rho (\overline{L(\alpha_1,... \alpha_n)}).
$$ By the remark before the theorem, 
we can think of each component of $L_k$ as colored by 
$a_{(i_k, \alpha_k)}$. Insert coupons
$u_{i,k}^*u_{i,k}=1 $  on $L_k$. By (3) of Prop. 2.3.1  and Lemma 2.3.2,
we can pull the coupon $u_{i,k}$ away from 
$u_{i,k}^*$ around the component $L_k$, and when
$u_{i,k}$ returns close to $ u_{i,k}^*$ from the other direction, 
using $u_{i,k}u_{i,k}^*=1$,   we get 
a link with components colored by $\sigma_{i_k} a_{\bar\alpha_k}$. 
We can effectively 
think of $L_k$ consists of two parallel parts with one part
colored by $\sigma_{i_k}$ and the other part colored by $a_{\bar\alpha_k}$.
By (3) of Prop. 2.3.1 , the braidings among $ \sigma_{i_k}$ and
$a_{\bar\alpha_l}$ are trivial, and by
using the naturality of the relative 
braiding  we can move all the strings
colored
by $\sigma_{i_k}$'s away from the link colored by $a_{\bar\alpha_k}$. 
The result is two links $L_1$, $L_2$, where $L_1$ colored by
 $\sigma_{i_k}$'s, and $L_1$ is identical to $L$, and 
$L_2$ is colored by $a_{\bar\alpha_k}$'s, but $ L_2$ is identical
to $L$ except all the overcrossing and the undercrossings of
 $L$ have been exchanged by (3) of Prop. 2.3.1.  
Note $L_1(\sigma_{i_1}, ..., \sigma_{i_n}) = 
\rho^{-1} \phi_{I^c_0} (L_1(i_1,..., i_n)) =
L(i_1,..., i_n)$, since $L_1(i_1,..., i_n)= L(i_1,..., i_n)\in {\Bbb
C}$. Similarly  $L_2 (a_{\bar\alpha_1},..., a_{\bar\alpha_n})
= L_2 (\bar\alpha_1, ..., \bar\alpha_n) $, and by Cor. 2.8.1 of
[Tu]   
$$L_2 (\bar\alpha_1, ..., \bar\alpha_n)
$$ is the same as 
$L_2'(  \alpha_1, ..., \alpha_n) $ where $L_2'$ is obtained from
$L_2$ by reversing the direction of each component. Hence $L_2'$ is 
exactly the negation of $L$ defined on Page 110 of [Tu]. It follows by 
Lemma 5.1.3 of [Tu] that 
$$
L_2'(  \alpha_1, ..., \alpha_n) = \bar L(\alpha_1, ..., \alpha_n)
$$ since $C(H)$ is a unitary modular category by Th. A.
So we have 
$$
 L((i_1,\alpha_1),..., ((i_n,\alpha_n))) =
L(i_1,...,i_n) \overline{L(\alpha_1,... \alpha_n)}
$$ and the proof is complete.
\enddemo
\qed
\par
\heading{\S3. Applications} \endheading
\subheading{\S3.1 $G=SU(N)_k, H={e}$}
In this section we show that the framed link invariants
from the category $C(G)$ with $G=SU(N)_k$ is the same as
the link invariants in  [W3]. This imply that $\tau_G(M)$ in this case
is the same as the 3-manifold invariants of [W3], [TW1] and [RT]. 
The proof is based on lemma 3.1.1 which slightly improves
(3) of Th. 3.8 of [X4]. \par
Let us recall the representations $\pi_n^{(N,N+k)}$ as defined on 
Page 368 of [W1]. This a $C^*$ representation of Hecke algebra
$H(n,q)$ which is an algebra with
generators $1,g_i, i=1,...,n-1$ and the following relations
$$
\align
g_i g_{i+1} g_i &=g_{i+1} g_{i} g_{i+1} \\
g_i g_j &= g_j g_i, \text{\rm if} \ |j-i|\geq 2\\
g_i^2 &= (q-1) g_i + q.
\endalign
$$ \par  
The minimal projections $q_{\tilde \lambda}$ of  $\pi_n^{(N,N+k)}(H_n(q))$
are labeled by $(N,N+k)$ Young diagrams $[x_1,...,x_N]$. The correspondence
between $q_{\tilde \lambda}$ and irreducible representations of
$LSU(N)_k$ is given by (cf. P.267 of [W2])
$$
\tilde \lambda=[x_1,...,x_N]\rightarrow\lambda= (x_1-x_N+1,...,x_{N-1}-x_N+1)
\tag 3.1.1
$$ 
Denote by $v$ the vector representation of $LG$. Consider the 
sequences of algebras
$$
End(v^n) \subset End(v^{n+1}),
$$ where the $\subset$ is the natural inclusion.
The minimal projection corresponding to $\lambda \prec f^n$ will
be denoted by $p_\lambda$. \par  

Denote by $A= \cup_n End(v^n)$. 
Define 
$$
h_i = q^{\frac{N+1}{2N}} v^{i-1} (c_{v,v}), i=1,2,...,
$$ where $q:= \exp(\frac{2\pi i}{N+k})$. Define a normalized 
trace Tr on $A$ by 
$$
Tr(f)= \frac{1}{d(v)^n} tr (f), \forall f\in End(v^n)
$$ where tr is defined as (1.7.1).  
\proclaim{Lemma 3.1.1}
(1)
The trace as defined above is a Markov trace on $A$, i.e., \par
$Tr(f h_{n}) = Tr(f) z, \forall f\in End(v^n)$, 
where $z= \frac{q-1}{1-q^{-N}}$. \par
(2)
$$
h_i^2 = (q-1) h_i + q , \forall i;
$$
(3) $1,h_1,...,h_{n-1}$ generates $End(v^n)$; \par 
(4) The map 
$$\psi_n: End(v^n) 
\rightarrow\pi^{(N,N+K)}(H_n(q)) 
$$ such that $\psi_n( h_i) =  \pi^{(N,N+K)} (g_i), i=1,2,...,n-1$
gives a trace-preserving $*$ isomorphism between the algebras
$End(v^n)$ and $\pi^{(N,N+K)}(H_n(q))$; \par
(5) Let $p_{2\Lambda_i}$ be the minimal projection in $End(v^i)$
corresponding
to $2\Lambda_i$, then 
$$\psi_i (p_{2\Lambda_i})
$$ is a minimal projection 
in $  \pi^{(N,N+K)}(H_i(q))$ corresponding to $2\Lambda_i$ via
(3.1.1). 
\endproclaim
\demo{Proof}
Ad (1): Since $C(G)$ is a ribbon category, we have
$$
tr(f h_n) = tr(f) q^{\frac{N+1}{2N}} \omega(v), \forall f\in End(v^n)
,$$ where $\omega(v)= \exp(\frac{2\pi i(N^2-1)}{2N})$. (1) now follows
from $d(v)= \frac{q^{N/2} - q^{-N/2}}{q^{1/2}-q^{-1/2}}$ and the 
definitions. \par
Ad (2): It is sufficient to show (2) for $i=1$. Note by (1.5.13)
$$
[v^2] = [3\Lambda_1]+ [2\Lambda_2]
.$$ Let $T_s\in  Hom(3\Lambda_1, v^2), T_a\in Hom(2\Lambda_2, v^2) $
be two isometries, uniquely determined up to a phase, and $p_s:= T_s T_s^*,
p_a:= T_a T_a^*$ ($s$ here stands for symmteric and $a$ stands for 
antisymmetric). Suppose $h =q^{\frac{N+1}{2N}} c_{v,v} = x p_a + y p_s$.
It is enough to show that $x=-1, y=q$. By Prop. 1.4.3 we have
$$
c_{v,v}^2 p_a = \omega(2\Lambda_2) \omega(v)^{-2} = q^{\frac{-N-1}{N}},
c_{v,v}^2 p_s = \omega(3\Lambda_1) \omega(v)^{-2} = q^{\frac{N-1}{N}},
$$ and so
$$
x^2= 1, y^2 = q^2.
$$ On the other hand by Spin-Statistical theorem (Th. 1.3.3)
$$
tr(h)= d(v) \omega(v), 
$$  and since $tr(h)= x d(2\Lambda_2) + y d(3\Lambda_1)$, 
we get the following equation 
$$
x d(2\Lambda_2) + y d(3\Lambda_1) = d(v) \omega(v)
,$$ and a little algebra using
$$
d(2\Lambda_2) = \frac{[N]_q[N-1]_q}{[2]_q}, 
d(3\Lambda_1) = \frac{[N]_q[N+1]_q}{[2]_q}
,$$ where $[m]_q:= \frac{q^{m/2} - q^{-m/2}}{q^{1/2} - q^{-1/2}}$, 
gives
$$
(x+1) [N-1]_q + (y-q)[N+1]_q = 0
.$$  It is then easy to  rule out the three cases $x=1, y=\pm q; x=-1, y=-q$. 
and we must have $x=-1, y=q$. \par
Denote by $B_n$ the subalgebra of $End(v^n)$ generated by 
$1,h_1,...,h_{n-1}$. 
From (1), (2) and (b) of Th. 3.6 of [W1]  the map
$\psi_n: B_n \rightarrow \pi^{(N,N+K)}(H_n(q))$
defined by $ \psi_n (h_i) =  \pi^{(N,N+K)}(g_i), i=1,2,...,n-1$ 
is a trace-preserving $*$ isomorphism.  So to prove (4) it is enough to
show (3), i.e., $B_n= End(v^n)$.  By Th. 3.6 of [W1], the simple ideals 
of $B_n$ are given by $(N,K+N)$ diagrams with $n$ boxes, which
are in one-to-one correspondence with the irreducible descendants of
$v^n$. The inclusion matrix of $B_n \subset B_{n+1}$, given by (2.14) on 
Page 369 of [W1], is precisely the same as the inclusion matrix of 
$End(v^n)\subset End(v^{n+1})$ by 
(1.5.13). Since $B_1= End(v)\equiv {\Bbb C}$, this
proves $B_n =End(v^n)$ and completes the proof of (3) and (4). \par
Ad (5): (5) is
obviously true for $i=1$, so let us assume $i>1$.  
First we claim that $h_m p_{2\Lambda_i} = - p_{2\Lambda_i},
k=m,..,{i-1}$. Note since $2\Lambda_i$ appears in $v^i$ once and only 
once, 
and by (1) of lemma 1.5.1, the only other possible
subsector of  $v^i$ with multiplicity 1 is $(i+1)\Lambda_1$.
By the fusion rule (1.5.13)
we can choose co-isometries $T_a^* \in Hom(v^2, 2\Lambda_2), T^*:
Hom(2\Lambda_2 v^{i-2}, 2\Lambda_i)$ such that
$$
p_{2\Lambda_i} =T_a TT^* T_a^*, 
$$ and it follows from (2) that
$h_1 p_{2\Lambda_i} = - p_{2\Lambda_i}. $
Hence   $\psi_i ( p_{2\Lambda_i})$ is minmal projection of 
$\pi^{(N,N+k)}(H_i(q))$ corresponding to the one dimensional simple 
ideal with the property 
$$ 
\pi^{(N,N+k)}(h_m) \psi_i (p_{2\Lambda_i}) 
= - \psi_i (p_{2\Lambda_i}), m=1,2,..., i-1.
$$
Similarly by using (2) one can show that if 
 $(i+1)\Lambda_1$ appear in  $v^i$ and $p_{(i+1)\Lambda_1}$ is the
corresponding minimal projection, then 
$$ \pi^{(N,N+k)}(h_m) \psi_i (p_{(i+1)\Lambda_1}) 
= q \psi_i (p_{(i+1)\Lambda_1}), m=1,2,..., i-1.
$$
Note that there are at most two one dimensional simple ideals in
$\pi^{(N,N+k)}(H_i(q))$ by (1) of lemma 1.5.1, it is then easy to
see that
$\psi_i ( p_{2\Lambda_i})$ corresponds to 
$q_{\tilde \lambda}$ where $\tilde \lambda = [1^i]$
(it is
one of the trivial $*$ representations defined on P. 370 of [W1])
. Note via (3.1.1)
$[1^i]$ corresponds to $2\Lambda_i$, so (5) is proved.    
\enddemo
\qed
\par  
(5) of Lemma 3.1.1 may seem to be redundant. However
note that the minimal projections $p_\lambda$ and
$q_{\tilde \lambda}$ as described before the lemma 
are all related to the  
irreducible representations of
$LSU(N)_k$, but it is not obvious that $ \psi(p_\lambda)$
should correspond to $\lambda$ via (3.1.1).  One must keep in 
mind that our construction of invariants  is  different from that
of [W2] even though it will  be shown in Prop. 3.1.2 that
they lead to the same invariants,  and (5) of lemma 3.1.1 is
used in the proof. \par  
Now assume that $L$ is  an oriented  framed link  in $S^3$ with 
$n$ components colored by $\lambda_1,..., \lambda_n$, and 
$M_L$ a 3-manifold obtained by surgery on $L$.  Denote by 
$L^Q(\lambda_1,..., \lambda_n)$ 
(resp. $\tau^Q(M_L)$)
the link invariant (resp. the 3-manifold invariant) defined by using the 
quantum group associated with $SU(N)$ with the deformation parameter 
$q= \exp(\frac{2\pi i}{N+k})$ (cf. [RT], [TW]), and let
$L^c(\lambda_1,..., \lambda_n)$ (resp. $\tau^c(M)$) 
be the  link invariant (resp. the 3-manifold invariant)  via cabling
as (*) on P. 247 of [W2]. 
\proclaim{Proposition 3.1.2}
(1)
$$
L(\lambda_1,..., \lambda_n) = L^c(\lambda_1,..., \lambda_n) = 
L^Q(\lambda_1,..., \lambda_n);
$$
(2)
$$
\tau(M) = \tau^c(M)= \tau^Q(M)
$$ for any closed oriented 3-manifold $M$. In fact it is also true if
$M$ is replaced by a pair $(M,\Omega)$ where $\Omega$ is a colored 
ribbon graph in $M$ as defined on Page 82 of [Tu].  
\endproclaim
\demo{Proof}
$$ L^c(\lambda_1,..., \lambda_n) = 
L^Q(\lambda_1,..., \lambda_n)
$$ follows from Cor. 3.3.3 of [W2], so we just have to show
$$
L(\lambda_1,..., \lambda_n) = L^c(\lambda_1,..., \lambda_n)
.$$ Since both sides transform in the same way under the change of 
orientation of a componenet of $L$ (the invararint is unchanged if one
change the orientation of a component and change the color of that
component by its dual), we can assume that $L$ is represented as 
the closure of a tangle $T_L$ as in the proof of (3) of Lemma 1.7.5.
Now let $p_{\lambda_i}$ be the minimal projection in 
$End(v^{m_i})$ corresponding to $\lambda_i$, and let 
$\psi_{m_i}(p_{\lambda_i})$ be a minimal projection 
in $\pi^{(N,k+N)} (H_i(q))$ where $\psi_{m_i}$ is defined as in
(4) of Lemma 3.1.1. Suppose that $\psi_{m_i}(p_{\lambda_i})$ corresponds
to $a(\lambda_i)$ which is a representation of $LSU(N)_k$ via (3.1.1). 
By using (4) of Lemma 3.1.1, the proof of (3) of Lemma 1.7.5 now applies
verbatim, and we have
$$
L(\lambda_1,...,\lambda_n) = L^c(a(\lambda_1),..., a(\lambda_n))
.$$ To finish the proof we just have to show that $a=id$. Note by 
(5) of lemma 3.1.1 $a(\Lambda_i) = \Lambda_i, i=1,...,N-1.$ On the
other hand since
$$
L(\lambda_1,...,\lambda_n) = L^c(a(\lambda_1),..., a(\lambda_n)) =
L^Q(a(\lambda_1),..., a(\lambda_n))
,$$ we must have
$$
S_{\lambda,\mu} = S_{a(\lambda),a(\mu)}, 
$$ since the elements of $S$ matrices are special cases of colored
link invariants
where the link is oriented downward, and by \S3.5 of 
[TW] the $S$ matrix computed from type $A$ quantum group at 
$q=  \exp(\frac{2\pi i}{N+k})$ agrees with the $S$ matrix defined at the 
end of \S1.5. So we have
$
S_{\lambda\mu} = S_{a(\lambda)a(\mu)},$ 
and $a(\Lambda_i) = \Lambda_i, i=1,...,N-1$. By lemma 1.5.1 $a=id$
and the proof of (1) is complete. (2) now follows immediately from 
definitions. 
\enddemo
\qed
\par
\subheading{\S3.2 Symmetry principle}
The symmetry principle in the case of $N=2$ first appeared in
[MK] and later generalized to $N>2$ case in [KT]. The proof in
[KT] uses the subtle property of $R$ matrix at roots of unity.
Here we give a proof of the
symmetry principle of our invariants $L(i_1,...,i_n)$, and by
Prop. 3.1.1 this also gives a slightly different proof of the 
results of [KT1].  This result will be used in \S3.6. \par
Note that 
$$
\Delta_{\sigma(\lambda)} - \Delta_\lambda - 
\frac{((N-1) k - 2\tau(\lambda))}{2N} \in {\Bbb Z}, \tag 3.2.1
$$  which follows from a simple computation using definitions. 
It follows from (3.2.1) that
$$
\Delta_{\sigma(\lambda)} - \Delta_\lambda - \Delta_{\sigma(1)}
+ \frac{\tau(\lambda)}{N} \in {\Bbb Z}, \tag 3.2.2
$$
$$
\Delta_{\sigma^{-1}(\lambda)} - \Delta_\lambda - \Delta_{\sigma^{-1}(1)}
- \frac{\tau(\lambda)}{N} \in {\Bbb Z}, \tag 3.2.3
$$ Also note from the definitions that
$$
\tau(\lambda) + \tau(\bar \lambda)\equiv 0 \text{\rm mod} (N) \tag 3.2.4
$$ 
For an oriented framed link $L$ with $n$ components $L_i, i=1,2,...,n$
the linking matrix is given by $(L_i\cdot L_j)$ where $L_i\cdot L_i$ is
defined to be the framing. 
\proclaim{Lemma 3.2.1}
Let $L$ be an oriented framed knot. Then
$$
L(\sigma(\lambda))= \exp(2\pi i (\Delta_{\sigma(\lambda)} -
\Delta_\lambda) L\cdot L) L(\lambda)
$$  
\endproclaim
\demo{Proof}
First assume that $L$ is presented as the closure of a counter-clockwise
braid, i.e., $L$ is the closure of a tangle $T_L$ with all bands 
oriented downwards, and we will also assume that $T_L$ consists of 
only twists and braidings.  
 Choose the blackboard 
framing so that $L\cdot L= b_+ + t_+-b_--t_-$, where $b_+$ (resp. $b_-$) 
are the number of positive (resp. negative) crossings in $T_L$, and
$t_+$ (resp.$t_-$) are the number of positive (resp. negative) 
twistings in $T_L$. \par 
Since $[\sigma(\lambda)] = [\sigma(1) \lambda]$
, by (3) of Lemma 1.7.4, we can replace the color 
$\sigma(\lambda)$ by $ \sigma(1) \lambda$, and by Cor. I.2.8.3 of
[Tu], we can obtain the same invariant $L(\sigma(\lambda))$ by cutting
the band $L$ along its core and coloring two newly emerging annuli 
denoted by $L_1, L_2$,  with
color $\lambda$ and $\sigma(1)$ respectively. Now for each positive
crossing or 
twisting of $L$ we can change the positive crossing of $\sigma(1)$ and
$\lambda$ to negative crossings of  $\sigma(1)$ and
$\lambda$ 
provided we multiply the resulting expression by 
$$
\exp(2\pi i (\Delta_{\sigma(\lambda)} - \Delta_\lambda - \Delta_{\sigma(1)})
),$$ by Prop. 1.4.3. Similarly for each negative
crossing or twisting of $L$ we can change the negative  crossing of  
$\lambda$ and $\sigma(1)$  to positive  crossings of  $\sigma(1)$ and
$\lambda$ 
provided we multiply the resulting expression by 
$$
\exp(-2\pi i (\Delta_{\sigma(\lambda)} - \Delta_\lambda - \Delta_{\sigma(1)})
),$$ by Prop. 1.4.3.  
Afterwards  it is
easy to see that we can then pull $L_1$ and $L_2$
apart, and  $L(\sigma(\lambda))$ is equal to
$$
\exp(2\pi i (\Delta_{\sigma(\lambda)} - \Delta_\lambda - \Delta_{\sigma(1)})
L\cdot L) L(\sigma (1)) L(\lambda).   
$$ Note that for each positive crossing of $L$ colored by 
$\sigma (0)$, the corresponding braiding operator is a scalar multiple of
identity since $End(\sigma(1)^2)\equiv {\Bbb C}$, and this scalar is equal
to $\omega(\sigma(1))= \exp(2\pi i \Delta_{\sigma(1)})$ 
by Th. 1.3.3 since
$d(\sigma(1))=1$.  Similarly  for each negative  crossing of $L$ colored by 
$\sigma (1)$, the corresponding braiding operator is a scalar multiple of
identity  and this scalar is equal
to $\omega(\sigma(1))^{-1}= \exp(-2\pi i \Delta_{\sigma(1)})$. 
So 
$$
L(\sigma (1))= \exp(2\pi i \Delta_{\sigma (1)} L\cdot L).
$$
It follows
that
$$
L(\sigma(\lambda))= \exp(2\pi i (\Delta_{\sigma(\lambda)} -
\Delta_\lambda) L\cdot L) L(\lambda)
$$ when $L$ is oriented as the closure of  counter-clockwise braids. 
Now suppose $L$ is oriented as the closure of clockwise braids. 
Denote by $L^{op}$ the knot obtained from $L$ by changing its direction.
From the above proof we have
$$
L^{op}(\bar\lambda)= \exp(2\pi i (\Delta_{\bar\lambda} -
\Delta_{\sigma^{-1}(\bar\lambda)}) L^{op}\cdot L^{op}) 
L^{op}(\sigma^{-1}(\bar\lambda))
,$$ and it follows that
$$
L(\sigma(\lambda))= \exp(2\pi i (\Delta_{\sigma(\lambda)} -
\Delta_\lambda) L\cdot L) L(\lambda)
$$ by using $\overline{\sigma(\lambda)} = \sigma^{-1}(\bar\lambda)$, 
$\Delta_\lambda = \Delta_{\bar \lambda}$  and $L^{op}\cdot L^{op} = 
L\cdot L$.  \par
\qed
\enddemo
\par
\proclaim{Proposition 3.2.2 (Symmetry principle)}\par
Suppose $L$  has $k$ components $L_1,...L_{k}$. Then
$$
\align
&L(\lambda_1,...,\sigma(\lambda_i),...,\lambda_{k})
= L( \lambda_1,...,\lambda_{k}) \times \\
&\exp(2\pi i ((\Delta_{\sigma(\lambda_i)}- \Delta_{\lambda_i}) L_i\cdot L_i
+ \sum_{j\neq i}(\Delta_{\sigma(\lambda_j)}- \Delta_{\lambda_j} - 
\Delta_{\sigma(1)})L_j\cdot L_i)) 
\endalign
$$ 
\endproclaim
\demo{Proof}
Without loss of generality we assume that $i=1$.
As in the proof of Lemma 3.2.1  
we can obtain the same invariant  by cutting
the band $L_1$ along its core and coloring two newly emerging annuli 
denoted by $N_2, N_1$,  with
color $\lambda$ and $\sigma(1)$ respectively. By using Lemma 3.2.1, 
$L(\sigma(\lambda_1), \lambda_2,...,\lambda_{k})$, up 
to  multiplication by 
$\exp(2\pi i (\Delta_{\sigma(\lambda_1)}- \Delta_{\lambda_1}) L_1\cdot L_1)
$ 
is the same as the invariant of a  link $L'$ 
where $N_1$ colored by  $\sigma(1)$
is unknotted and split from $N_2$ colored by $\lambda_1$,
and the rest of $L'$ is the same as $L$. 
Now suppose that $N_1$ and $L_j, j=2,...k$ are all oriented anti-clockwise.
The crossings of $N_1$ and $L_j$ occur in right or left handed pairs 
which algebraically sum to $L_1\cdot L_j$. 
For each righthanded pair of  crossings of $N_1$ and $L_j$, the associated 
operator is a scalar multiple of identity since 
$End(\sigma(1) \lambda_j) \equiv {\Bbb C}$, and the scalar is
$$
\exp(2\pi i (\Delta_{\sigma(\lambda_j)}- \Delta_{\lambda_j} - 
\Delta_{\sigma(1)}  ))
$$ by  Prop. 1.4.3, and similarly 
for each lefthanded pair of  crossings of $N_1$ and $L_j$, the associated 
operator is a scalar multiple of identity and the scalar is
$$
\exp(-2\pi i (\Delta_{\sigma(\lambda_j)}- \Delta_{\lambda_j} - 
\Delta_{\sigma(1)}  ))
,$$ so all told we have shown the proposition in the case that
all $L_i$ are oriented anti-clockwise. 
To finish the proof
we just have to show the phase factor of the equation in the proposition
remain the same under  changing  the  orientations of
components $L_j$ and in the same time  changing the color 
$\lambda_j$ to its dual $\bar\lambda_j$ when $j>1$, and 
changing the orientation of $L_1$ and in the same time changing 
the color $\lambda_1$ to $\bar\lambda_1$, $\sigma$ to $\sigma^{-1}$. 
For $j>1$ this
is equivalent to 
$$
(\Delta_{\sigma(\lambda_j)} -\Delta_{\lambda_j} - 
\Delta_{\sigma(1)}) + 
(\Delta_{\sigma(\bar \lambda_j)} -\Delta_{\bar \lambda_j} - 
\Delta_{\sigma(1)}) \in {\Bbb Z}
,$$ which follows from  immediately from (3.2.2), (3.2.3) and (3.2.4).
For $j=1$ this is equivalent to 
$$
(\Delta_{\sigma^{-1}(\lambda_j)} -\Delta_{\lambda_j} - 
\Delta_{\sigma^{-1}(1)}) + 
(\Delta_{\sigma( \lambda_j)} -\Delta_{ \lambda_j} - 
\Delta_{\sigma(1)}) \in {\Bbb Z}
,$$ which also 
follows from  immediately from (3.2.2), (3.2.3), (3.2.4)
and 
$$
\overline{\sigma(\lambda)} = \sigma^{-1}(\bar\lambda), 
\Delta_\lambda = \Delta_{\bar \lambda} ,L_1^{op}\cdot L_1^{op} = 
L\cdot L.
$$
\enddemo
\qed
\par
\subheading{\S3.3 Level-rank duality}
An almost immediate consequence of Th. B is a result of
[KT] relating two three manifold invariants. 
Our proof below is different from [KT].  
Let us first
prepare some notations.
The conformal inclusion  $SU(n)_m \times
SU(m)_n \subset  SU(nm)_1$ has been used in
[X1] (also. cf. [X7] and [X8]) and the branchings rules are given in [ABI].
  Let $\pi^0$ be the vacuum representation of $LSU(nm)$ on Hilbert space
$H^0$.  The decomposition of $\pi^0$ under $L(SU(m) \times SU(n))$ is
known, see, e.g. [ABI].  To describe such a decomposition, let us
prepare some notation.  We shall use $\dot S$ to denote the
$S$-matrices of $SU(m)$, (see \S2.2), and $\ddot S$ to denote the $S$-matrices
of $SU(n)$.  The level $n$ (resp. $m$) weight of $LSU(m)$ (resp. $LSU(n)$)
will be denoted by $\dot \lambda$ (resp. $\ddot \lambda$). \par
We start by
describing $\dot P_+^n$ (resp. $\ddot P_+^m$), i.e. the highest weights
of level $n$ of $LSU(m)$ (resp. level $m$ of $LSU(n)$).  
$\dot P_+^n$ is the set of weights
$$
\dot \lambda = \tilde k_0 \dot \Lambda_0 + \tilde k_1 \dot \Lambda_1 +
\cdots + \tilde k_{m-1} \dot \Lambda_{m-1}
$$
where $\tilde k_i$ are non-negative integers such that
$$
\sum_{i=0}^{m-1} \tilde k_i = n
$$
and $\dot \Lambda_i = \dot \Lambda_0 + \dot \omega_i$, $1 \leq i \leq m-1$,
where $\dot \omega_i$ are the fundamental weights of $SU(m)$.
 
Instead of $\dot \lambda$ it will be more convenient to use
$$
\dot \lambda + \dot \rho = \sum_{i=0}^{m-1} k_i \dot \Lambda_i
$$
with $k_i = \tilde k_i + 1$ and $\overset m-1 \to{\underset i=0 \to \sum}
k_i = m + n$ as in the 
notation of \S1.5.   Due to the cyclic symmetry of the extended Dykin diagram  
of $SU(m)$, the group $\Bbb Z_m$ acts on $\dot P_+^n$ by
$$
\dot \Lambda_i \rightarrow \dot \Lambda_{(i+\sigma)\mod m}, \quad
\sigma \in \Bbb Z_m.
$$
Let $\Omega_{m,n} = \dot P_+^n / \Bbb Z_m$.  Then there is a natural
bijection between $\Omega_{m,n}$ and $\Omega_{n,m}$ (see \S2 of
[ABI]).
 
We shall parameterize the bijection by a map
$$
\beta : \dot P_+^n \rightarrow \ddot P_+^m
$$
as follows.  Set
$$
r_j = \sum^m_{i=j} k_i, \quad 1 \leq j \leq m
$$
where $k_m \equiv k_0$.  The sequence $(r_1, \ldots , r_m)$ is decreasing,
$m + n = r_1 > r_2 > \cdots > r_m \geq 1$.  Take the complementary
sequence $(\bar r_1, \bar r_2, \ldots , \bar r_n)$ in $\{ 1, 2, \ldots ,
m+n \}$ with $\bar r_1 > \bar r_2 > \cdots > \bar r_n$.  Put 
$$
S_j = m + n + \bar r_n - \bar r_{n-j+1}, \quad 1 \leq j \leq n.
$$
Then $m + n = s_1 > s_2 > \cdots > s_n \geq 1$.  The map $\beta$ is
defined by
$$
(r_1, \ldots , r_m) \rightarrow (s_1, \ldots , s_n).
$$
The following lemma summarizes what we will .  For the proof,
see Lemma 3, 4 of [ABI].

\proclaim{Lemma 3.3.1} \rm{(1)}  Let $\dot Q$ be the root lattice of
$SU(m), \
\dot \Lambda_i, \ 0 \leq i \leq m-1$ its fundamental weights and
$\dot Q_0 = (\dot Q + \dot \Lambda_0) \cap \dot P_+^n$.  Then for
each $\dot \lambda \in \dot Q_0$, there exists a unique $\ddot \lambda
\in \ddot P_+^m$ with $\ddot \lambda = \sigma \beta(\dot \lambda)$
for some $\sigma \in \Bbb Z_n$ such that $H_{\dot \lambda} \otimes
H_{\ddot \lambda}$ appears once and only once in $H^0$.   
Moreover, $H^0$, as representations
of $L(SU(m) \times SU(n))$, is a direct sum of all such $H_{\dot \lambda}
\otimes H_{\ddot \lambda}$.
 
\rm{(2)}  $\underset \dot \lambda \in \dot Q_0 \to \sum \dot S_{\dot
\lambda \dot 1}^2 = \frac{1}{m}$.
 
\rm{(3)}  $\dot S_{\dot \lambda\dot 1} 
= \left( \frac{n}{m} \right)^{\frac{1}{2}}
\ \ddot S_{\sigma \beta (\dot \lambda) \ddot 1}$.
\endproclaim
We shall denote the bijection as in (1) of 
lemma 3.3.1 by $\dot \lambda \in \dot Q_0\rightarrow
\ddot \lambda = F(\dot \lambda) \in \ddot Q_0$. \par
Consider the coset $H=SU(m)_n \subset G=SU(nm)_1$. By using the above
decompositions, it is proved (cf. Lemma 3.2 of [X1]) that the coset
conformal precosheaf 
 is  the conformal precosheaf  of $H_1=SU(n)_m$, and it follows by
definition that $F(\dot \lambda)$ is isomorphic with 
$(1,\dot \lambda)$ as representations of the conformal precosheaf  of 
$H_1=SU(n)_m$, where $1$ denotes the identity sector of $G$,
$ \dot \lambda\in \dot Q_0$. By Th. B we have
$$
L(F(\dot \lambda_1),..., F(\dot \lambda_n))
=\overline{L(\dot \lambda_1,..., \dot \lambda_n)} L(1,...,1) .
$$ 
But $L(1,...1)=1$ so we have
$$
L(F(\dot \lambda_1),..., F(\dot \lambda_n))=
\overline{L(\dot \lambda_1,..., \dot \lambda_n)}.
$$
When   $m$ and $n$ are relatively prime, 
denote by 
$$\tau_{SU(n)_m/{\Bbb Z}_n}(M)$$ 
and  $$
\tau_{SU(m)_n/{\Bbb Z}_m}(M)$$ 
the three manifold invariants associated with $SU(n)_m/{\Bbb Z}_n$
and $SU(m)_n/{\Bbb Z}_m$ respectively  as defined on Page 226 of  [MW] (our 
notation $H_m/{\Bbb Z}_n$  corresponds to their $PSU(n)$ at level $m$)
or as on Page 260 of [KT2]. \par  
Since $C_H + C_{H_1} = mn-1$,
$$
\align
\sum_{\dot \lambda \in\dot Q_0} d(\dot \lambda)^2 =
\sum_{ \dot \lambda \in\dot Q_0} \frac{\dot S_{\dot \lambda\dot 1}^2}{
\dot S_{\dot 1 \dot 1}^2} & = \frac{1}{ m \dot S_{\dot 1}^2}\\
& = \frac{1}{ n \dot S_{\ddot 1}^2} \\
& = \sum_{\ddot \lambda \in\ddot Q_0} d(\ddot \lambda)^2
\endalign
$$ by  lemma 3.3.1, 
it follows by above  definitions   that
we have proved the following proposition:
\proclaim{Proposition 3.3.2}
If $m$ and $n$ are relatively prime, then
$$
\tau_{SU(n)_m/{\Bbb Z}_n}(M) = \overline{\tau_{SU(m)_n/{\Bbb Z}_m}(M)}
$$ for any closed oriented 3-manifold.
\endproclaim
Note by (2) of Prop. 3.1.2 the above proposition gives a different proof
of the main result Th. 4.2.7  of [KT2] which is based on symmetries
of certain Boltzmann weights.
\subheading{\S3.4 Simple current extensions of $\hat su(N)$} \par
Recall that the finite set of irreducible representations of $LSU(N)$ at
level $k$  is given by 
$$
 P_{++}^{h}
= \bigg \{ \lambda \in P \mid \lambda
= \sum _{i=1, \cdots , N-1}
\lambda _i \Lambda _i , \lambda _i \geq 1\, ,
\sum _{i=1, \cdots , n-1}
\lambda _i < h \bigg \}
$$
where $P$ is the weight lattice of $SU(N)$ and $\Lambda _i$ are the
fundamental weights and $h=N+k$.   
Recall that  this set admits a $\Bbb {Z}_N$ automorphism generated
by
$$
\sigma_1: \lambda=( \lambda_1, \lambda_2,..., \lambda_{N-1}) \rightarrow
\sigma_1( \lambda) =
(h-\sum_{j=1}^{N-1} \lambda_j, \lambda_1,...,\lambda_{N-2})
.$$  Note  $\tau(\lambda):\equiv \sum_i (\lambda_i-1)i
 \text{\rm mod} (N)$ and $Q$ is  the root lattice of $\widehat{SL(N)}$
(cf. \S1.3 of [KW]). Also note that $\lambda\in Q$ iff
$\frac{1}{N}\tau(\lambda) \in {\Bbb Z}$.
\par
Suppose  $N=mq$, with $m,q$ positive integers.
Assume level $k$ satisfies:
$$
kq \in 2m{\Bbb Z} \ \ \text{\rm if} N\in 2{\Bbb Z}; \ \
kq \in m{\Bbb Z} \ \ \text{\rm if} N\in 2{\Bbb Z}+1, \tag 3.4.1
$$ then there is an extension by the simple current $\sigma^q(1)$,
realized as standard net of subfactors (cf. \S2.1) as in Prop. 6.4
of [BE3].  The relation between condition (3.4.1) and orbifold subfactors
in the string algebra framework was noticed in [X5]. \par

Let $G_1:=SU(N)/{\Bbb Z}_m$. We shall denote the above nets of extensions
as described in Prop. 6.4 of [BE3] 
simply by $\hat G \subset \hat G_1$. 
We use $\hat G_1$ since it is likely that the net 
as described in Prop. 6.4 of [BE3] are related to the representations 
of loop group $LG_1$ which is not connected. The representation theory 
of such loop groups should be very close to that of [PS] but has not 
been developed yet. 
Note that the index of the 
inclusion  $\hat G \subset \hat G_1$ is $m$ (cf. \S6 of [BE]).    
Now we can apply the formulas in \S2.1 to this inclusion. Note that
this is not a coset construction,
but we shall see that Cor. 1.7.3  applies to the theory
determined by $\hat G_1$, and we can use the inclusion $\hat G \subset \hat 
G_1$
to calculate the 3-manifold invariants of $G_1$ in terms of $G$.\par
\proclaim{Lemma 3.4.1}
(1)
The conformal precosheaf  $\hat G_1$ is strongly additive and
$\mu$-rational. Its $\mu$-index, denoted by $\mu_{G_1}$ is given by
$$
\mu_{G_1} = \frac{\mu_{G}}{m^2}.
$$
(2) 
Let $ X$ be the finite set of all irreducible sectors of  $\hat G_1$, 
and let $D_-:= \sum_x d_x^2 \omega_x^{-1}$. Then 
$$
D_- = \sqrt{\mu_{G_1}} \exp(-\frac{\pi i C_G}{4})
,$$ where $C_G$ is the central charge of $LSU(N)$ at level $k$.
\endproclaim
\demo{Proof}
Ad (1):
Since $\hat G \subset \hat G_1$ has index $m$,  
$\hat G_1$ is strongly additive by lemma 2.1.4.  The vacuum representation of 
$\hat G_1$, denoted by
$\pi^0$,  decomposes into a finite number of irreducible representations
of $\hat G$, and the same argument showing that $\hat G$ is split, using
the asymptotics of the generator of the rotation proved in Th. B of [KW]
(cf. [BAF] or P. 21 of [X7]) now applies verbatim, and shows that
$\hat G_1$ is split.  Note
$$
\pi^0 (\hat G (E))'' \subset \pi^0 (\hat G_1 (E))'' \subset 
\pi^0 (\hat G_1 (E^c))' \subset \pi^0 (\hat G (E^c))'
,$$ where $E\subset S^1$  is  
union of two open intervals $I_1, I_2\subset S^1$ such that
the intersection of their closure in $S^1$ is empty.
Since $\pi^0 (\hat G (E))'' \subset \pi^0 (\hat G (E^c))'$
has finite index (cf. Th. 3.5 of  [X7]), it follows that
$$
\pi^0 (\hat G_1 (E))'' \subset 
\pi^0 (\hat G_1 (E^c))' 
$$ also has finite index, and so $\hat G_1$ has finite $\mu$ index.
Now applying   Prop. 21 of [KLM]  
gives the formula of the $\mu$-index. \par
Ad(2): By (1) and Th. 30 of [KLM] the 
set $X$ of all irreducible sectors of  $\hat G_1$ is finite. 
Apply the formalism of \S1.5 to this set $X$. 
By Cor. 32 of
[KLM], the $Y$-matrix as defined (1.5.0) is non-degenerate, 
and 
$$
D_- = \sqrt{\mu_{G_1}} \dot C^3
$$ for some $\dot C \in {\Bbb C}$ with $|\dot C |=1$. By
the remark after the proof of Prop. 3.1 of [X3]
(note in this case $\ddot C^3$ there is $1$  
since the coset $\hat G \subset
\hat G_1$ is trivial and $C^3$ there is our $k_G^{-3}$ )  , 
we have
$$
k_G^{-3} = \dot C^3, 
$$ and this finishes the proof of (2) by (1.7.5).
\enddemo
\qed
\par
\proclaim{Lemma 3.4.2}
(1) 
$$
\langle a_{\lambda},  a_{\mu} \rangle
=  \sum_{0\leq j\leq m-1} \delta(\lambda, \sigma^{jq} 
(\mu));
$$
(2)
$$
\langle a_{\lambda}, \tilde a_{\mu} \rangle
=  \delta^m (\tau(\lambda)) \sum_{0\leq j\leq m-1} \delta(\lambda, \sigma^{jq} 
(\mu))
$$ where $\delta^m (\tau(\lambda))$ is defined to be $1$ if
$\tau(\lambda) \equiv 0 \text{\rm mod} (m)$ and $0$ otherwise.
\endproclaim
\demo{Proof}
Both equations follow from similar equations in Th. 6.9 of 
[BE3], using the dictionary in \S2.1.
\enddemo
\qed
\par
Denote the  set $\lambda, \tau(\lambda)\equiv 0 \text{\rm mod}(m)$ 
 by $P_m$. Since   $ \tau(\sigma^q) = qk \equiv 0 \text{\rm
mod} (m)$
by (3.4.1),  $P_m$  admits a  ${\Bbb Z}_m$ action
generated by $\sigma^q$ , 
and decomposed into disjoint union of orbits, denoted by
$O_1,..., O_p$. Consider an orbit $O_s$. Suppose 
$O_s =\{ \lambda_s, \sigma^q (\lambda_s),...\sigma^{f(s)q} (\lambda_s) 
\}$.  Denote irreducible sectors of $a_\lambda, \lambda\in O_s$ 
by $\sigma_t(O_s), 1\leq t\leq g( O_s)$. By (2) of the above lemma 
$$
\langle \sigma_t(O_s), \sigma_{t'}(O_{s'}) \rangle =0
$$ if $s\neq s'$. Note that these  $\sigma_t(O_s)$ are in one-to-one
correspondence with some of the irreducible sectors of the extended net
$\hat G_1$, and with a slightly abuse of notations we will 
use $\sigma_t(O_s)$ to denote the irreducible sector of  $\hat G_1$.  
We will see in the following that they are in fact all
of the  irreducible sectors of the extended net
$\hat G_1$. 
We will (it is always possible ) choose 
an involution of our labels $s\rightarrow \bar s$ such that
$$[a_{\lambda_{\bar s}}] = [\overline{a_{\lambda_s}}]
,$$ and 
an involution of our labels $t\rightarrow \bar t$ such that
$$
[\overline{\sigma_t(O_s)}] = [\sigma_{\bar t}(O_{\bar s})]
.$$\par
Now consider the covariant irreducible representation $\pi^{\sigma_t(O_s)}$
of 
$\hat G_1$.  Consider the following analogue of Jones'
s basic construction:
$$
\pi^{\sigma_t(O_s)} (\hat G (I))'' \subset \pi^{\sigma_t(O_s)} (\hat G_1
(I))'' \subset \pi^{\sigma_t(O_s)} (\hat G_1 (I^c))'' 
\subset \pi^{\sigma_t(O_s)} (\hat G (I^c))''.
$$ 
Note that the minimal index of 
$$
\pi^{\sigma_t(O_s)} (\hat G_1
(I))'' \subset \pi^{\sigma_t(O_s)} (\hat G_1 (I^c))''
$$ is 
$d(\sigma_t(O_s))^2$, and the minimal index of 
$$
\pi^{\sigma_t(O_s)} (\hat G (I))'' \subset \pi^{\sigma_t(O_s)} (\hat G_1
(I))''
$$

and 
$$
\pi^{\sigma_t(O_s)} (\hat G_1
(I))'' \subset \pi^{\sigma_t(O_s)} (\hat G_1 (I^c))'' 
$$ are the same and equal to $m$. When restricting to 
$\hat G$, $\pi^{\sigma_t(O_s)}$ decomposes into a finite number of 
irreducible representations, and by our choices, each of such 
irreducible representation belong to the orbit $O_s$, and the 
multiplicity of any such $\sigma^{jq}(\lambda_s)$ is given by
$$
\langle \sigma_t(O_s), a_{ \sigma^{jq}(\lambda_s)}
\rangle = \langle \sigma_t(O_s), a_{\lambda_s}
\rangle
$$ by (2) of lemma 3.4.2. By the additivity of statistical dimensions
(cf. [L3]) the minimal index of the inclusion 
$$
\pi^{\sigma_t(O_s)} (\hat G (I))''  
\subset \pi^{\sigma_t(O_s)} (\hat G (I^c))''
$$ 
is 
$$
f(s)^2 d(\lambda_s)^2
(\langle \sigma_t(O_s), a_{\lambda_s}
\rangle)^2.
$$

On the other hand  by using the multiplicative 
properties of statistical dimensions (cf. [L3]) the minimal index of
$$
\pi^{\sigma_t(O_s)} (\hat G (I))''  
\subset \pi^{\sigma_t(O_s)} (\hat G (I^c))''
$$
is $m^2 d(\sigma_t(O_s))^2$, so
we have:
$$
d(\sigma_t(O_s)) = \frac{f(s)}{m} d(\lambda_s)
\langle \sigma_t(O_s), a_{\lambda_s}
\rangle \tag 3.4.2
$$ 
\proclaim{Lemma 3.4.3}
The irreducible sectors $\sigma_t(O_s), 1\leq t\leq g(O_s)$ above 
are all the irreducible sectors of $\hat G_1$.
\endproclaim
\demo{Proof}
It is sufficient to show that
$$
\sum_{O_s, 1\leq t\leq g(O_s)} d(\sigma_t(O_s))^2 = \mu_{G_1}
$$ by Th. 30 of [KLM]. By (3.4.2) we have
$$
\align
\sum_{O_s, 1\leq t\leq g(O_s)} d(\sigma_t(O_s))^2 & = 
\sum_{O_s, 1\leq t\leq g(O_s)}  \frac{f(s)^2}{m^2} d(\lambda_s)^2
\langle \sigma_t(O_s), a_{\lambda_s}
\rangle^2  \\
&= \sum_{O_s} \frac{f(s)^2}{m^2} d(\lambda_s)^2 \sum_{ 1\leq t\leq g(O_s)}
\langle \sigma_t(O_s), a_{\lambda_s}
\rangle^2 \\
&= \sum_{O_s} \frac{f(s)^2}{m^2} d(\lambda_s)^2 \langle a_{\lambda_s},
 a_{\lambda_s}
\rangle \\
&= \sum_{O_s} \frac{f(s)}{m} d(\lambda_s)^2 \\
&= \frac{1}{m}\sum_{ \lambda\in P_m}  d(\lambda)^2
\endalign
$$
where we used (cf. (2) of lemma 3.4.2)
$$
\sum_{ 1\leq t\leq g(O_s)}
\langle \sigma_t(O_s), a_{\lambda_s}\rangle = \langle a_{\lambda_s},
 a_{\lambda_s}
\rangle = \frac{m}{f(s)}
.$$ Let us show 
$$
\sum_{ \lambda\in P_m}  d(\lambda)^2 = \frac{1}{m}
\sum_{ \lambda \in P_{++}^h}  d(\lambda)^2 
.$$ By definitions
$$
\align
\sum_{ \lambda\in P_m}  d(\lambda)^2 &= \sum_{\lambda:
\tau(\lambda)\equiv 0 \text{\rm mod} (m)}  d(\lambda)^2 \\
&=  \sum_{1\leq i \leq q} \sum_{\lambda:
\tau(\lambda)\equiv mi \text{\rm mod} (N)}  d(\lambda)^2 \\
&= \sum_{1\leq i \leq q} \sum_{\lambda:
\lambda\equiv \lambda_i \text{\rm mod} (Q)}  d(\lambda)^2 \\
&= \sum_{1\leq i \leq q}  \frac{1}{N} \mu_G \\
&=\frac{1}{m} \mu_G  
\endalign
$$ where for each $1\leq i \leq q$ we choose a $\lambda_i$
such that $\tau(\lambda_i) \equiv mi \ \text{\rm mod} (N)$, 
and in the sixth $=$
we have used Cor. 2.7 of [KW]. It follows that
$$
\sum_{O_s, 1\leq t\leq g(O_s)} d(\sigma_t(O_s))^2 = \frac{1}{m^2} \mu_G
,$$ and by (1) of lemma 3.4.1 the proof is complete.
\enddemo
\qed
\par
Note that combine the above lemma and Cor. 32 of [KLW] gives a 
different proof
of Th. 6.12 of [BE3] (also cf. [EK3]). \par
Now let us construct a unitary modular category $C(G_1)$ as in 
\S1.7. We choose the set of objects $\rho_i, i\in I$ to be
$\sigma_t(O_s), 1\leq t \leq g(O_s)$ and $\sigma(O_s)$ such that
$$
[\sigma(O_s)] = \sum_{1\leq t \leq g(O_s)} 
\langle \sigma_t(O_s),
a_{\lambda_s} \rangle
[\sigma_t(O_s)].
$$ 
Note that since  $\sigma_t(O_s)$ has univalence
$\omega(\lambda_s)$ independent of $t$, so 
$\sigma(O_s)$ is of uniform valence with univalence $\omega(\lambda_s)$. 
By (1) of Cor. 1.7.3 
$C(G_1)$ is an abelian unitary ribbon category. Now using 
lemma 3.4.1, lemma 3.4.3 and Cor. 32 of [KLM] we conclude that 
$C(G_1)$ is a unitary modular category. \par
Note that $[\sigma(O_s)] = [a_{\lambda_s}]$ by our definition. Choose 
unitaries
$$
u(O_s)\in Hom (a_{\lambda_s},\sigma(O_s)).
$$ Define
$$
r(O_s):= \sigma(O_{\bar s}) (u^*(O_s)) u^*(O_{\bar s})  
\bar\rho (R_{\lambda_s}), 
\bar r(O_s):=  \sigma(O_{ s}) (u^*(O_{\bar s}))u^*(O_s)
 \bar\rho( \bar R_{\lambda_s}), 
$$ then we have the following analoge of (3) of Prop. 2.3.1 and 
Lemma 2.3.2 :
\proclaim{Lemma 3.4.4}
(1)
$$
u^*(O_t) \sigma(O_t) (u^*(O_s)) \sigma_{ \sigma(O_s) \sigma(O_t)} = 
\bar\rho{\sigma_{\lambda_s,\lambda_t}}  u^*(O_t)\sigma(O_t) (u^*(O_s))   
$$ where $\sigma_{xy}$ is the braiding operator 
between $x$ and $y$ as defined in \S1.4 ;\par
(2)
$$
\bar r(O_s)^* \sigma(O_s)( r(O_s)) = id_{ \sigma(O_s)},
r(O_s)^* \sigma(O_{\bar s})(\bar r(O_s)) =  id_{ \sigma(O_{\bar s})};
$$
(3)
$$
u^*(O_s) \theta_{\sigma(O_s)} = \bar\rho(\theta_{\lambda_s})u^*(O_s) ;
$$
(4)
$$ 
L(O_{s_1}, ..., O_{s_n}) =L(\lambda_{s_1}, ..., \lambda_{s_n}).   
$$
\endproclaim
\demo{Proof}
(1) follows from lemma 2.4. (2) follows from definitions. (3)
follows since $\omega(\sigma(O_s)) =\omega(\lambda_s)$. \par
To prove (4), 
insert  
coupons 
$$u(O_{s_k})^*u(O_{s_k}) =1,
$$ on the $k$-th componenet of $L, k=1,...,n$ as in the 
beginning of the proof of Th. B. By (1), (2) and (3) we   
 can pull the coupon $u_{i,k}$ away from 
 $ u_{i,k}^*$  around the component $L_k$, and when
$u_{i,k}$ returns close to $ u_{i,k}^*$ from the other direction,
using $u_{i,k}u_{i,k}^*=1$,   we get
a link identical to $L$, but with the 
$k$-th components colored by $a_{\lambda_{s_k}}, k=1,...,n$.
Since 
$$
L(a_{\lambda_{s_1}},..., a_{\lambda_{s_n}})
= \bar\rho(L(\lambda_{s_1},...,\lambda_{s_n}))
=L(\lambda_{s_1},...,\lambda_{s_n}),
$$ the proof of (4) is complete.
\enddemo
\qed
\par
Note
$$
\align
\sum_{1\leq t \leq g(O_s)} d(\sigma_t(O_s)) L(\sigma_t(O_s),...)
&=\frac{1}{m} \sum_{1\leq t \leq g(O_s)} f(s) \langle \sigma_t(O_s),
a_{\lambda_s} \rangle d(\lambda_s) L(\sigma_t(O_s),...) \\
&= \frac{1}{m} f(s)d(\lambda_s) L(\sigma(O_s),...), \tag 3.4.3
\endalign
$$ where in the second $=$ we have used (2) of Lemma 1.7.4
and the colors indicated by $...$ are kept fixed in the summation.  So we 
have
$$
\align
&\sum_{1\leq t_j \leq g(O_{s_j}), j=1,...n} d(\sigma_t(O_{s_1})) 
 ...  d(\sigma_t(O_{s_n}))
L(\sigma_{t_1}(O_{s_1}),\sigma_{t_2}(O_{s_2}), ..., \sigma_{t_n}(O_{s_n}) )
\\
&=\frac{1}{m^n}  f(s_1) f(s_2)...  f(s_n)d(\lambda_{s_1})...d(\lambda_{s_n})
L(\sigma(O_{s_1}),...\sigma(O_{s_n})) \\
&= \frac{1}{m^n}  f(s_1) f(s_2)...  f(s_n)d(\lambda_{s_1})...d(\lambda_{s_n})
L(\lambda_{s_1},...,\lambda_{s_n} ) \\
&= \frac{1}{m^n}\sum_{\lambda_j \in O_{s_j}, j=1,...,n}
d(\lambda_1) ...d(\lambda_n) L(\lambda_1,...,\lambda_n ) \tag 3.4.4 
\endalign
$$
where in the second and the third $=$ we have used (4) of Lemma 3.4.4.  
Apply the above to  the closed 3-manifold invariants  associated with
$C(G_1)$ we have
$$
\align
\tau_{G_1} (M_L)&= k_{G_1}^{3(b_-(L)- b_+(L))} \frac{1}{D_{G_1}^{n}}
\sum_{\sigma_{t_k}(O_{s_k}), k=1,...n} 
d(\sigma_{t_1}(O_{s_1}))... d(\sigma_{t_n}(O_{s_n})) \times \\
& L(\sigma_{t_1}(O_{s_1}), ..., \sigma_{t_n}(O_{s_n}))\\ 
& =  k_{G_1}^{3(b_-(L)- b_+(L))} \frac{1}{D_{G_1}^{n}}
\frac{1}{m^n} \sum_{s_j, \lambda_j \in O_{s_j}, j=1,...,n}
d(\lambda_1)...d(\lambda_n) L(\lambda_1 ,..., \lambda_n )\\
& =  k_{G}^{3(b_-(L)- b_+(L))} \frac{1}{D_{G}^{n}} 
\sum_{\lambda_1 \in P_m, ..., \lambda_n \in P_m}
d(\lambda_1) ...d(\lambda_n) L(\lambda_1 ,..., \lambda_n ) \\
&= k_{G}^{3(b_-(L)- b_+(L))} \frac{1}{D_{G}^{n}} \times \\
&\sum_{\lambda_j \in P_{++}^h, \tau(\lambda_j)\equiv 0 mod (m), j=1,...,n}
d(\lambda_1)...d(\lambda_n) L(\lambda_1 ,..., \lambda_n ), \tag 3.4.5
\endalign
$$ where we used lemma 3.4.1 that 
$k_{G_1}^{3}= k_{G}^{3}, D_{G} = m D_{G_1} $,  and 
 $D_{G}$ and $k_{G}^3$ are given by (1.7.4) and (1.7.5) respectively. 
\par
Let us consider a special case of (3.4.5) when $N=2$, $m=2$. Choose
$k=4k_1, k_1\in {\Bbb N}$ so that (3.4.1) is satisfied.  Note that there are 
$k_1+2$ simple objects in $C(G_1)$, 2 of them coming from the 2 irreducible
subsectors of $a_{k_1}$ by (2) of lemma 3.4.2, where as in \S3.5, we 
use half integers to denote the irreducible representations of $SU(2)$.  
It is easy to see
that our invariants $\tau_{G_1}$ in this case is exactly the invariants
discussed on Page 234 of [MW] when $l=k+2$ is even but not divisible by
$4$.  What is new here is a description of a unitary modular
tensor category 
(consequently a  3 dimensional Topological Quantum Field Theory
in the sense of   Chap. 4 of [Tu])
$C(G_1)$  underlying this invariant.  
\par

\subheading{\S3.5 Parafermion cosets}
Consider the coset $H\subset G_m$ with $G=SU(n)$ and $H$ is the Cartan
subgroup of $G$, a $n-1$ dimensional torus. This coset is an
example of parafermion cosets (cf. [DJ]). It is proved in \S4.4
of [X1] that
this coset verifies all the assumptions of Th.1, and hence we have
a 3 manifold invariant $\tau_{G/H}$. Also by \S4.4 of [X1]
the assumption of Th. B is also satisfied in this case, hence we
can effectively calculate these invariants. We will do this for the
case of $G=SU(2)$ in the following.\par
We will label the representations of $LSU(2)$ at level $k$ by $i$
such that $0\leq 2i\in {\Bbb Z} \leq k$. The representations of the
subalgebra $LU(1)$ at level $2k$ is labeled by $0\leq \alpha\in
{\Bbb Z} \leq 2k-1$. The fusion ring of  $LU(1)$ at level $2k$ is
isomorphic to the group ring of 
${\Bbb Z}_{2k}$ and the conformal dimension of
$\alpha$ is given by $\Delta_\alpha = \exp (2\pi i \frac{\alpha^2}{4k})$.
The central charge of this coset (cf. (1.7.6)) is given by
$C_{G/H}= \frac{2k-2}{k+2}$. \par
Note $(i,\alpha)\in exp$ iff $i+\frac{\alpha}{2} \in{\Bbb Z}$.  
By (2.6.12) of [KW], we know that
$(i,\alpha)\succ (0,0)$ iff $(i,\alpha)= (0,0) \ or \ (k/2, k)$.
On the other hand since $[a_{(k/2, k)}]=[\sigma_{k/2} a_{k}]$
has statistical dimension 1, it follows that
$$
[(k/2, k)]=[(0,0)]
.$$ This leads to the only selection rule that 
$$
[(i,\alpha)] =[(k/2-i, \alpha+ k)]. \tag 3.5.1
$$  
It is then easy to calculate the sum of the index of all 
different sectors $[(i,\alpha)]$, and the result is
$\frac{1}{2} k \mu_G$. But by Lemma 2.2 of [X3],
we have
$$
\mu_{G/H} = \frac{d(G/H)^4 \mu_G}{\mu_H},
$$ and since $ d(G/H)^4=k^2, \mu_H\geq 2k,   \frac{1}{2} k \mu_G\leq
\mu_{G/H}, $ we must have
$$
\mu_H= 2k, \mu_{G/H}= \frac{1}{2} k \mu_G
,$$ and by [KLM] this shows that $(i,\alpha)\in exp$ subject to
(1) are all the irreducible sectors of the coset $G/H$. \par
By Th. B, we have
$$
L((i_1,\alpha_1),..., (i_n,\alpha_n))
= L(i_1,...,i_n) \overline{ L(\alpha_1,..., \alpha_n)}.
$$ 
Note that the fusion ring generated by  $\alpha$'s  is isomorphic to 
the group ring of ${\Bbb Z}_{2k}$, and easier argument than that  in the 
proof of Prop. 3.2.2 shows that (also cf. [MOO]) 
$$
L(\alpha_1,..., \alpha_n) = \exp (\frac{2\pi i}{4k} \alpha L\cdot 
 \alpha L)
,$$ where $\alpha L:=\sum_i \alpha_i L_i$, and
$$
\alpha L\cdot 
 \alpha L := \sum_{1\leq i,j\leq n} \alpha_i  \alpha_j L_i\cdot L_j
.$$ \par
$L(i_1,...,i_n)$ is the same as the framed link invariants in [KM]
by Prop. 3.1.2.\par
The three manifold invariant is given by
$$
\align
\tau_{G/H} (M_L) & = k_{G/H}^{3(b_-(L) - b_+(L))} D_{G/H}^{-n}
2^{-n}
\sum_{ i_j+ \frac{\alpha_j}{2}\in 
{\Bbb Z}, j=1,...,n} d({i_1})... d({i_n}) \times
\\
& L(i_1,...,i_n)
\exp (-\frac{2\pi i}{4k} \alpha L\cdot 
 \alpha L) \\ 
& = k_{G/H}^{3(b_-(L) - b_+(L))} (D_{G}D_H)^{-n}
\sum_{ i_j+ \frac{\alpha_j}{2}\in 
{\Bbb Z}, j=1,...,n} d({i_1})... d({i_n}) \times
\\
& L(i_1,...,i_n)
\exp (-\frac{2\pi i}{4k} \alpha L\cdot 
 \alpha L),   \tag 3.5.2
\endalign
$$
where $D_G= \frac{1}{\sqrt{\frac{2}{k+2}} \sin(\frac{\pi}{k+2})},
D_H= \sqrt{2k}$, and $k_{G/H}=\exp(2\pi i \frac{2k-2}{k+2})$, which
follows from (1.7.3) and (1.7.4). The $2^{-n}$ factor appears
in the first $=$ due to (3.5.1).  
Note that  the summation above  is over all those
representation labels given at the beginning of this 
section subject to the indicated constraints. \par
When $k=2k_1+1$ is odd, we can choose a unique representative of $(i,\alpha)$
in the equivalent class obtained by relation (3.5.1)
such that $2i$ is even, and in  this case $\alpha$ is also even
since $i+\frac{\alpha}{2} \in {\Bbb Z}$. So
$$
\align
\tau_{G/H} (M_L) & = k_{G/H}^{3(b_-(L) - b_+(L))} D_{G/H}^{-n}
 \sum_{i_j\in{\Bbb Z}, j=1,...,n
} d({i_1})... d({i_n}) L(i_1,...,i_n) \times \\
& \sum_{ \alpha_j\in 2{\Bbb Z},
j=1,...,n} \exp (-\frac{2\pi i}{4k} \alpha L\cdot 
 \alpha L) \\
&=2^{n/2} k_{G}^{3(b_-(L) - b_+(L))}D_{G}^{-n}
\sum_{ i_j\in{\Bbb Z}, ..., j=1,...,n} 
d({i_1})... d({i_n}) L(i_1,...,i_n) \\
&\times 2^{n/2}
k_{H}^{3(b_-(L) - b_+(L))}D_{H}^{-n} \sum_{
\alpha_j\in 2{\Bbb Z}, j=1,...,n} 
\exp (-\frac{2\pi i}{4k} \alpha L\cdot 
 \alpha L),
\endalign 
$$ note that 
$$
2^{n/2} k_{G}^{3(b_-(L) - b_+(L))}D_{G}^{-n}
\sum_{ i_1\in{\Bbb Z}, 
j=1,...,n} d({i_1})... d({i_n}) L(i_1,...,i_n)
$$ is $\tau_{SO(3)}(M)$ (cf. [MW] or [KT2]) , and
$$  
2^{n/2}
k_{H}^{3(b_-(L) - b_+(L))}D_{H}^{-n} \sum_{
 \alpha_j\in 2{\Bbb Z},j=1,...,n } 
\exp (-\frac{2\pi i}{4k} \alpha L\cdot 
 \alpha L)
$$ is also a 3 manifold invariant 
as defined on P. 545 of [MOO] with $N=k$,
and we denote it by  $\tau_{U(1)/{\Bbb Z}_2}$. Hence
$$
\tau_{G/H} (M) = \tau_{SO(3)}(M) \overline{\tau_{U(1)/{\Bbb Z}_2}(M)}. 
\tag 3.5.3
$$ Note that $\tau_{U(1)/{\Bbb Z}_2}$ depends only on the linking matrix,
and it follows from Page 521 of [KM] that it is homotopy invariant,
determined in fact by the first Betti number of $M_L$ and the linking
pairing on $Tor H_1 (M_L)$.  In fact  $\tau_{U(1)/{\Bbb Z}_2}$ is 
expressed in terms of classical invariants in [MOO].\par
However when $k$ is even, it is not clear if $ \tau_{G/H}$ can
be factorized into the products of other invariants. One may try to
express $\tau_G (M) \overline{\tau_H(M)}$ in terms of $ \tau_{G/H}$, but this
does not seem to be possible for general $M$. To illustrate this let
us consider the case that $k=2k_1$ with $k_1$ odd, and $L$ has 
only one component with framing $n=2n_1,$ and $n_1$ is odd.  Consider
the transformation $\alpha\rightarrow\alpha+ k_1$ which changes the
parity of $\alpha$, and  since
$$
\exp (\frac{2 \pi i}{4k} (\alpha+ k_1)^2n) =
\exp (\frac{2 \pi i}{4k} \alpha^2n) (-1)^{\alpha n_1} 
\exp(\frac{\pi ik_1n_1}{2})
,$$ we have
$$
\sum_{\alpha\in 2{\Bbb Z}} \exp (\frac{2 \pi i}{4k}\alpha^2n)
= -\exp(\frac{\pi ik_1n_1}{2})
\sum_{\alpha\in 2{\Bbb Z}+1} \exp (\frac{2 \pi i}{4k}\alpha^2n)
.$$ Denote the Gauss sum
$$
\sum_{0\leq \alpha\leq 2k-1} \exp (\frac{2 \pi i}{4k}\alpha^2n)
=  \sum_{0\leq \alpha\leq 4k_1-1} \exp (\frac{2 \pi i}{4k_1}\alpha^2n_1)
$$ by $G(n_1, 4k_1)$, note that  $G(n_1, 4k_1)\neq 0$ (cf. [Lang]).
From the above we have
$$
\sum_{\alpha\in 2{\Bbb Z}} \exp (\frac{2 \pi i}{4k}\alpha^2n)
=\frac{1}{(1- \exp(\frac{-\pi ik_1n_1}{2}))}G(n_1, 4k_1)
,$$ and
$$
\sum_{\alpha\in 2{\Bbb Z}+1 } \exp (\frac{2 \pi i}{4k}\alpha^2n)
=\frac{1}{(1- \exp(\frac{\pi ik_1n_1}{2}))}G(n_1, 4k_1)
.$$ So
$$
\align
\tau_{G/H}(M_L) &= \frac{1}{2} k_{G/H}^{3(b_-(L) - b_+(L))} D_{G/H}^{-1}
(\sum_{j, 2j\in{\Bbb Z}} L(j) - 
 \sum_{j, 2j\in{\Bbb Z}+1} L(j) \exp (\frac{\pi i k_1 n_1}{2}))
\times \\
& \frac{\overline{ G(n_1, 4k_1)}} 
{ (1- \exp(\frac{\pi ik_1n_1}{2}))}.
\endalign
$$ Compared to
$$
\tau_{G}(M_L) \bar \tau_{H}(M_L) = 
k_{G/H}^{3(b_-(L) - b_+(L))} D_{G}^{-1} D_{H}^{-1}
(\sum_{j, 2j\in{\Bbb Z}} L(j) + 
 \sum_{j, 2j\in{\Bbb Z}+1} L(j)) 
{\overline{G(n_1, 4k_1)}}, 
$$  
expressing $ \tau_{G} \bar \tau_{H}$  in terms of 
$\tau_{G/H}$ would
imply a relation between 
$$
\sum_{j, 2j\in{\Bbb Z}} L(j) - 
 \sum_{j, 2j\in{\Bbb Z}+1} L(j) \exp (\frac{\pi i k_1 n_1}{2})
$$ and
$$
\sum_{j, 2j\in{\Bbb Z}} L(j) + 
 \sum_{j, 2j\in{\Bbb Z}+1} L(j) 
,$$ which seems unlikely since when $k$ is even the symmetry
$j\rightarrow k/2-j$ preserves the
parity of $2j$, and does not relate the above two sums, so the
symmetry principle as on Page 513 of [KM] (also cf. 
Prop. 3.2.2 ) does not help. 
Also since in this case $k=2k_1$ is not divisible by $4$, the
$SO(3)$ invariants as discussed in the end of \S3.4  
does not exist. It is not clear if $\tau_{G/H}$ is related to any
of the previous known invariants in a simple way.
\par
\subheading{\S3.6 Diagonal cosets of type A}
We consider the coset 
$$H:=SU(N)_{m'+m''}\subset 
G:=SU(N)_{m'} \times SU(N)_{m''},
$$ where
the embedding $H\subset G$ is diagonal. 
We
use $i$ (resp.$\alpha$) to denote the irreducible positive
energy representations of $LG$ (resp.$LH$).  To
compare our notations with that \S2.7 of 
[KW], note that our $i$ is $(\Lambda',\Lambda'')$
of [KW] , and our $\alpha$ is $\Lambda$ of [KW]. 
We will  identify $i=(\Lambda',\Lambda'')$ and $\alpha =\Lambda$
where $\Lambda',\Lambda''$, $\Lambda$ are the weights of 
$SL(N)$ at levels $m', m'', m'+m''$ respectively.
Suppose 
$$
i=({\Lambda_1}',{\Lambda_1}''), j=({\Lambda_2}',{\Lambda_2}''),
k=({\Lambda_3}',{\Lambda_3}''), \alpha= {\Lambda_1}, \beta={\Lambda_2},
\delta=\Lambda_3.
$$  

Then the fusion coefficients
$N_{ij}^k := N_{{\Lambda_1}'{\Lambda_2}'}^{\Lambda_3'}
N_{{\Lambda_1}''{\Lambda_2}''}^{{\Lambda_3}''}$
(resp. $N_{\alpha\beta}^\delta := N_{\Lambda_1\Lambda_2}^{\Lambda_3} $
)of $LG$ (resp. $LH$)
are given by Verlinde formula (cf. \S1.5).
Recall $\pi_{i,\alpha}$ are the covariant representations of
the coset $G/H$. The set of all $(i,\alpha):=
(\Lambda',\Lambda'',\Lambda)$
which appears
in the decompositions of $\pi^i$ of $LG$ with respect to $LH$
is denoted by $exp$. This set is determined on P.194 of
[KW] to be $(\Lambda',\Lambda'',\Lambda) \in exp$ iff
$\Lambda'+\Lambda''-\Lambda \in Q$, where $Q$ is the root lattice (cf. 
\S3.4). 
The  ${\Bbb Z_N}$ action on $(i,\alpha), \forall i,
\forall \alpha$  is denoted by $\sigma(i,\alpha):=
(\sigma (\Lambda'), \sigma (\Lambda''), \sigma (\Lambda))$,
$\sigma \in  {\Bbb Z_N}$. This is also known as diagram automorphisms
since they corresponds to the automorphisms of Dykin diagrams.  Note that
this  ${\Bbb Z_N}$ action preserves $exp$ and therefore induces a
${\Bbb Z_N}$ action on  $exp$.
\par
By Th.4.3 of [X1], this coset verifies the conditions of Th. A and
Th. B. Let us first calculate $\mu_{G/H}$.  By lemma 1.6.1 
$$
\mu_{G/H} =\frac{d(G/H)^4 \mu_G}{\mu_H}.
$$
By Th. 4.3 of [X1]:
$$
d(G/H)^2 = \sum_{\alpha \in Q} d_\alpha^2 = \frac{1}{N} \mu_H
,$$  so
$$
\mu_{G/H} =\frac{\mu_H \mu_G}{N^2}
.$$  A direct way of calculating 
the sum of index of all irreducible sectors
of $(i,\alpha) \in exp$ can be found  in (2) of Th. 2.3 of
[X2].  From (1) of Th. 2.3 we  have 
$$
k_{G/H}^3= (k_G/k_H)^3,
$$ and this also follows from (1.7.3). \par

To calculate $\tau_{G/H}$ we need to prepare some further
notations from [X3]. 
We define a vector space $W$ over $\Bbb {C}$ whose
orthonormal basis are denoted by
$i\otimes \alpha$ with $i=(\Lambda',\Lambda'') , \alpha=\Lambda$. $W$
is also a commutative ring with structure constants given by
$N_{ij}^k N_{\alpha\beta}^\delta$. Let
$V$ be the vector space over $\Bbb {C}$ whose basis are given by
the irreducible components of $\sigma_i a_{1\otimes \bar \alpha}$
(cf. \S4.3 of [X4]) .
Then $V= V_0 \oplus V_1$, where $V_0$ is a subspace of $V$ whose
 basis are given by
the irreducible components of $\sigma_i a_{1\otimes \bar \alpha}$
with $(i,\alpha)\in exp$, and $V_1$ is the orthogonal complement of
$V_0$ in $V$. The composition of sectors gives  $V$  a ring structure.
   By (1) of theorem 4.3 of [X4], the irreducible
subrepresentations of $(i,\alpha)$ of the coset are in one-to-one
correspondence with the basis of $V_0$ and this map is a ring
isomorphism by (1) of Prop.4.2 of [X4], and we will identify
 the irreducible
subrepresentations of $(i,\alpha)$ of the coset with the basis of
$V_0$ in the following when no confusion arises.  Note that $V_0$ is
a subring of $V$ and $V_0.V_1\subset V_1$. \par
Define a linear map $P: W\rightarrow V$ such that $P(i\otimes \alpha) =
\sigma_i a_{1\otimes  \bar \alpha}$. By Th. 4.3 of [X4]
$$
P(i\otimes \alpha) =
\sigma_i a_{1\otimes  \bar\alpha} = P(i'\otimes \bar\alpha') =
\sigma_{i'} a_{1\otimes  \alpha'}
$$ iff $\sigma^s(i)=i', \sigma^s(\alpha)=\alpha'$ for some $s\in
{\Bbb Z}$.  Also
$\langle P(i\otimes \alpha), P(j\otimes \beta) \rangle= 0
$ if  $P(i\otimes \alpha)\neq P(j\otimes \beta)$ by (*) of \S4.3 of
[X4].
Note $P(\sigma (1) \otimes \sigma (1)) = 1$ and
$P$ is  a ring homomorphism from  $W$ to $V$.
Define $W_0:= P^{-1} (V_0), W_1:= P^{-1} (V_1)$, then
$W= W_0 \oplus W_1$ since $exp$ is $\sigma$ invariant.  Note that
$i\otimes \alpha \in W_0$ iff $i- \alpha \in Q$.
Define the action of $Z_N$ on $W$ as $\sigma(i\otimes \alpha)
= \sigma(i) \otimes \sigma(\alpha)$.
\par
Assume
$\sigma^s(i\otimes \alpha) = i\otimes \alpha$ for some $i\otimes \alpha
\in W_0$,
and $0<s\leq N$ is the least positive integer with this property.  Let
$t=\frac{N}{s}$. By equation (*) on Page 30 of [X4] we have:
$$
\langle P(i\otimes \alpha), P(i\otimes \alpha) \rangle =t
.$$  
In particular when $t>1$, i.e., when the subgroup of ${\Bbb Z}_N$ 
which fix $(i\otimes \alpha)$ is nontrivial, 
the sector $ P(i\otimes \alpha)$ is not
irreducible. The question of  
 decomposing $P(i\otimes \alpha)$ when
$t>1$ into irreducible pieces is known as {\it Fixed point 
resolutions} which is discussed in [LVW], [FSS1] [FSS2]. 
The problem is answered in a mathematical framework in [X2].  The
result is (cf. Page 10 of [X2]) that 
there exists
$c_1,...,c_t \in V_0$  such that
$$
P(i\otimes \alpha)= \sum_{1\leq k\leq t} c_k, \ 
d(c_k) = \frac{1}{t}d(i) d(\alpha), k=1,...,t. \tag 3.6.0
$$  
Note that if $P^{-1}(P(i\otimes \alpha))= \{i_1\otimes \alpha_1,...
i_s\otimes \alpha_s \}$, then $st = N$. \par
Note we identify the covariant representations of the coset with
the basis of $P(W_0) = V_0$.
The univalence of $A:=P(i\otimes \alpha), i\otimes \alpha \in W_0$ are
given
by:
$\omega_A = \exp(2\pi i (\Delta_i - \Delta_\alpha))$,
 where
$\Delta_i, \Delta_\alpha$ are the conformal dimensions (cf. \S1.5, and
if
$i=(\Lambda',\Lambda''), \Delta_i:= \Delta_{\Lambda'} + \Delta_{\Lambda''}
$).
Note if $A\succ a$, then $ \omega_a = \omega_A$. \par
Under the action of $\sigma$, $W_0$ decomposes into disjoint 
orbits denoted by
$O_p := \{ (i_p\otimes \alpha_p), \sigma((i_p\otimes \alpha_p)),
..., \sigma^{f(O_p)} ((i_p\otimes \alpha_p)).$  
Note that all the elements of $O_p$ are mapped to the same 
element $ P(i_p\otimes \alpha_p)$ by $P$.
Notice that
$$
\langle  P(i_p\otimes \alpha_p),  P(i_p\otimes \alpha_p) \rangle
= \frac{N}{ f(O_p)},
$$ and so
$ P(i_p\otimes \alpha_p)$ decomposes into $g(O_p):=\frac{N}{ f(O_p)}$ 
irreducible
sectors denoted by $x_1(O_p), ..., x_{g(0_p)} (O_p)$ by 
(3.6.0). To summarize,
we have to each orbit $O_p$ associate irreducible
sectors denoted by $x_1(O_p), ..., x_{g(0_p)} (O_p)$, and
$$
\langle x_j(O_p),  x_i(O_{p'}) \rangle = 0
$$ if $O_p\neq O_{p'}$. Moreover the statistical dimension 
$d(x_j(O_p)) = \frac{f(p)}{N} d(i_p)d(\alpha_p)$ by (3.6.0).  Use Th. B
the calculation of $\tau_{G/H}$ is now rather similar to the 
calculation of $\tau_{G_1}$, and in fact essentially the same argument
as (3.4.3) and (3.4.4)  with some change of notations gives the following
$$
\tau_{G/H}(M_L) 
= (\frac{k_G}{k_H})^{3(b_-(L)-b_+(L))} ({D_{G}D_H})^{-n}
\sum_{i_j-\alpha_j \in Q, 1\leq j \leq n} 
L(i_1,...,i_n) \overline{L(\alpha_1,...,\alpha_n)} \tag 3.6.1
$$
where $k_G, k_H$ are given by (1.7.5), $D_G, D_H$ are given by
(1.7.4), $Q$ is the root lattice (note $\lambda\in Q$ iff $
\frac{\tau(\lambda)}{N}\in {\Bbb Z}$), and 
$$
L(i_1,...,i_n)= L(\Lambda(1)',...,\Lambda(n)')  
L(\Lambda(1)'',...,\Lambda(n)'') 
$$ if $i_j = ( \Lambda(j)',\Lambda(j)'), j=1,...,n$. \par
Let us first consider a special case of (3.6.1) when $m''=1$. This
coset is called Coset $W_N$ algebras with critical parameters in
[X1] due to its close relations with $W$ algebras (cf. [FKW] and 
[BBSS]).                 
\proclaim{Proposition 3.6.1}
Suppose $m''=1$. Then there exists a function $F$ which only depends on
the linking matrix of $L$ such that
$$
\tau_G (M_L) \overline{ \tau_H (M_L)} = F \tau_{G/H} (M_L)
.$$
\endproclaim
\demo{Proof}
To save some writing only in this proof we shall denote write 
$i_j=(s_j,t_j)$ and $I_j:=(s_j, t_j; \alpha_j), j=1,2,...,n$. The 
link invariants $L(i_1,...,i_n) \overline{L(\alpha_1,...,\alpha_n)}$
will be denoted by $L(I_1, I_2,..., I_n)$.  Define
$\tau(I_j):= \tau(s_j) +\tau(t_j) - \tau(\alpha_j), 
d(I_j)=d(s_j)d(t_j) d(\alpha_j),$  
$\sigma(I_j):= (\sigma(s_j), \sigma(t_j); \sigma(\alpha_j)),$ and
$\sigma'(I_j):= (s_j, \sigma(t_j); \alpha_j), j=1,2,...,n.$
By symmetry principle Prop. 3.2.2
$$
L(I_1,..., \sigma(I_i),..., I_n) = 
\exp(-2\pi i \sum_{j\neq i } \frac{\tau(I_j)}{N} L_j
\cdot L_i + 2 \pi i \frac{\tau (I_i)}{N} L_i\cdot L_i)  
L(I_1, I_2,..., I_n) \tag 3.6.2
$$

$$
L(I_1,...,\sigma'(I_i),...,I_n)
= \exp(2\pi i (\frac{N-1}{2N} - \frac{\tau(t_i)}{N}) L_i\cdot L_i-
2\pi i \sum_{j\neq i} \frac{\tau(t_j)}{N} L_j\cdot L_i) \tag 3.6.3
$$                                  

\par
Consider the following summation
$$
\sum_{\tau(I_i) \equiv x_i mod (N)} d(I_1)...d(I_n)
L(I_1,...,...,I_n). \tag 3.6.4
$$  Since $\tau(\sigma (I_i)) -\tau(I_i) \equiv 0 \text{\rm mod} (N)$, the set 
which is summed over above is invariant under the action of 
$\sigma$. By using (3.6.2) it is easy to see that if
$$
-\sum_{j\neq i}  x_j L_j\cdot L_i + x_i L_i \cdot L_i  
$$ is not divisible by $N$ for some $i$,   
then the summation in (3.6.4) is $0$.  Let us assume that
$$
\sum_{j\neq i}  
-x_j L_j\cdot L_i + x_i L_i \cdot L_i \equiv 0 mod (N) , \forall i\tag 3.6.5 
$$
By using (3.6.2) repeatedly it is easy to see that under (3.6.5) any term
$L(I_1,...,I_n)$ in the summation (3.6.4) has the property that
$$
L(\sigma^{y_1}(I_1),..., \sigma^{y_n}(I_n)) =L(I_1,...,I_n), \forall y_1,
..., y_n.
$$    
Using this  we have:
$$
\align
\sum_{\tau(I_i) \equiv x_i mod (N)} L(I_1,...,...,I_n)
&= N^n \sum_{ t_i =\sigma^{x_i}(1), \tau(I_i) \equiv x_i \text{\rm mod} (N)}
L(I_1,...,...,I_n) \\
&=  N^n \sum_{ t_i = 1, \tau(I_i) \equiv 0 \text{\rm mod} (N)}
L({\sigma'}^{x_1}(I_1),...,...,{\sigma'}^{x_n}(I_n)) \\
&=  N^n f({x_i}, L)  \sum_{ t_i = 1, \tau(I_i) \equiv 0 \text{\rm mod} (N)}
L(I_1,...,...,I_n),
\endalign
$$ where $f({x_i}, L)$ depends only on the linking matrix of $L$
and numbers $ \{x_i \}$ and in the last step we have repeatedly used
(3.6.3). Also keep in mind that the action of $\sigma$ on the label
$t_i$ is transitive and $ \tau(\sigma^{x_i}(1))=i$.  
 When $\{x_i \}$ does not verify (3.6.5), we
define $f({x_i}, L)=0$. \par
So we have
$$
\align
\tau_G(M_L) \bar \tau_H(M_L) & =
(\frac{k_G}{k_H})^{3(b_-(L) - b_+(L))} D_G^{-n} D_H^{-n}
\sum_{I_i} d(I_1)... d(I_n) L(I_1,...,...,I_n) \\
&= (\frac{k_G}{k_H})^{3(b_-(L) - b_+(L))} D_{G/H}^{-n} N^{-n} \times \\
&\sum_{0\leq x_i\leq N-1} \sum_{\tau(I_i)\equiv x_i mod (N)}
 d(I_1)... d(I_n) L(I_1,...,...,I_n) \\
&= (\frac{k_G}{k_H})^{3(b_-(L) - b_+(L))} D_{G/H}^{-n} N^{-n}
\sum_{0\leq x_i\leq N-1}  N^n f({x_i}, L) \times \\  
&\sum_{t_i = 1,\tau(I_i)\equiv 0 mod (N)}
d(I_1)... d(I_n)L(I_1,...,...,I_n) \\
&=
\sum_{0\leq x_i\leq N-1}f({x_i}, L) \tau_{G/H}(M_L).
\endalign
$$ The proposition now follows by choosing 
$$
F= \sum_{0\leq x_i\leq N-1}f({x_i}, L). 
$$
\enddemo  
\qed
\par
\proclaim{Corollary 3.6.2} 
Under the condition of the previous proposition, there exists
a 3-manifold invariant $F_1$ such that
$$
\tau_G(M_L) \overline{ \tau_H(M_L)} = F_1 (M_L)\tau_{G/H}(M_L)
.$$ 
\endproclaim
\demo{Proof}
By the previous proposition we can define
$$
F_1 (M_L):= \frac{\tau_G(M_L) \overline{ \tau_H(M_L)}}{\tau_{G/H}(M_L)}
,$$ if ${\tau_{G/H}(M_L)}\neq 0$, and 
$$
F_1 (M_L):=0
,$$ if  ${\tau_{G/H}(M_L)}= 0$. $F_1$ is a  3-manifold invariant
since $\tau_G(M_L), 
\bar \tau_H(M_L) $ and \par
${\tau_{G/H}(M_L)}$ are 3-manifold invariants.
\enddemo
\qed
\par
It seems to be possible to choose $F_1$ in Cor. 3.6.2  to be a homotopy
invariant, but the above proof does not show this. \par
Let us consider the simplest case  $H=SU(2)_2\subset G=SU(2)_1\times
SU(2)_1$. 
This coset has central charge $C_{G/H} = 1/2$ and corresponds to
the critical Ising model.
There are three simple objects in $C(G/H)$, denoted by $1:=
(0,0;0), x:=(0,1/2;1/2) , y:=(0,0;1)$, with
univalence  $\omega (x)= \exp(2\pi i \frac{1}{16})$.
As sectors we have $[x^2]= [1]+[y], [xy]=[x], [y^2]=1$. Note that
$C(G/H)$ has the same number of simple objects with identical fusion
rules as $C(SU(2)_2)$, but with one different univalence
$\omega (1/2) =  \exp(2\pi i\frac{3}{16})$.  
We can evaluate $\tau_{G/H}$ by using
the formula (3.6.1)  but instead let us take the following short cut
by using
the similarity with $C(SU(2)_2)$. 
Consider the sequence of algebras
$$          
Hom(x^n, x^n) \subset Hom(x^{n+1}, x^{n+1}) ...
.$$ By a similar argument as in the proof of lemma 3.1.1 (2)  we have
If 
$$h_n := \exp(\frac{\pi i}{8} x^{n-1}(c_{x,x})),
$$ then
$$
h_n^2 = (q-1) h_n + q
,$$ with $q=\exp(\frac{2\pi i}{4})$. 
The Bratteli diagram of
the algebras
$$
Hom(x^n, x^n) \subset Hom(x^{n+1}, x^{n+1}) ...
$$ is determined by the fusion rules and it is then easy to see
that this sequence of algebras are a special case of the sequence
of algebras analyzed by Jones in \S5.2 of [J1] (with $\tau= 1/2$).
Now an almost identical argument as in the proof of  Cor. 4.11 of [KM]
shows that
$$
L(x,x,...,x) =
\sqrt{2} \exp(\frac{2\pi i}{16} L\cdot L ) \tilde V_L
,$$ where $\tilde V_L$ is the modified Jones polynomial
as in Cor. 4.11 of [KM]. Now the proof of Th. 7.1 of [KM] 
applies verbatim, with the constant $c= \exp(\frac{-6\pi i}{16})$ on
page 524 of [KM] replaced by $ \exp(\frac{-2\pi i}{16})$. The result
is a formula for $\tau_{G/H}(M)$ as in Th. 7.1 of [KM],
with their $c= \exp(\frac{-6\pi i}{16})$ replaced by
$ \exp(\frac{-2\pi i}{16})$. Let us record this result in the
following proposition:
\proclaim{Proposition 3.6.3}
Let $M$ be a closed, oriented 3-manifold, $H=SU(2)_4\subset G=SU(2)_2
\times SU(2)_2$. Then
$$
\tau_{G/H}(M) = \sum_\theta {c'}^{\mu (M_\theta)}
,$$ where $c'= \exp(\frac{-2\pi i}{16})$, and $\mu (M_\theta)$
is the $\mu$-invariant of the spin structure $\theta$ on $M$
and the sum is taken over all spin structures.
\endproclaim   
It is also interesting to compare the above coset with
the simplest parafermion coset $H_1:=U(1)_4\subset G_1:=SU(2)_2$. This
coset has central charge $1/2$ with three simple objects. 
In the notations of \S3.5 these three simple objects are given
by 
$x_1:= (1/2, 1), y_1:= (1,0), 1:=(0,0)$ with
fusion rules
$$
[x_1^2]= [1] +[y_1], [x_1y_1] = [x_1], [y_1^2]=[1]
,$$ and the conformal dimension (modulo integers) of $x_1$
is $\frac{1}{16}$.  Now exactly the same argument as above
shows that the coset $H_1:=U(1)_4\subset G_1:=SU(2)_2$ gives the
same 3-manifold invariant as in Prop. 3.6.3. \par
In general there is no clear relation between 
$ \tau_G(M_L)\bar \tau_H(M_L)$ and ${\tau_{G/H}(M_L)}$
as in Prop. 3.6.1 by inspecting (3.6.1) and using the symmetry
principle. This is especially true in the case when the action of
$\sigma$ has fixed points on $exp$. 
For example the coset $SU(2)_{2(k+l)}\subset SU(2)_{2k} \times 
SU(2)_{2l}$ has a fixed point $(k/2,l/2; \frac{k+l}{2})$. Any of the
diagram automorphisms of the three groups $SU(2)_{2(k+l)}, 
SU(2)_{2k}$ and $SU(2)_{2l}$ preserves the parity of the representation
labels.
One runs into a problem
similar to the example at the end of \S3.5 if 
$ \tau_G(M_L)\bar \tau_H(M_L)$ can be expressed in terms of 
${\tau_{G/H}(M_L)}$.  
\par 
The simplest case when there is a fixed point under the action of
$\sigma$ is the coset $H=SU(2)_4 \subset G=SU(2)_2\times SU(2)_2$. 
In this case there is a unique fixed point $(1/2,1/2;1)$. 
The corresponding unitary modular category $C(G/H)$ has 
$13$ simple objects: $11$ from the $11$ orbits of 
${\Bbb Z}_2$ action on the set $(i_1, i_2,\alpha)$ with $
i_1 + i_2 -\alpha \in {\Bbb Z}$, where each of this orbit contains
two elements, and $2$ from the fixed point  $(1/2,1/2;1)$ by (3.6.0). 
Note from (3.6.1) that 
in this case the three manifold invariants 
are expressed in terms of Jones polynomial
at sixth and fourth roots of unity via cabling (cf. [KM]), 
and since Jones polynomial
at sixth and fourth roots of unity 
can be expressed in terms of classical topological
invariants, it is an interesting question to see if one can do the
same with our $\tau_{G/H}$ along the lines of  [KMZ]. \par
\subheading{\S3.7 A ``Maverick" Coset}
Let us consider a ``Maverick" Coset  $H=SU(2)_8 \subset G=SU(3)_2$  
(cf. [DJ1]) considered at the 
end of [X1].
This coset is also considered in   [X2] as the first counter-example
to a hypothesis of Kac-Wakimoto (cf. [KW], [X8]) 
which is not a conformal inclusion.
This coset verifies Th. A, so there is 
a unitary modular category $C(G/H)$ as constructed in \S1.7.  However
as shown in [X1] 
this example does not verify the condition of Th. B, so we need to
find a way to calculate $\tau_{G/H}$.\par
As shown in [X1], $C(G/H)$  has 6 simple objects, denoted by the
following :
$1=(00,0), x=(00,4), y=(10,2), \bar y= (01,2), 
z=(10,4), \bar z = (01,4)
$. The notation for representation is slightly different from \S1.5. Here the 
representation of $\hat su(3)_2$ is denoted by two integers 
$(\lambda_1\lambda_2), \lambda_1+ \lambda_2 \leq 2 
$ such that the vacuum representation is given by 
$(00)$ and the representation of $\hat su(2)_8$ is labelled by an integer
$k, 0\leq k \leq 8$ 
such that $0$ denotes  the vacuum representation. \par 

The nontrivial fusion rules are
$$
[x^2]=[1]+[x], [y\bar y]=[1]+[x], [z^3]=[1], [y]= [xz]
.$$ The central charge of this coset is $0.8$, which is 
the same as the diagonal coset $H_1=SU(3)_2 \subset G_1=SU(3)_1 \times
SU(3)_1$. Furthermore the diagonal coset $H_1=SU(3)_2 \subset G_1=SU(3)_1 \times
SU(3)_1$ also has $6$ irreducible sectors labeled in the following
$$
 1=(00,00; 00), x_1 =(00,00; 11), y_1=(00, 01; 01), \bar y= (00,10;10), 
$$
$$
z_1=(00, 20; 20), \bar z_1 = (00,02;20).
$$  where we use the same notations 
for representations of $SU(3)$ as above. 
It is an easy exercise to show that the map
$1\rightarrow 1, [x]\rightarrow [x_1], [y]\rightarrow [y_1],
 [z]\rightarrow [z_1]$ is a ring 
isomorphism, and 
$$
\omega(y) = \omega (y_1) =\exp (2\pi i \frac{1}{15}).
$$ Now notice that $y$ (resp. $y_1$) is a generating element in 
 $C(G/H)$ (resp. $C(G_1/H_1) $) , 
and by a  simpler argument as in lemma 3.1.1  $Hom(y^n,y^n)$ 
(resp. $Hom(y_1^n,y_1^n) $) 
is 
isomorphic to  Wenzl's representation (cf. [W1])
$\pi^{(3,5)}(H_n(q))$ with $q=\exp (\frac{2\pi i}{5})$. 
Hence $C(G/H)$ and $C(G_1/H_1) $ are compatible as defined before 
lemma 1.7.5.
By Lemma 1.7.5 we have shown 
$$
\tau_{G/H}(M) =  \tau_{G_1/H_1}(M).
$$ 
Finally  note that the parafermion coset (cf. \S3.5) $U(1)_6\subset SU(2)_3$
also has central charge $0.8$ and $6$ simple objects. A simple 
calculation shows that the element $(1/2,1)$ is a generating element
with univalence $\exp(2\pi i \frac{1}{15})$. Now similar argument
as above shows that the modular category associated with 
the coset $U(1)_6\subset SU(2)_3$ is compatible with 
$C(G/H)$
or $C(G_1/H_1) $ above, and so by  Lemma 1.7.5 give the same
3-manifold invariants. We can therefore use (3.5.3) to calculate 
the three manifold invariant associated with 
the coset  $H=SU(2)_8 \subset G=SU(3)_2$. \par
For more such ``Maverick" Cosets, see [DJ2] and [FSS2]. It will 
be interesting if one can calculate the corresponding 
3-manifold invariants by 
a similar identification as above. \par 
\subheading{3.8 A Question}
An interesting question is to investigate the perturbative aspects
of $\tau_{G/H}$ similar to the case of $\tau_{G}$ (cf.  
[Ro], [Gar], [O]). Note that 3-d Chern-Simons action
and 2-d gauged WZW action exist for the coset (cf. [MS], [KS]), 
it remains to see if one can draw any conclusion on the calculations
of $\tau_{G/H}$ similar to that of [Ro].\par
\heading References \endheading
\roster
\item"{[A]}" H. H. Andersen, {\it Tensor product of quantized tilting
modules,} Comm. Math. Phys. 149 (1992), 149-159.
\item"{[ABI]}" D. Altschuler, M. Bauer and C. Itzykson, {\it The
branching rules of conformal embeddings},  Comm.Math.Phys.,
 {\bf 132}, 349-364
(1990).  
\item"{[BAF]}" D.Buchholz, C.D'Antoni and K.Fredenhagen,
{\it local factorizations.  The universal structure of
local algebras}, Comm.Math.Phys.,  111, 123-135 (1987).
\item"{[BBSS]}" F.A. Bais, P. Bouwknegt, K. Schoutens, M. Surridge,
{\it Coset construction for extended Virasoro algebras,} Nucl. Phys. B
304, 371-391, (1988).
\item"[BE1]" J. B\"{o}ckenhauer, D. E. Evans,
{\it Modular invariants, graphs and $\alpha$-induction for
nets of subfactors. I.},
Comm.Math.Phys., {\bf 197}, 361-386, 1998.
\item"[BE2]" J. B\"{o}ckenhauer, D. E. Evans,
{\it Modular invariants, graphs and $\alpha$-induction for
nets of subfactors. II.},
Comm.Math.Phys., {\bf 200}, 57-103, 1999.
\item"[BE3]" J. B\"{o}ckenhauer, D. E. Evans,
{\it Modular invariants, graphs and $\alpha$-induction for
nets of subfactors. III.},
Comm.Math.Phys., {\bf 205}, 183-228, 1999.
\item"[BE3]" J. B\"{o}ckenhauer, D. E. Evans,
{\it Modular invariants, graphs and $\alpha$-induction for
nets of subfactors. III.},
Comm.Math.Phys., {\bf 205}, 183-228, 1999. Also see hep-th/9812110.
\item"[BEK1]" J. B\"{o}ckenhauer, D. E. Evans, Y. Kawahigashi,
{\it On $\alpha$-induction, chiral generators and modular invariants
for subfactors}, Comm.Math.Phys., {\bf 208}, 429-487, 1999. Also
see math.OA/9904109.
\item"[BEK2]" J. B\"{o}ckenhauer, D. E. Evans, Y. Kawahigashi,
{\it Chiral structure of modular invariants for subfactors},
Comm.Math.Phys., {\bf 210}, 733-784, 2000.    
\item"[C]" A. Connes, {\it Noncommutative Geometry}, Academic Press, Inc.
1994.
\item"[CM]" T. D. Cochran, P. Melvin, {\it Quantum Cyclotomic orders
of 3-manifolds}, to appear in Topology.
\item"[DHR]" S. Doplicher, R. Haag and J.E. Roberts, {\it Local
observables and particle statistics I,} Comm. Math. Phys. 23, 
199-230 (1971).
\item"{[DL]}" C. Dong and J. Lepowsky, {\it Generalized vertex algebras
and relative vertex operators,} Progress in Mathematics, 112 (1993)
\item"[DJ1]" D. Dunbar and K. Joshi,
{\it Characters for Coset conformal field theories and Maverick
examples}, Inter. J. Mod. Phys. A, Vol.8, No. 23 (1993), 4103-4121.
\item"[DJ2]"  D. Dunbar and K. Joshi, {\it Maverick examples of 
Coset conformal field theories}, 
Mod. Phys. Letters A, Vol.8, No. 29 (1993), 2803-2814.
\item"[EK1]" D. E. Evans and Y. Kawahigashi, {\it On Ocneanu's theory
of asymptotic inclusions for subfactors, topological quantum field theories
and quantum doubles,} Int. J. Math. vol. 6 , no. 2 (1995) 205-228. 
\item"[EK2]" D. E. Evans and Y. Kawahigashi, {\it Quantum Symmetries
on Operator Algebras}, Oxford University press, 1998.
\item"[EK3]" D. E. Evans and Y. Kawahigashi, {\it Orbifold subfactors
from Hecke algebras II. Quantum doubles and braiding,} Comm. Math. Phys.
196 (1998) 331-361. 
\item"{[FG]}"  J. Fr\"{o}hlich and F. Gabbiani, {\it Operator algebras and
CFT}, Comm. Math. Phys., {\bf 155}, 569-640 (1993).  
\item"{[FKW]}" E. Frenkel, V. Kac and M. Wakimoto, {\it Characters and
Fusion rules for W-algebras via Quantized Drinfeld-Sokolov Reductions},
RIMS-861, 1992.
\item"{[Fre]}"  K.Fredenhagen, {\it Generalizations of the theory
of superselection sectors}, in "The algebraic theory of
superselection sectors", D.Kastler ed., World Scientific, Singapore 1990. 
\item"{[FRS]}" K.Fredenhagen, K.-H.Rehren and B.Schroer
,\par
{\it Superselection sectors with braid group statistics and
exchange algebras. II}, Rev. Math. Phys. Special issue (1992), 113-157.  
\item"{[FSS1]}" J. Fuchs, B. Schellekens and C. Schweigert,
{\it The resolution of field identification fixed points in diagonal coset
theories,} Nucl. Phys. B 461 (1996) 371, hep-th/9509105.
\item"{[FSS2]}" J. Fuchs, B. Schellekens and C. Schweigert,
{\it From Dynkin diagram symmetries to fixed point structures,}
Commun. Math. Phys. 180 (1996) 39 
\item"{[Gar]}" S. Garoufalidis, {\it On finite type 3-manifold invariants I},
J. Knot. Th. Ram. 5 (1996), no.4, 441-461.
\item"{[GKO]}" P. Goddard, A. Kent  and D. Olive,  { \it
Unitary representations of
Virasoro and super-Virasoro algebras}
Comm. Math. Phys. 103 (1986), No. 1, 105-119.  
\item"{[GL]}"  D.Guido and R.Longo, {\it  The conformal  spin and
statistics theorem},  \par
Comm.Math.Phys., {\bf 181}, 11-35 (1996) 
\item"{[J1]}" V. F. R.  Jones, 
{\it Index for subfactors}, Invent. Math. 72 (1983), no.1, 1-25.
\item"{[J2]}" V. F. R.  Jones, 
{\it Hecke algebra representations of braid groups and link polynomials,}
Ann. Math. 126 (1987), 335-388.
\item"{[KLM]}" Y. Kawahigashi, R. Longo and M. M\"{u}ger,
{\it Multi-interval Subfactors and Modularity of Representations in
Conformal Field theory}, Preprint 1999, see also math.OA/9903104.
\item"{[KM]}" R. C. Kirby and P. Melvin, 
{\it The 3-manifold invariants of 
Witten and \par
Reshetikhin-Turaev for sl(2,C),}
Invent. Math. 105 (1991), 473-545. 
\item"{[KMZ]}" R. C. Kirby, P. Melvin and Zhang, Xingru,
{\it Quantum invariants at sixth root of unity}, Comm. Math. Phys. 151 (1993),
no.3, 607-617.
\item"{[K]}" H. Kosaki, {\it Extension of Jones' theory on index to
arbitary factors,} J. Funct. Anal. 66, 123-140 (1986).
\item"{[KS]}" D. Karabali and H. J. Schnitzer, {\it BRST quantization
of the gauged WZW action and coset conformal field theories,}
Nucl. Phys. B 329 (1990), no. 3, 649-666.
\item"{[KT1]}" T. Kohno and T. Takata, {\it Symmetry for Witten's
3-manifold invariants for $Sl(n, {\Bbb C})$,}
J. Knot Th. Ramif. 2 (1993), 149-169. 
\item"{[KT2]}" T. Kohno and T. Takata, {\it Level-rank duality for Witten's
3-manifold invariants,}
Adv. Stud. Pure Math. 24 (1996), 243-264.
\item"{[KW]}"  V. G. Kac and M. Wakimoto, {\it Modular and conformal
invariance constraints in representation theory of affine algebras},  
Advances in Math., {\bf 70}, 156-234 (1988).
\item"{[Lang]}" S. Lang, {\it Algebraic Number Theory}, GTM 110, 
Springer-Verlag, 1994.
\item"{[L1]}"  R. Longo, {\it Index of subfactors and statistics of
quantum fields}, I, Comm. Math. Phys., {\bf 126}, 217-247 (1989.
\item"{[L2]}"  R. Longo, {\it Index of subfactors and statistics of
quantum fields}, II, Comm. Math. Phys., {\bf 130}, 285-309 (1990).
\item"{[L3]}"  R. Longo, {\it Minimal index and braided subfactors,
} J.Funct.Analysis {\bf 109} (1992), 98-112.
\item"{[LR]}"  R. Longo and K.-H. Rehren, {\it Nets of subfactors},
Rev. Math. Phys., {\bf 7}, 567-597 (1995).  
\item"{[LVW]}" W. Lerche, C. Vafa and N. P. Warner, 
{\it Chiral rings in N=2 superconformal theories,}
Nucl. Phys. B324
(1989) 427.  
\item"{[MOO]}" H. Murakaimi, T. Ohtsuki and M. Okada, {\it
Invariants of 3-manifolds derived from linking matrices of framed links,}
Osaka. J. Math., 29 (1992), 545-572. 
\item"{[MS]}" G. Moore and N. Serberg, {\it Taming the conformal zoo},
Lett. Phys. B  {\bf 220}, 422-430, (1989).
\item"{[MW]}" G. Masbaum and H. Wenzl, {\it Integral modular categories
and integrality of quantum invariants at roots of unity of prime order,}
J. reine angew. Math. 505 (1998), 209-235. 
\item"{[O]}" T. Ohtsuki, T., {\it A filtration of the set of integral
homology 3 spheres}, Doc. Math. 1998, Extra Volume II, 473-482.
\item"{[PP]}" M.Pimsner and S.Popa,
{\it Entropy and index for subfactors}, \par
Ann. \'{E}c.Norm.Sup. {\bf 19},
57-106 (1986). 
\item"[PS]" A. Pressly and G. Segal, {\it Loop Groups,} O.U.P. 1986.
\item"[Reh]" Karl-Henning Rehren, {\it Braid group statistics and their
superselection rules} In : The algebraic theory of superselection
sectors. World Scientific 1990 
\item"[Ro]" L. Rozansky, {\it Trivial connection contribution to Witten's
 invariants and finite type invariants of rational homology 
spheres}, Comm. Math. Phys. 183 (1997), no. 1, 23-54.
\item"[RT]" N. Y. Reshetikhin and V. G. Turaev, {\it Invariants of
3-manifolds via link polynomials and quantum groups,}
Invent. Math. 103 (1991), 547-597.
\item"[SY]" A. N. Schellekens and S. Yankielowicz, Nucl. Phts. B 324,
67, (1990).
\item"[Tu]" V. G. Turaev, {\it Quantum invariants of knots and
3-manifolds,} Walter de Gruyter, Berlin, New York 1994. 
\item"[TW1]" V. G. Turaev and  H. Wenzl, {\it Quantum 3-manifold
invariants associated with classical simple Lie algebras,}
Int. J. Math. 4 (1993), 323-358.
\item"[TW2]" V. G. Turaev and  H. Wenzl, {\it Semisimple and modular
categories from link invariants,}
Math. Ann. 309 (1997), 411-461.
\item"{[Wa]}"  A. Wassermann, {\it Operator algebras and Conformal
field theories III},  Invent. Math. 133 (1998), 467-538.
\item"[W1]" H.Wenzl, {\it Hecke algebras of type $A_n$ and subfactors,}
Invent. Math. 92 (1988), 349-383.
\item"[W2]" H.Wenzl, {\it Braids and invariants of 3-manifolds,}
Invent. Math. 114 (1993), 235-275.
\item"[W3]" H.Wenzl, {\it $C^*$ tensor categories from quantum groups,}
Journal AMS 11 (1998) no.2 261-282.
\item"[Wi]" E. Witten, {\it Quantum field theory and the Jones polynomial,}
Comm. Math. Phys. 121 (1989), 351-399.
\item"[X1]" F.Xu, {\it Algebraic coset conformal field theories},
 36 pages, q-alg/9810053.
\item"[X2]" F.Xu, {\it Algebraic coset conformal field theories II}, \par
 25 pages,  math.OA/9903096.
\item"[X3]" F.Xu, {\it  On a conjecture of Kac-Wakimoto, }
18 pages, Math.RT/9904098. 
\item"[X4]" F.Xu, {\it   New braided endomorphisms from conformal
inclusions, } \par
Comm.Math.Phys. 192 (1998) 349-403.
\item"[X5]" F.Xu, {\it Standard $\lambda$-lattices from quantum groups,}
  Invent.Math. Vol.134, 455-487 (1998). 
\item"[X6]" F.Xu, {\it Orbifold constructions in subfactors, } \par
Comm.Math.Phys. 166 (1994),237-253 
\item"[X7]" F.Xu, {\it Jones-Wassermann subfactors for 
Disconnected Intervals}, 40 pages,
 q-alg/9704003.
\item"[X8]" F.Xu, {\it
Applications of Braided endomorphisms from Conformal inclusions,}
 Inter. Math. Res. Notice., No.1, 5-23 (1998),  also see q-alg/9708013,
and Erratum, Inter. Math. Res. Notice., No.8, (1998).
\item"{[Y]}"  S. Yamagami, {\it A note on Ocneanu's approach to Jones
index theory}, Internat. J. Math., {\bf 4}, 859-871 (1993). 
\endroster 
\enddocument